\documentclass[12pt,a4paper]{article}
\usepackage[T2A]{fontenc}
\usepackage[cp1251]{inputenc}
\usepackage{amsmath,amssymb,amsthm,amsfonts,amscd}

\begin{document}
\sloppy
\date{}

\title{Geometric approach to stable
homotopy groups of spheres II. The
Kervaire invariant}

\author{Petr M. Akhmet'ev \thanks{This paper was presented at the Conference 
"The Kervaire Invariant and Stable Homotopy Theory"' Edinburgh, 25-29 April 2011 }}

\sloppy \theoremstyle{plain}
\newtheorem{theorem}{Theorem}
\newtheorem*{main*}{Main Theorem}
\newtheorem*{theorem*}{Theorem}
\newtheorem{lemma}[theorem]{Lemma}
\newtheorem{proposition}[theorem]{Proposition}
\newtheorem{corollary}[theorem]{Corollary}
\newtheorem{conjecture}[theorem]{Conjecture}
\newtheorem{problem}[theorem]{Problem}

\theoremstyle{definition}
\newtheorem{definition}[theorem]{Definition}
\newtheorem{remark}[theorem]{Remark}
\newtheorem*{remark*}{Remark}
\newtheorem*{example*}{Example}
\newtheorem{example}[theorem]{Example}
\def\aa{\dot{a}}
\def\i{{\bf i}}
\def\j{{\bf j}}
\def\k{{\bf k}}
\def\hh{{\bf \dot{h}}}
\def\e{{\bf e}}
\def\f{{\bf f}}
\def\h{{\bf h}}
\def\dd{{\dot{d}}}
\def\bb{{\dot{b}}}
\def\Z{{\Bbb Z}}
\def\R{{\Bbb R}}
\def\RP{{\Bbb R}\!{\rm P}}
\def\N{{\bf N}}
\def\C{{\Bbb C}}
\def\A{{\bf A}}
\def\D{{\bf D}}
\def\Q{{\mathbf Q}}
\def\QQ{{\dot{\mathbf Q}}}
\def\E{{\bf H}}
\def\F{{\bf F}}
\def\J{{\bf J}}
\def\JJ{{\dot{\bf J}}}
\def\G{{\bf G}}
\def\I{{\bf I}}
\def\II{{\dot{\bf I}}}
\def\H{{\bf E}}
\def\fr{{\operatorname{fr}}}
\def\st{{\operatorname{st}}}
\def\mod{{\operatorname{mod}\,}}
\def\cyl{{\operatorname{cyl}}}
\def\dist{{\operatorname{dist}}}
\def\sf{{\operatorname{sf}}}
\def\dim{{\operatorname{dim}}}
\def\dist{\operatorname{dist}}

\maketitle

\begin{abstract}

A solution to the Kervaire invariant problem is presented. We
introduce the concepts of abelian structure on skew-framed
immersions, bicyclic structure on $\Z/2^{[3]}$--framed immersions, and
quaternionic-cyclic structure on $\Z/2^{[4]}$--framed immersions.
Using these concepts, we prove that for sufficiently large $n$,
$n=2^{\ell}-2$, an arbitrary skew-framed immersion in Euclidean
$n$-space $\R^n$ has zero Kervaire invariant. Additionally, for
$\ell \ge 12$ (i.e., for $n \ge 4094$) an arbitrary skew-framed
immersion in Euclidean $n$-space $\R^n$ has zero Kervaire
invariant if this skew-framed immersion admits a compression of
order $16$.
\end{abstract}

\section{Self-intersections of immersions and the
Kervaire invariant}

The Kervaire invariant one problem was for many years an open problem in
algebraic topology. For algebraic approaches to the problem see
Snaith [S], Barratt-Jones-Mahowald [B-J-M] and
Cohen-Jones-Mahowald [C-J-M]. Recently, Hill, Hopkins, and Ravenel
obtained a solution of this problem for all dimensions, except
$n=126$ (see [H-H-R]). We consider an alternative geometric approach.
For a different geometric approach see Carter [C1], [C2].

The proof is based on a paper by P.J. Eccles [E1], identifying the
Kervaire invariant with the number of multiple points of an
immersion using the Kahn-Priddy map $MO(1) = P^{\infty} \to S^0$.
Since all Steenrod squares are non-zero in the mapping cone of
this map, an element $\alpha\in \pi^S_{2k-2}MO(1)$ is detected by
$Sq^k$ on $a_{k-1}$ if and only if its image is detected by the
secondary operation coming from $Sq^kSq^k$ ($k$ a power of $2$ and
$n+2=2k$). By W. Browder's result [B] this gives the geometric
interpretation of the Kervaire invariant. Namely, Peter Eccles
showed that the Kervaire invariant can be interpreted as the parity
of the number of $(2k-2)$-fold points of a corresponding
immersion. But with our approach we do not directly make use of
this.

Using the equivalence of the James-Hopf invariant and the
Steenrod-Hopf invariant this is also equivalent to the stable Hopf
invariant $j_2(\alpha) \in \pi^S_{2k-2}D_2MO(1)$ having Hurewicz
image $a_{k-1}^2$.  This element represents the double point
manifold of an immersion corresponding to $\alpha$.  Forgetting
the additional structure on the double point manifold maps
$j_2(\alpha)$ to an element $\beta \in \pi^S_{2k-2}MO(2)$ with
Hurewicz image $a_{k-1}^2$.

If $\alpha$ corresponds to an immersion $f$  then $\beta$ corresponds to an immersion $g$, as in the following paragraph.  By
dualizing to cohomology and using the Thom isomorphism one can check that
$h(\beta) = a_{k-1}^2$ if and only if the characteristic number
$\langle w_2^{k-1},[N] \rangle = 1$ which is our interpretation of the Kervaire invariant.
More detailed information can be found in [A-E]. As far as I know, there is no explicitly geometric
proof of  Eccles' Theorem on the Kervaire invariant.

Consider a smooth immersion $f: M^{n-1} \looparrowright \R^n$, $n=
2^{\ell} -2$, $\ell>1$  in general position and having codimension
1. We denote by $g: N^{n-2} \looparrowright \R^n$ the immersion of
the manifold of self-intersections.

\begin{definition}
The {\it Kervaire invariant} of the immersion $f$ is defined by the
formula

\begin{eqnarray}\label{arf}
 \Theta^{sf}(f) = \langle \eta_N^{\frac{n-2}{2}}; [N^{n-2}] \rangle,
\end{eqnarray}
 where $\eta_N = w_2(N^{n-2})$ denotes the second normal Stiefel-Whitney class of
the manifold  $N^{n-2}$.
\end{definition}

The Kervaire invariant is an invariant of the regular cobordism class of the
immersion $f$. Moreover, the Kervaire invariant determines a homomorphism
\begin{eqnarray}\label{1}
\Theta^{sf}: Imm^{sf}(n-1,1) \to  \Z/2,
\end{eqnarray}
the cobordism group (and cobordism groups mentioned below) of immersions is defined in [A1].
The normal bundle $\nu_g$ of the immersion $g:
N^{n-2} \looparrowright \R^n$ is a 2-dimensional
bundle over $N^{n-2}$, which is naturally equipped with a $\D$-framing, where $\D$
denotes the dihedral group of order $8$.
The classifying map of this bundle
(and its corresponding characteristic class) are denoted by  $\eta_N: N^{n-2} \to K(\D,1)$.
The pair $(g,\eta_N)$ represents an element of the cobordism group $Imm^{\D}(n-2,2)$.
The passage from $f$ to $(g,\eta_N)$ gives rise to a well defined
homomorphism
\begin{eqnarray}\label{2}
\delta^{\D}: Imm^{sf}(n-1,1) \to Imm^{\D}(n-2,2).
\end{eqnarray}

The cobordism group $Imm^{sf}(n-k,k)$ generalizes the cobordism group $Imm^{sf}(n-1,1)$.
The new group is defined as the cobordism group of triples $(f,\Xi,\kappa_M)$,
where $f: M^{n-k} \looparrowright \R^n$ is an immersion of a compact closed manifold; moreover there is given
a morphism of bundles (a bundle map) $\Xi:
\nu_f \cong k \kappa_M$, which is invertible,
i.e., which is a fiberwise isomorphism, and which is called a skew-framing,
where $\nu_f$ denotes the normal bundle of the immersion $f$ and $\kappa_M$ is a given
line bundle over $M^{n-k}$, whose characteristic class is also denoted by $\kappa_M \in
H^1(M^{n-k};\Z/2)$.  The relation of cobordism on the set of triples is the standard
one (see section 1 of [A1] for more details for the definition of the cobordism relation).

The group $Imm^{\D}(n-2,2)$ is generalized in the following way. We
shall define cobordism groups $Imm^{\D}(n-2k,2k)$. Each element of the group
$Imm^{\D}(n-2k,2k)$ is represented by a triple $(g,\Psi,\eta_N)$, where $g: N^{n-2k} \looparrowright \R^n$
is an immersion, $\Psi$  is a dihedral framing of codimension $2k$, i.e., a fixed isomorphism
$\Xi: \nu_g \cong k \eta_N$, and where $\eta_N$ is a $2$-dimensional bundle over $N^{n-2k}$
with structure group $\D$, $\tau$ is the universal 2-bundle over $K(\D,1)$.
The characteristic mapping of this bundle, and also
the corresponding characteristic Euler class (respectively, the universal characteristic
Euler class) will be denoted also by $\eta_N: N^{n-2k} \to K(\D,1)$, $\eta_N \in H^2(N^{n-2k};\Z/2)$
(respectively, $\tau \in H^2(K(\D,1);\Z/2)$).

The mapping $\eta_N$ is also called characteristic
for the bundle $\nu_g$, since $\nu_g
\cong k \eta_N$.

We define the Kervaire homomorphism ($\ref{1}$) (see $(\ref{arf})$) as the composition of the
homomorphism ($\ref{2}$) and a homomorphism
\begin{eqnarray}\label{3}
\Theta^{\D} : Imm^{\D}(n-2,2) \to \Z/2, \quad
\Theta^{\D}(g,\Psi,\eta_N)=\langle \eta_N^{\frac{n-2}{2}}; [N^{n-2}] \rangle.
\end{eqnarray}
The homomorphism $(\ref{3})$ is called the Kervaire invariant of a $\D$--framed immersion.

The Kervaire homomorphism can be defined in more general situations
by means of a direct generalization of the homomorphisms
($\ref{1}$) and ($\ref{3}$):
\begin{eqnarray}\label{4a}
\Theta^{sf}_{k}: Imm^{sf}(n-k,k) \to \Z/2, \qquad \Theta^{sf}_k: =
\Theta^{\D}_k  \circ \delta^{\D}_k.
\end{eqnarray}
\begin{eqnarray}\label{44}
\Theta^{\D}_{k} : Imm^{\D}(n-2k,2k) \to \Z/2, \qquad
\Theta^{\D}_{k}[(g,\Psi,\eta_N)] = \langle
\eta_N^{\frac{n-2k}{2}}; [N^{n-2k}] \rangle.
\end{eqnarray}
For $k = 1$ the new homomorphism $(\ref{4a})$  coincides with the homomorphism
$(\ref{3})$ already defined; moreover the following diagram, in which the homomorphisms
$J^{sf}$ and $J^{\D}$ were defined in the first part  of the paper [A1] (Proposition 2), is
commutative:
\begin{eqnarray}\label{5}
\begin{array}{ccccc}
Imm^{sf}(n-1,1) & \stackrel  {\delta^{\D}}{\longrightarrow} &
Imm^{\D}(n-2,2) & \stackrel{\Theta^{\D}}{\longrightarrow} & \Z/2  \\
\downarrow J^{sf}_k & &  \downarrow J^{\D}_{k}  &  &  \vert \vert \\
Imm^{sf}(n-k,k) & \stackrel{\delta^{\D}_{k}}{\longrightarrow} &
Imm^{\D}(n-2k,2k) &
\stackrel{\Theta^{\D}_k}{\longrightarrow} & \Z/2.  \\
\end{array}
\end{eqnarray}

We shall need to generalize formula ($\ref{44}$) for immersions
with framing of a more general form. Denote by $\Z/2^{[s]}$ the
wreath product of $2^{s-1}$ copies of the cyclic group $Z/2$. This
group is a subgroup of the orthogonal group $O(2^{s-1})$, and can
be defined in the following way:

Transformations in $Z/2[s]$ leave invariant the collection of $(s-1)$ sets
$\Upsilon_{s}$,
$\Upsilon_{s-1}$, $\dots$, $\Upsilon_2$ of coordinate subspaces. The set of subspaces $\Upsilon_i$,
$2 \le i \le s$
consists of the $2^{i-1}$ coordinate subspaces $(Lin(\e_1, \dots, \e_{2^{s-i}}), \dots, Lin(\e_{2^{s-1}-2^{s-i}+1},
\dots, \e_{2^{s-1}}))$, spanned by the orthonormal basis vectors.
The blocks of basis vectors are disjoint and all of the same size.

In particular, in this new notation the dihedral group $\D$ will be denoted
by $\Z/2^{[2]}$. This group is defined as the subgroup of orthogonal transformations
of the plane, carrying the set $\Upsilon_2 =
\{Lin(e_1),Lin(e_2)\}$ of lines into itself. In this paper
we shall make use of the groups $\Z/2^{[s]}$ for $2 \le s \le 5$. By definition, there
is an inclusion $\Z/2^{[s]} \subset \Z/2 \wr
\Sigma(2^{s-1})$, which coincides with the inclusion of
a 2-Sylow subgroup of the symmetric group $\Sigma(2^{s})$. E.g., the dihedral group
$\Z/2^{[2]}$ is a 2-Sylow subgroup of $\Sigma(4)$.
\[  \]

Consider an immersion $g: N^{n-k2^{s-1}} \looparrowright \R^n$ in general position and of
codimension $k2^{s-1}$. We say that the immersion $g$ is $\Z/2^{[s]}$-framed (with multiplicity
$k$), if an isomorphism $\Psi: \nu_g \cong k \eta_N$ is given between the normal
bundle  $\nu_g$ of the immersion $g$ and the Whitney sum of $k$ copies of a $2^{s-1}$-
dimensional bundle $\eta_N$ with structure group $\Z/2^{[s]}$.

The bundle  $\eta_N$ is classified by a mapping $\eta_N:
N^{n-k2^{s-1}} \to K(\Z/2^{[s]},1)$.
(The corresponding characteristic class is also denoted by $\eta_N$.) The characteristic
class of the universal $2^{s-1}$-dimensional $\Z/2^{[s]}$-bundle over $K(\Z/2^{[s]},1)$
is denoted by $\tau_{[s]}$. Therefore, $\eta_N^{\ast}(\tau_{[s]})=\eta_N$. The mapping $\eta_N$ is also called a
characteristic map for the bundle $\nu_g$, since $\nu_g \cong k \eta_N$.

The set of all possible triples $(g,\Psi,\eta_N)$, as described above, generate the
cobordism group $Imm^{\Z/2^{[s]}}(n-k2^{s-1},k2^{s-1})$. In some considerations we use
an additional index in the notation, connected with the structure group. For
example, a representative of the group $Imm^{\Z/2^{[2]}}(n-2k,2k)$ will sometimes
be denoted by  $(g_{[2]},\Psi_{[2]}, \eta_{N_{[2]}})$ and so on.

The manifold of self-intersections of an arbitrary $\Z/2^{[s]}$-framed immersion
admits a natural $\Z/2^{[s+1]}$-framed immersion. Thus, the manifold of self-intersections
yields a triple $(h, \Lambda,
\zeta_L)$, where $h: L^{n-k2^{s}} \looparrowright \R^n$ is an immersion,
$\Lambda: \nu_h \cong k \zeta_L$ and $\zeta_L: L^{n-k2^{s}} \to
K(\Z/2^{[s+1]},1)$ is the classifying map of the
$2^{s}$-dimensional bundle $\zeta_L$. We therefore obtain a homomorphism
\begin{eqnarray}\label{6}
\delta^{\Z/2^{[s+1]}}_{k} : Imm^{\Z/2^{[s]}}(n-k2^{s-1},k2^{s-1})
\to
 Imm^{\Z/2^{[s+1]}}(n-k2^{s},k2^{s}),
\end{eqnarray}
$s \ge 1$, assigning to the normal cobordism class $[(g,\Psi,\eta_N)]$ the normal cobordism
class $[(h, \Lambda,
\zeta_L)]$.

In this formula the positive integer $k$ indicates the multiplicity of the framing and, for $s=1$,
$k$ is equal to the codimension of the immersion. In this case
 the index
$\Z/2^{[s]}$ is replaced by the index  $sf$.

A subgroup $i_{[s+1]}: \Z/2^{[s]} \subset \Z/2^{[s+1]}$ is defined as the subgroup of transformations
of the subspace $Lin(\e_1, \dots, \e_{2^{s-1}})=\R^{2^{s-1}} \subset
\R^{2^{s}}$, generated by the first $2^{s-1}$ basis
vectors $\{\e_1, \dots, \e_{2^{s-1}}\}$, and acting as the identity on the remaining basis vectors.

A subgroup of index 2
\begin{eqnarray}\label{9}
\bar i_{[s+1]}: \Z/2^{[s]} \times \Z/2^{[s]} \subset \Z/2^{[s+1]}
\end{eqnarray}
is defined as the subgroup of transformations leaving invariant each subspace
in the set $\Upsilon_2$.

The subgroup ($\ref{9}$) induces a double covering $\pi_{[s+1]}
: K(\Z/2^{[s]} \times \Z/2^{[s]},1) \to K(\Z/2^{[s+1]},1)$.
The characteristic mapping $\zeta_L: L^{n-k2^{s}} \to
K(\Z/2^{[s+1]},1)$ induces
a double covering $\pi_{[s+1],L} :
\bar L^{n-k2^{s}} \to L^{n-k2^{s}}$ from the covering $\pi_{[s+1]}$
over the classifying space. The double covering $\pi_{[s+1],L}$ can be defined geometrically,
namely it coincides with the canonical double covering of the
manifold $L^{n-k2^{s}}$ of points of self-intersection of the $\Z/2^{[s]}$-framed immersion
$(g,\Psi,\eta_N)$ (see [A1], section 1, formula (3)).

The projection $p_{[s]}: \Z/2^{[s]} \times \Z/2^{[s]} \to
\Z/2^{[s]}$ onto the first factor induces
a mapping
 $p_{[s]}: K(\Z/2^{[s]} \times \Z/2^{[s]},1) \to
 K(\Z/2^{[s]},1)$.

For the manifold of self-intersections $(h,\Lambda,\zeta_L)$ of an arbitrary $\Z/2^{[s]}$--
framed immersion $(g, \Psi,
\eta_N)$, we consider the double covering $\bar \zeta_L: \bar
L^{n-k2^{s}}_{[s]} \to K(\Z/2^{[s]} \times \Z/2^{[s]},1)$  over the classifying mapping $\zeta_L$,
 which is induced from
the covering $\pi_{[s+1],L}$. This covering coincides with the canonical double covering
over the classifying mapping $\zeta_L:
L^{n-k2^s} \to K(\Z/2^{[s+1]},1)$, which is defined by geometric considerations. The characteristic class
$(p_{[s]}
\circ \bar \zeta_L)^{\ast}(\tau_{[s]}\times \tau_{[s]}) \in H^{2^{s}}(\bar
L^{n-k2^{s}}_{[s]};\Z/2)$,
 $\tau_{[s]} \in H^{2^{s-1}}(K(\Z/2^{[s]},1);\Z/2)$
coincides with the characteristic class $\pi_{[s+1],L} \circ \bar
\zeta_{L}(\tau_{[s+1]} )$.

We define the mapping $i_{tot}=i_{[s]}\circ \dots \circ
i_{[3]}$ from the tower:
\begin{eqnarray}\label{pp}
 K(\D,1) \stackrel{i_{[3]}}{\longrightarrow} K(\Z/2^{[3]},1)
\stackrel{i_{[4]}}{\longrightarrow}  \dots
\stackrel{i_{[s]}}{\longrightarrow} K(\Z/2^{[s]},1)
\stackrel{i_{[s+1]}}{\longrightarrow} K(\Z/2^{[s+1]},1).
\end{eqnarray}
There is defined a tower of canonical double coverings
\begin{eqnarray}\label{pipi}
\bar L^{n-k2^{s}}_{[2]} \stackrel{\pi_{[3]}}{\longrightarrow} \bar
L^{n-k2^{s}}_{[3]} \stackrel{\pi_{[4]}}{\longrightarrow} \dots
\stackrel{\pi_{[s]}}{\longrightarrow} \bar L^{n-k2^{s}}_{[s]}
\stackrel{\pi_{[s+1]}}{\longrightarrow} L^{n-k2^s}.
\end{eqnarray}
This tower of coverings is endowed with characteristic mappings to the diagram
($\ref{pp}$). There is defined a sequence of characteristic classes
\begin{eqnarray}\label{bzbz}
\bar \zeta_{[2],L}^{\ast}(\tau_{[2]}) \in H^2(\bar L^{n-k2^{s}}_{[2]};\Z/2), \dots,
\zeta_{[s+1],L}^{\ast}(\tau_{[2]}) \in H^{2^{s}}( L^{n-k2^{s}};\Z/2).
\end{eqnarray}
Each element in this sequence is induced from the characteristic class of
the corresponding universal space in ($\ref{pp}$). We denote by
\begin{eqnarray}\label{tot}
\pi_{tot}=\pi_{[s]}\circ \dots \circ \pi_{[3]} : \bar L^{n-k2^{s}}_{[2]} \to L^{n-k2^s}
\end{eqnarray}
the covering defined as the composition of the coverings in the
diagram $(\ref{pipi})$.  The tower of coverings $(\ref{pipi})$ and
the sequence of characteristic classes $(\ref{bzbz})$ are defined not only for a $\Z/2^{[s+1]}$--
framed immersion which occur as the parametrization of a manifold of self-intersection of a suitable
$\Z/2^{[s]}$--framed immersions, but also for an arbitrary $\Z/2^{[s+1]}$--framed immersion.

\begin{definition}
The Kervaire invariant $\Theta_{\Z/2^{[s+1]}}^k$ of an arbitrary $\Z/2^{[s+1]}$--
framed immersion $(h,\Lambda,\zeta_L)$ is defined by the following formula:
\begin{eqnarray}\label{thetaz2^3}
 \Theta^{\Z/2^{[s+1]}}_k  (h,\Lambda,\zeta_L) = \langle (\bar \zeta_{[2],L}^{\ast}(\tau_{[2]}))^{\frac{n-k2^{s}}{2}};[\bar
L_{[2]}] \rangle,
\end{eqnarray}
where $[\bar L_{[2]}]$ denotes the fundamental class of the covering manifold in the
sequence ($\ref{pipi}$).
\end{definition}

The invariant just constructed defines a homomorphism
$\Theta^{\Z/2^{[s]}}_k: Imm^{\Z/2^{[s]}}(n,n-k2^{d-1}) \to \Z/2$,
which is included in the following commutative diagram:

\begin{eqnarray}\label{77}
\begin{array}{ccc}
Imm^{\Z/2^{[s]}}(n-k2^{s-1},k2^{s-1}) & \stackrel{\Theta^{\Z/2^{[s]}}_k}{\longrightarrow} & \Z/2  \\
\downarrow             \delta^{\Z/2^{[d+1]}}_k  &  &  \| \\
Imm^{\Z/2^{[s+1]}}(n-k2^{s},k2^{s}) &
\stackrel{\Theta^{\Z/2^{[s+1]}}_{k}}{\longrightarrow} & \Z/2.  \\
\end{array}
\end{eqnarray}

For an arbitrary $s \ge 2$ the diagram
$(\ref{77})$ is commutative. For $s=2$ this is proved
in [A1], Lemma 11, for $s>2$ the proof is analogous.

In section 2 we define the concept of an  $\I_{b \times
\bb}$--structure (abelian structure) on a skew-framed immersion
representing an element of the cobordism group $Imm^{sf}(n-k,k)$.
It is proved in Theorem $\ref{th6}$ that under an appropriate dimensional
restriction and modulo elements of odd order, an arbitrary
cobordism class of skew-framed immersions admits an $\I_{b \times
\bb}$-structure. The hypothesis of this theorem presupposes the
existence of a compression of characteristic classes of
skew-framed immersions, see Definition $\ref{5}$. The compression theorem
$\ref{comp}$ will be proved somewhere else.

Let us define a positive integer $\sigma$ (c.f. (1)[A1]) by
the following formula:
\begin{eqnarray}\label{sigma}
\sigma = \left[\frac{\ell}{2} \right]-1.
\end{eqnarray}
In particular, for  $\ell=12$, $\sigma=5$. In section 3 we
formulate the concept of an $\I_a \times \II_a$--structure
(bicyclic structure) on a $\Z/2^{[3]}$--framed immersion. In
Corollary $\ref{cor13}$ it is proved
that, under the conditions of Theorem $\ref{th6}$ (in this theorem
the natural number $n$ can be taken to be greater or equal to
$254$ or greater, if the Compression Theorem for $q=16$ is satisfied; this condition corresponds to the following: the adjoint element in the stable homotopy group of spheres belongs to the image of the suspension of the order $17$, i.e. is borne on the sphere of the dimension $S^{2^{\ell}-18}$) an
arbitrary element of the group
\begin{eqnarray}\label{ImD4}
\mathrm{Im}(\delta^{\Z/2^{[3]}}_{k} \circ \delta^{\Z/2^{[2]}}_{k}:
Imm^{sf}(n-k,k) \to Imm^{\Z/2^{[3]}} (n-4k,4k)),
\end{eqnarray}
\begin{eqnarray}\label{check}
k=\frac{n-m_\sigma}{16}, \quad m_{\sigma}=2^{\sigma}-2, \sigma \ge
5, n \ge 254
\end{eqnarray}
is represented by a $\Z/2^{[3]}$--framed immersion with bicyclic structure. For
such an immersion the Kervaire invariant can be evaluated in terms of an
$\I_a \times \II_a$--characteristic class of the manifold of self-intersections.

In section 4 we formulate the concept of a $ \Q \times
\Z/4$--structure (quaternionic-cyclic structure or briefly quaternionic structure) on a
$\Z/2^{[4]}$--framed immersion. In Corollary $\ref{cor19}$ it is
proved that, under the conditions of Theorem 6, an arbitrary
element of the group

$$\mathrm{Im}(\delta^{\Z/2^{[4]}}_{k} \circ
\delta^{\Z/2^{[3]}}_{k} \circ \delta^{\Z/2^{[2]}}_{k} :
Imm^{sf}(n-k,k) \to$$
\begin{eqnarray}\label{ImZZZ}
 Imm^{\Z/2^{[4]}}( n-8k,8k))
\end{eqnarray}
is represented by a $\Z/2^{[4]}$--framed immersion with
quaternionic-cyclic structure. For such an immersion the Kervaire
invariant can be evaluated in terms of a $\Q \times
\Z/4$--characteristic class of the manifold of self-intersections.

The following two diagrams explain the plan of the proof. In the
first diagram the structure groups $\Z/2^{[s]}$, $2 \le s \le 5$
of parameterizations of self-intersections of immersions with
structure group $\Z/2^{[s-1]}$  are given, as well as the names of
the structures on immersions corresponding to each subgroup:

\begin{eqnarray}\label{A3}
\begin{array}{cccc}
\I_{b \times \bb} & {\subset} & \Z/2^{[2]} & \frac{Abelian}{structure} \\
\downarrow & & i_{[3]} \downarrow &\\
\H_{b \times \bb} &  \subset &  \Z/2^{[3]} & \frac{cyclic-Abelian}{structure}\\
\downarrow &  & i_{[4]} \downarrow  &\\
\J_a \times \JJ_a &   \subset &  \Z/2^{[4]} & \frac{bicyclic}{structure} \\
 \downarrow & &  i_{[5]} \downarrow &\\
\Q \times \Z/4 &  \subset &  \Z/2^{[5]}& \frac{quaternionic-cyclic}{structure} \\
\end{array}
\end{eqnarray}

In the following diagram the natural homomorphisms of cobordism groups
of immersions that will be used are shown, and the Kervaire invariants on
each of these groups are indicated:

\begin{eqnarray}\label{A4}
\begin{array}{ccccc}
Imm^{\Z/2^{[2]}}(n-2,2) &
\stackrel{J^{\Z/2^{[2]}}_{k}}{\longrightarrow} &
Imm^{\Z/2^{[2]}}(n-2k,2k) & \stackrel{\Theta^{\Z/2^{[2]}}_{k}}{\longrightarrow} & \Z/2 \\
\downarrow \delta^{\Z/2^{[3]}} & &  \downarrow \delta^{\Z/2^{[3]}}_{k}  & & \parallel \\
Imm^{\Z/2^{[3]}}(n-4,4) &
\stackrel{J^{\Z/2^{[3]}}_{k}}{\longrightarrow} &
Imm^{\Z/2^{[3]}}(n-4k,4k) &
\stackrel{\Theta^{\Z/2^{[3]}}_{k}}{\longrightarrow} & \Z/2 \\
\downarrow \delta^{\Z/2^{[4]}} & &  \downarrow \delta^{\Z/2^{[4]}}_{k}  && \parallel \\
Imm^{\Z/2^{[4]}}(n-8,8) &
\stackrel{J^{\Z/2^{[4]}}_{k}}{\longrightarrow} &
Imm^{\Z/2^{[4]}}(n-8k,8k)  &
\stackrel{\Theta_{\Z/2^{[4]}}^{k}}{\longrightarrow} & \Z/2 \\
\downarrow \delta^{\Z/2^{[5]}} & &  \downarrow \delta^{\Z/2^{[5]}}_{k}  && \parallel \\
Imm^{\Z/2^{[5]}}(n-16,16) &
\stackrel{J^{\Z/2^{[5]}}_{k}}{\longrightarrow} &
Imm^{\Z/2^{[5]}}(n-16k,16k) &
\stackrel{\Theta^{\Z/2^{[5]}}_{k}}{\longrightarrow} & \Z/2 \\
\end{array}
\end{eqnarray}

In view of the commutativity of this diagram, it suffices to show that the
Kervaire invariant defined in the last row of the diagram is zero. This is
proved with the help of the concept of biquaternionic structure.

The structure group of a $\Z/2^{[5]}$--framed immersion that
parameterizes the manifold of self-intersection points of a
$\Z/2^{[4]}$--framed immersion contains the subgroup $\Q \times
\Z/4$, the direct product of the quaternion group $\Q$ of order 8
and the cyclic group $\Z/4$ of order 4.

The cohomology group $H^{4i}(\Q;\Z)$, for $i \ge 1$, (see  [At],
section 13) contains a characteristic class of order 8. In the
part 1 (see [A1]) the Hopf invariant of a corresponded skew-framed
immersion is estimated using this characteristic class. The Main
Theorem (see section 6) is proved analogously. 
\[  \]

The author thanks Prof. M. Mahowald (2005), Prof. R. Cohen (2007),
and Prof. P. Landweber  for discussions, Prof. A.A. Voronov for an
invitation to lecture at the University of Minnesota (2005), and Prof. V.
Chernov for an invitation to lecture at Dartmouth College (2009).
This paper was begun at the seminar of M.M. Postnikov in 1998. The
paper is dedicated to the memory of Prof. Yu.P. Solov'ev. The compression
theorem was proved in the seminar of A.S. Mishchenko.

\section{Geometric control of the manifold of
\mbox{self-intersections} of a skew-framed immersion}

In this and the following sections, we shall make use of the
cobordism groups $Imm^{sf}(n-k,k)$, $Imm^{\Z/2^{[2]}}(n-2k,2k)$.
It is well known that if the first argument in parentheses,
denoting the dimension of the immersed manifold, is strictly
positive, then the indicated group is a finite group.

The dihedral group $\Z/2^{[2]}$ is defined by its presentation:
$$\{a,b \,\vert \, a^4 = b^2 = e, [a,b]=a^2\}.$$

This group is a subgroup of the orthogonal group $O(2)$, namely,
the group of orthogonal transformations of the standard Euclidean
plane $Lin(\e_1,\e_2)$, preserving the unordered pair of lines generated by a basis of
orthogonal unit vectors $\{\e_1,\e_2\}$. The element $a$ is
represented by a rotation of the plane through an angle
$\frac{\pi}{2}$. The element $b$ is represented by a reflection of
the plane with respect to the line $l_1 = Lin(\e_1 +
\e_2)$, generated by the vector $\e_1 +
\e_2$.

Consider the subgroup $\I_{b \times \bb} = \I_b \times \II_b \subset \Z/2^{[2]}$
 of
the dihedral group, generated by the elements $\{b,ba^2\}$. Notice
that this is an elementary $2$-group of rank $2$. This is the
subgroup of $O(2)$ consisting of transformations preserving
individually each of the lines $l_1$, $l_2$ in the directions of
the vectors $\f_1=\e_1+\e_2$, $\f_2=\e_1-\e_2$ respectively. The
cohomology group $H^1(K(\I_b \times \II_b,1);\Z/2)$ is also an
elementary $2$-group with two generators. We now describe these
generators.

Let us define the cohomology classes
\begin{eqnarray}\label{kappab}
\kappa_b \in H^1(K(\I_b \times \II_b,1);\Z/2)), \quad \kappa_{\bb} \in H^1(K(\I_b \times \II_b,1);\Z/2)).
\end{eqnarray}
We denote by $p_{b}: \I_{b \times \bb} \to \I_b$ the homomorphism,
whose kernel consists of the reflection with respect to the
bisector of the second coordinate angle and the identity. Define
$\kappa_{b} = p_{b}^{\ast}(t_{b})$,
 where $e \ne t_{b} \in H^1(K(\I_{b},1);\Z/2)
\cong \Z/2$. Denote by $p_{\bb}: \I_b \times \II_b \to \II_b$ the
homomorphism, whose kernel consists of the reflection with respect
to the bisector of the first coordinate angle and the identity,
or equivalently whose kernel consists of the composition of the
central symmetry and the symmetry with respect to the second
coordinate angle and the identity. Define $\kappa_{\bb} =
p_{\bb}^{\ast}(t_{\bb})$, where $e \ne t_{\bb} \in
H^1(K(\II_{b},1);\Z/2) \cong \Z/2$.

%The cohomology class $\kappa_b,
%\kappa_{\bb} \in H^2(K(\I_{b \times \bb},1);\Z/2))$ is well defined.

We next define a group
 $\I_{b \times \bb} \int_{\chi^{[2]}} \Z$ and an epimorphism $\I_{b \times \bb} \int_{\chi^{[2]}} \Z \to \Z/2^{[2]}$.
Consider the automorphism
\begin{eqnarray}\label{chidd}
\chi^{[2]}: \I_{b \times \bb} \to \I_{b \times \bb}
\end{eqnarray}
determined by the external conjugation of the subgroup
$\I_{b \times \bb} \subset \Z/2^{[2]}$ by the element $ab \in \Z/2^{[2]}$,
this element correspond to the reflection with respect to the line  $Lin(\e_1)$..

Define the automorphism (the notation is similar)
\begin{eqnarray}\label{chi2}
\chi^{[2]}: \Z/2^{[2]} \to \Z/2^{[2]},
\end{eqnarray}
by the reflection in the line $Lin(\e_1)$.
It is easy to see that the inclusion $\I_{b \times \bb} \subset \Z/2^{[2]}$ commutes with $(\ref{chidd})$ and
$(\ref{chi2})$ such that the corresponding diagram is commutative.

Define the group
\begin{eqnarray}\label{chi2Z}
\I_{b \times \bb} \int_{\chi^{[2]}} \Z
\end{eqnarray}
as the factorgroup of the group $\I_{b \times \bb} \ast \Z$ (the free product of the groups $\I_{b \times \bb}$ and $\Z$) by the relation $zxz^{-1}=\chi^{[2]}(x)$,
where $z \in \Z$ is the standard generator,
$x \in \I_{b \times \bb}$ is an arbitrary element. The group $(\ref{chi2Z})$ is a particular example of a semi-direct product of groups
$A \rtimes_{\phi} B$, $A=\I_{b \times \bb}$, $B=\Z$, where $\phi: B \to \mathrm{Aut}(A)$ is a homomorphism and the set $A \times B$ is equipped with the
operation $(a_1,b_1) \ast (a_2,b_2) \mapsto (a_1 \phi_{b_1}(b_2),b_1b_2)$.

The classifying space
 $K(\I_{b \times \bb} \int_{\chi^{[2]}} \Z,1)$ is a skew-product of the standard circle  $S^1$ and the space $K(\I_{b \times \bb},1)$, the map $K(\I_{b \times \bb},1) \to K(\I_{b \times \bb},1)$,
of the shift in the cyclic covering, associated with the cyclic covering over  $K(\I_{b \times \bb} \int_{\chi^{[2]}} \Z,1)$
is induced by the automorphism $\chi^{[2]}$.
Denote the standard  fibration by 
\begin{eqnarray}\label{pbb}
p_{b \times \bb}: K(\I_{b \times \bb}\int_{\chi^{[2]}}\Z,1) \to S^1.
\end{eqnarray}

Take a marked point
 $pt_{S^1}   \in S^1$.
and define the subspace 
\begin{eqnarray}\label{inclK}
K(\I_{b \times \bb},1) \subset K(\I_{b \times \bb} \int_{\chi^{[2]}} \Z,1) 
\end{eqnarray}
as the inverse image of the marked point
$pt_{S^1}$ by the mapping $(\ref{pbb})$.

Let us consider the homology groups
$H_{i}(K(\I_{b\times \bb} \int_{\chi^{[2]}} \Z,1);\Z/2)$,
$H_{i}(K(\I_{b \times \bb} \int_{\chi^{[2]}} \Z,1);\Z)$
(below the coefficients  $\Z/2$ are omitted).

 The standard basis of the group $H_{i}(K(\I_{b \times \bb} \int_{\chi^{[2]}} \Z,1);\Z)$
is sufficiently complicated and we do not need its description. The basis of the group $H_{i}(K(\I_{b \times \bb},1);\Z)$ is described by
the K\"unneth formula:
\begin{eqnarray}\label{Kunnd}
0 \to \bigoplus_{i_1+i_2=i} H_{i_1}(K(\I_b,1);\Z) \otimes H_{i_2}(K(\II_b,1);\Z) \longrightarrow
H_i(K(\I_{b \times \bb},1);\Z)
\end{eqnarray}
$$\longrightarrow \bigoplus_{i_1+i_2=i-1} Tor^{\Z}(H_{i_1}(K(\I_b,1);\Z),H_{i_2}(K(\II_b,1);\Z) \to 0.$$

The standard basis of the group  $H_{i}(K(\I_{b \times \bb},1))$ is following: 
$$x\otimes y / (x \otimes y)+(y \otimes x),$$
where $x \in H_{j}(K(\I_b,1))$, $y \in H_{i-j}(K(\II_b,1))$.

In particular, in the case of odd $i$
the group $H_{i}(K(\I_{b \times \bb},1);\Z)$ contains the fundamental classes of the following submanifolds  $\RP^{i} \times pt \subset \RP^i \times \RP^i
\subset K(\I_b,1) \times K(\II_b,1)$, $pt \times \RP^{i} \subset \RP^i \times \RP^i
\subset K(\I_b,1) \times K(\II_b,1)$.
Denote the corresponding elements by
\begin{eqnarray}\label{tdtdd}
t_{b,i} \in H_{i}(K(\I_{b \times \bb},1);\Z),
t_{\bb,i} \in H_{i}(K(\I_{b \times \bb},1);\Z).
\end{eqnarray}

Let us define the analogous homology groups 
with local coefficients system:
\begin{eqnarray}\label{Z2lok}
H_{i}(K(\I_{b \times \bb} \int_{\chi^{[2]}} \Z,1);\Z/2[\Z/2]),
\end{eqnarray}
\begin{eqnarray}\label{Zlok}
H_{i}(K(\I_{b \times \bb} \int_{\chi^{[2]}} \Z,1);\Z[\Z/2]).
\end{eqnarray}

The following epimorphism 
\begin{eqnarray}\label{epid}
p_{b \times \bb}: \I_{b \times \bb} \int_{\chi^{[2]}} \Z \to \Z,
\end{eqnarray}
is well defined by the formula $x \ast y \mapsto y$, $x \in \I_{b \times \bb}$, $y \in \Z$..
The following epimorphism
\begin{eqnarray}\label{epid2}
\I_{b \times \bb} \int_{\chi^{[2]}} \Z \to \Z \to \Z/2
\end{eqnarray}
is defined by the formula
$p_{b \times \bb} \pmod{2}$.

Let us define the group
$(\ref{Zlok})$. Consider the group ring  $\Z[\Z/2] = \{a + bt\}$, $a,b \in \Z$, $t \in \Z/2$.
The generator $t \in \Z[\Z/2]$ is represented by the involution  $\chi^{[2]}: K(\I_{b \times \bb} \int_{\chi^{[2]}} \Z,1) \to
K(\I_{b \times \bb} \int_{\chi^{[2]}} \Z,1)$, the restriction of this involution on the subspace
$K(\I_{b \times \bb},1) \subset K(\I_{b \times \bb} \int_{\chi^{[2]}} \Z,1)$ 
is the reflection, which is induced by the permutation automorphism
$\I_b \times \II_b \to \II_{\bb} \times \I_b$.
Consider the following local system of the coefficients
$\rho_t: \Z/2[\Z/2] \to Aut(K(\I_{b \times \bb},1) \subset K(\I_{b \times \bb} \int_{\chi^{[2]}} \Z,1)$,
which identifies the chain  $(a + bt)\sigma$, with a support on a simplex  $\sigma \subset K(\I_{b \times \bb} \int_{\chi^{[2]}} \Z,1)$ with the chain  $(at + b)\chi^{[2]}(\sigma)$. The group  $(\ref{Zlok})$ is
well defined. 
The group $(\ref{Z2lok})$ is defined analogously.

The description of the groups
 $(\ref{Z2lok})$, $(\ref{Zlok})$ are sufficiently complicated and we will not use this description.
Let us define a subgroup 
\begin{eqnarray}\label{DZlok}
D_{i}(\I_{b \times \bb};\Z[\Z/2]) \subset H_{i}(K(\I_{b \times \bb} \int_{\chi^{[2]}} \Z,1);\Z[\Z/2])
\end{eqnarray}
by the formula: $D_{i}(\I_{b \times \bb};\Z[\Z/2]) = \mathrm{Im}(H_i(K(\I_{b \times \bb},1);\Z[\Z/2]) \to
H_{i}(K(\I_{b \times \bb} \int_{\chi^{[2]}} \Z,1);\Z[\Z/2]))$, where the homomorphism
$H_i(K(\I_{b \times \bb},1);\Z[\Z/2]) \to
H_{i}(K(\I_{b \times \bb} \int_{\chi^{[2]}} \Z,1);\Z[\Z/2])$ is induced by the inclusion of the subgroup
$(\ref{incl})$. The natural epimorphism
\begin{eqnarray}\label{epilok}
H_{i}(K(\I_{b \times \bb},1);\Z[\Z/2]) \to D_{i}(\I_{ b \times \bb};\Z[\Z/2])
\end{eqnarray}
is well defined.

The following natural homomorphism
 $H_{i}(K(\I_{b \times \bb},1);\Z)) \otimes \Z[\Z/2] \to
H_{i}(K(\I_{b \times \bb},1);\Z[\Z/2])$, is an isomorphism by the universal coefficients formula. 
This calculation do not use the structure of 
$\Z[\Z/2]$--module, and use the additive isomorphism   $\Z[\Z/2] \equiv \Z \oplus \Z$
and additivity of the functor
$Tor^{\Z}$. Analogously, $H_{i}(K(\I_{b \times \bb},1);\Z/2)) \otimes \Z/2[\Z/2] \equiv
H_{i}(K(\I_{b \times \bb},1);\Z/2[\Z/2])$.

The subgroup
$(\ref{DZlok})$ is generated by the following elements:
$X + Yt$, $X,Y \in H_{i}(K(\I_{b \times \bb},1);\Z)$.
The following equivalence relation determines the equivalence of the two representatives:
$X \equiv \chi^{[2]}_{\ast}(X)t$, where the automorphism
\begin{eqnarray}\label{ch2}
\chi^{[2]}_{\ast}: H_{i}(K(\I_{b \times \bb},1);\Z) \to H_{i}(K(\I_{b \times \bb},1);\Z)
\end{eqnarray}
is induced by the automorphism
$(\ref{chidd})$. The automorphism
 $(\ref{chidd})$ induces also the automorphism
\begin{eqnarray}\label{ch2lok}
\chi^{[2]}_{\ast}: H_{i}(K(\I_{b \times \bb},1);\Z[\Z/2]) \to H_{i}(K(\I_{b \times \bb},1);\Z[\Z/2]).
\end{eqnarray}

The subgroup
\begin{eqnarray}\label{DZ2lok}
D_{i}(\I_{b \times \bb};\Z/2[\Z/2]) \subset  H_{i}(K(\I_{b \times \bb} \int_{\chi^{[2]}} \Z,1);\Z/2[\Z/2])
\end{eqnarray}
is defined analogously to
 $(\ref{DZlok})$. The description of this subgroup is more simple, because
 this subgroup is generated by the elements
 $X + tY$,
$X = x \otimes y, Y = x' \otimes y'$, where $\chi^{[2]}_{\ast}(x \otimes y) = y \otimes x$.
%This follows from the fact: the kernel of the homomorphism 
% $(\ref{epilok})$ is generated by elements
%$X - \chi^{[2]}_{\ast}(X)t$.

%Define the subgroup  $D_{i}(\I_a \times \II_a;\Z/2[\Z/2])$
%analogously to the subgroup
%$(\ref{DZlok})$. The  description of this subgroup is simpler, because we have
%$X = x \otimes y, Y = x' \otimes y'$, $\chi^{[2]}_{\ast}(x \otimes y) = y \otimes x$.

Define the homomorphism
\begin{eqnarray}\label{delta}
\Delta^{[2]}: H_{i}(K(\I_{b \times \bb},1);\Z[\Z/2]) \to H_{i}(K(\I_{b \times \bb},1);\Z)
\end{eqnarray}
by the formula $\Delta^{[2]}(X + Yt) = X + \chi^{[2]}_{\ast}(Y)$.
Let us prove that the homomorphism $(\ref{delta})$ is naturally factorized to the homomorphism
\begin{eqnarray}\label{deltaD}
\Delta^{[2]}: D_{i}(\I_{b \times \bb};\Z[\Z/2]) \to H_{i}(K(\I_{b \times \bb},1);\Z),
\end{eqnarray}
which we denote the same. This follows from the following observation: the kernel of the homomorphism
$(\ref{epilok})$ is generated by elements
$X - \chi^{[2]}_{\ast}(X)t$.

Define the analogously homomorphisms with $\Z/2$--coefficients:
\begin{eqnarray}\label{deltaD2}
\Delta^{[2]}: D_{i}(\I_{b \times \bb};\Z/2[\Z/2]) \to H_{i}(K(\I_{b \times \bb},1)).
\end{eqnarray}
Let us consider the following composition of the homomorphisms:
\begin{eqnarray}\label{compos}
H_{i}(K(\I_{b \times \bb},1);\Z) \to H_{i}(K(\I_{b \times \bb},1);\Z[\Z/2]) \to
\end{eqnarray}
$$
D_{i}(\I_{b \times \bb};\Z[\Z/2]) \to
H_{i}(K(\I_{b \times \bb},1);\Z),$$
where the left homomorphism is the natural inclusion, the middle is the homomorphism
 $(\ref{epilok})$, the right homomorphism in this composition is  $(\ref{deltaD})$.
It is easy to check that the composition coincides with the identity homomorphism.

Let us define the forgetful homomorphism 
\begin{eqnarray}\label{forg}
forg_{\ast}: H_{i}(K(\I_{b \times \bb} \int_{\chi^{[2]}} \Z,1);\Z[\Z/2]) \to H_{i}(K(\I_{b \times \bb} \int_{\chi^{[2]}} \Z,1);\Z).
\end{eqnarray}
This homomorphism is induced by the forgetful mapping of the local coefficient system, analogously to the homomorphism
 $\Delta^{[2]}$.

Let us consider the homomorphism 
 $(\ref{epid})$ and let us consider the mapping  $(\ref{pbb})$, which is associated with this homomorphism. 
 Assume  $\theta \in H^1(S^1;\Z[\Z/2])$ is the generator.
Take the element $p_{b \times \bb}^{\ast}(\theta) \in H^1(K(\I_{b \times \bb} \int_{\chi^{[2]}} \Z,1);\Z[\Z/2])$ and
define the obstruction $(\ref{ox})$ by the formula $o(x) = x \cap p_{b \times \bb}^{\ast}(\theta)$.
It is easy to check, that  $\chi^{[2]}_{\ast}(o(x)) = o(x)$, where the automorphism  $\chi^{[2]}_{\ast}$ 
is defined by the formula  $(\ref{ch2lok})$. 

Consider the subgroup $(\ref{DZlok})$ and the epimorphism $(\ref{epilok})$.
Take an element $x \in H_{i}(K(\I_{b \times \bb} \int_{\chi^{[2]}} \Z,1);\Z/2[\Z/2])$. Let us describe the total obstruction
\begin{eqnarray}\label{ox}
o(x) \in D_{i-1}(\I_{b \times \bb};\Z/2[\Z/2])
\end{eqnarray}
for the following inclusion: $x \in D_i(\I_{b \times \bb};\Z/2[\Z/2])
\subset H_{i}(K(\I_{b \times \bb} \int_{\chi^{[2]}} \Z,1);\Z/2[\Z/2])$. Let us consider the generator $\theta \in H^1(S^1;\Z[\Z/2])$.
Take the element $p_{b \times \bb}^{\ast}(\theta) \in H^1(K(\I_{b \times \bb} \int_{\chi^{[2]}} \Z,1);\Z/2[\Z/2])$ and
define the obstruction $(\ref{ox})$ by the formula $o(x) = x \cap p_{b \times \bb}^{\ast}(\theta)$.

From the definition of the obstruction $(\ref{ox})$ is obvious that the value of this obstruction one can calculate as following: apply to a singular cycle with local coefficient, which represents a prescribed homology class, the homomorphism  $(\ref{forg})$, and then intersects the  singular cycle of the space 
$K(\I_{b \times \bb} \int_{\chi^{[2]}} \Z,1)$ with the subspace  $(\ref{inclK})$.

Evidently, $\chi^{[2]}_{\ast}(o(x)) = o(x)$, where the automorphism $\chi^{[2]}_{\ast}$ is given by the formula $(\ref{ch2lok})$.

In the case $i=2s$ basis elements
$y \in D_{i-1}(\I_{b \times \bb};\Z[\Z/2])$, which satisfies the equation $\chi^{[2]}_{\ast}(y) = y$, are the following:

1. $y=r$, $r= t_{b,2s-1} + t_{\bb,2s-1}$,  where $t_{d,2s-1},t_{\bb,2s-1} \in H_{2s-1}(K(\I_{b \times \bb},1))$ are given by the formula  $(\ref{tdtdd})$.

2. $y=z(i_1,i_2)$, where $z(i_1,i_2) = tor(r_{b,i_1}, r_{\bb,i_2}) + tor(r_{b,i_2}, r_{\bb,i_1})$, $i_1,i_2 \equiv 1 \pmod{2}$,
$i_1 + i_2 = 2s-2$,
$tor(r_{b,i_1}, r_{\bb,i_2}) \in Tor^{\Z}(H_{i_1}(K(\I_b,1);\Z),H_{i_2}(K(\II_b,1);\Z))$.
The  elements  $\{z(i_1,i_2)\}$ are in the kernel of the homomorphism
\begin{eqnarray}\label{coeff}
D_{i-1}(\I_{b \times \bb};\Z[\Z/2]) \to D_{i-1}(\I_{b \times \bb};\Z/2[\Z/2]),
\end{eqnarray}
given by the homomorphism $\Z \to \Z/2$ of coefficients.

Elements
 $R, \quad Z(i_1,i_2) \in H_{2s}(K(\I_{b \times \bb} \int_{\chi^{[2]}} \Z,1);\Z[\Z/2])$,
which satisfies the relation  $o(R)=r$, $o(R)=r$, $o(Z(i_1,i_2))=z(i_1,i_2)$
are defined the following way. Each element are given by cycle, which is represented by the the product of the corresponding
 $2s-1$--cycle $f: C_{2s-1} \to K(\I_{b \times \bb},1)$ with the circle.
The mapping of the cycle into  $K(\I_{b \times \bb} \int_{\chi^{[2]}} \Z,1)$
is determined by the composition of the mapping
$f \times id: C_{2s-1} \times S^1 \to K(\I_{b \times \bb},1) \times S^1$
with the standard $2$-sheeted covering
$K(\I_{b \times \bb},1) \times S^1 \to K(\I_{b \times \bb} \int_{\chi^{[2]}} \Z,1)$.

By means of the obstruction $(\ref{ox})$ we shall prove the following lemma.

\begin{lemma}\label{obstr}
The group $H_{2s}(K(\I_{b \times \bb}\int_{\chi^{[2]}}\Z,1);\Z[\Z/2])$
is isomorphic to the direct sum of the subgroup
$ D_{2s}(\I_{b \times \bb};\Z[\Z/2])$ and the subgroup generated by the elements
$R, \quad \{Z(i_1,i_2)\}$. The elements in the subgroup
$\bigoplus_{i_1+i_2=2s} H_{i_1}(K(\I_b,1);\Z) \otimes H_{i_2}(K(\II_b,1);\Z) \subset D_{2s}(\I_b \times \II_b;\Z[\Z/2])$, and the elements
 $R$ generate the image  $\mathrm{Im}(A)$ of the homomorphism
\begin{eqnarray}\label{A}
A: H_{2s}(K(\I_{b \times \bb}\int_{\chi^{[2]}}\Z,1);\Z[\Z/2]) \to
\end{eqnarray}
$$H_{2s}(K(\I_{b \times \bb}\int_{\chi^{[2]}}\Z,1);\Z/2[\Z/2]),$$
which is induced by the reduction of the local coefficients system modulo $2$.
\end{lemma}

\subsubsection*{Proof of Lemma $\ref{obstr}$}
Let $x \in H_{2s}(K(\I_{b \times \bb}\int_{\chi^{[2]}}\Z,1);\Z[\Z/2])$ is not in the subgroup $ D_{2s}(\I_{b \times \bb};\Z[\Z/2])$.
Then $o(x) \ne 0$ (see  $(\ref{ox})$ for $i=2s$. Therefore  $x$ as an element in a residue class with respect to the considered subgroup is expressed by means of the elements
$R, \quad Z(i_1,i_2)$. The subgroup of all values of the obstruction $(\ref{ox})$
 is a direct factor in  $H_{2s}(K(\I_{b \times \bb}\int_{\chi^{[2]}}\Z,1);\Z[\Z/2])$,
 because this obstructions are realized by a linear combination elements
$R, \quad \{Z(i_1,i_2)\}$. The elements  $Z(i_1,i_2)$ belong to $\mathrm{Ker}(A)$ and, therefore these elements are
not a generators of $\mathrm{Im}(A)$.
Lemma $\ref{obstr}$ is proved.
\[  \]

Define the inclusion
 $i_{\I_d, \I_{b \times \bb}}: \I_d \subset \I_b \times \II_b = \I_{b \times \bb}$,
as the diagonal inclusion.
The subgroup coincides with the kernel of the homomorphism
\begin{eqnarray}\label{omega2}
\omega^{[2]}: \I_b \times \II_b \to \Z/2,
\end{eqnarray}
this homomorphism is given by the formula  $(x \times y) \mapsto xy$.

Define the epimorphism
\begin{eqnarray}\label{Phi2}
\Phi^{[2]}: \I_{b \times \bb} \int_{\chi^{[2]}} \Z \to \Z/2^{[2]}
\end{eqnarray}
by the formula:
$\Phi^{[2]}(z)=ab$, $z \in \Z$ is the generator (the element $ab$ corresponds to the reflection of the first vector of the standard base with respect to the standard representation $\Z^{[2]} \subset O(2)$), the restriction $\Phi^{[2]} \vert_{\I_b \times \II_b \times \{0\}}: \I_b \times \II_b \subset \Z/2^{[2]}$ is the standard inclusion.
Therefore $\Phi^{[2]} \vert_{\I_{b \times \bb} \times \{1\}}: \I_b \times \II_b = \I_{b \times \bb} \subset \Z/2^{[2]}$
is the conjugated inclusion by the exterior automorphism in the subgroup $\I_b \times \II_b = \I_{b \times \bb} \subset \Z/2^{[2]}$.

Define
$$ (\Phi^{[2]})^{\ast}(\tau_{[2]}) = \tau_{b \times \bb}, $$
where
$\tau_{b \times \bb} \in H^2(K(\I_b \times \II_b \int_{\chi^{[2]}} \Z),1))$, $\tau_{[2]} \in H^2(K(\Z^{[2]},1)$.

\begin{definition}\label{ab}
Let an element $y 
\in Imm^{\Z/2^{[2]}}(n-2k,2k)$ be represented by a
$\Z/2^{[2]}$--framed immersion $(g,\Psi,\eta_N)$, $g: N^{n-2k}_{b \times \bb}
\looparrowright \R^n$. We say that the $\Z/2^{[2]}$-framed immersion
$(g,\Psi,\eta_N)$ is an $\I_{b \times \bb}$--immersion (abelian
immersion) if the following two conditions are satisfied.

1. The structure mapping $\eta_N: N^{n-2k}_{b \times \bb} \to
K(\Z/2^{[2]},1)$ is represented as the composition of a
mapping
\begin{eqnarray}\label{redu}
\eta_{b \times \bb,N}: N^{n-2k}_{b \times \bb} \to K(\I_{b \times
\bb} \int_{\chi^{[2]}} \Z,1)
\end{eqnarray}
and the mapping $\Phi^{[2]}: K(\I_{b \times
\bb} \int_{\chi^{[2]}} \Z,1)
\to K(\Z/2^{[2]},1)$.

2. Consider the submanifold
\begin{eqnarray}\label{NN}
N_{\eta^{7k}}^{n-16k} \subset N^{n-2k}_{b \times \bb},
\end{eqnarray}
which represents the Euler class
$[\eta_N^{7k}]^{op} \in H_{n-16k}(N^{n-2k}_{b \times \bb};\Z/2)$ of the bundle $7k\eta_N$
(the considered class is Poincar\'e dual to the cohomology class
$\eta_N^{7k} \in H^{7k}(N^{n-2k}_{b \times \bb};\Z/2)$).
It is required that the restriction of the mapping
 $(\ref{redu})$ on the submanifold  $(\ref{NN})$
 is given by the composition of the mapping
 \begin{eqnarray}\label{Neta7}
\eta_{Ab}: N_{\eta^{7k}}^{n-16k} \to K(\I_{b \times \bb},1)
\end{eqnarray}
with the standard inclusion
 $K(\I_{b \times \bb},1)
\subset K(\I_{b \times
\bb} \int_{\chi^{[2]}} \Z,1)$.
 \end{definition}

\begin{definition}\label{strucIdId}
Let a skew-framed immersion $(f,\Xi,\kappa)$, $f: M^{n-k}
\looparrowright \R^n$ represents an element $x \in
Imm^{sf}(n-k,k)$, where $n > 16k$. Let the $\Z/2^{[2]}$--framed
immersion $(g,\Psi,\eta_N)$, $g: N^{n-2k} \looparrowright \R^n$ be
the immersion of the manifold of self-intersections of the
immersion $f$, so $g$ represents the element
$y=\delta_{\Z/2^{[2]}}^{k}(x) \in Imm^{\Z/2^{[2]}}(n-2k,2k)$. We
say that the skew-framed immersion $(f,\Xi,\kappa_M)$ admits an
abelian structure  ($\I_{b \times \bb}$--structure) if
the self-intersection manifold $N^{n-2k}$ of the immersion  $f$ is decomposed
into two components (possibly, non-connected):
\begin{eqnarray}\label{abelcomp} 
N^{n-2k} = N^{n-2k}_{b \times \bb} \cup N^{n-2k}_{[2]}.
\end{eqnarray}
Each component is equipped by a  $\Z/2^{[2]}$--framed immersion into
$\R^n$. The following two conditions are  satisfies:

1. For the characteristic mapping  $\eta_N \vert _{N^{n-2k}_{b \times \bb}}$ of the $\Z/2^{[2]}$--framing over 
first component in the formula  $(\ref{abelcomp})$ is equipped by a reduction into the subspace  $K(\I_{b \times \bb} \int_{\chi^{[2]}} \Z,1) \subset K(\Z/2^{[2]},1)$, defined by the mapping
$(\ref{redu})$. For the mapping $(\ref{redu})$ the property 2 of Definition
$\ref{ab}$ is satisfied. Namely, the restriction of the considered mapping on a submanifold
$N_{\eta^{7k}}^{n-16k} \subset N^{n-2k}_{b \times \bb}$, which is determined by the formula
 $(\ref{NN})$, allows an additional reduction into the subspace 
$K(\I_{b \times \bb},1) \subset K(\I_{b \times
\bb} \int_{\chi^{[2]}} \Z,1)$, given by the formula $(\ref{Neta7})$.

2. The $\Z/2^{[2]}$--framed immersion, which is defined by the restriction of the immersion $g$ 
on the second component in the formula  $(\ref{abelcomp})$, has the trivial Kervaire invariant: the characteristic class
$(\ref{44})$ for the considered component is trivial.
\end{definition}

\subsubsection*{Example}
Assume a skew-framed immersion $(f,\Xi,\kappa_M)$, $f: M^{n-k}
\looparrowright \R^n$,  represents an element $x \in
Imm^{sf}(n-k,k)$, $n > 16k$.  Assume that a $\Z/2^{[2]}$--framed immersion $(g,\Psi,\eta_N)$, $g: N^{n-2k} \looparrowright \R^n$, which is the immersion of self-intersection manifold of $f$,  is a $\I_{b \times
\bb}$--immersion in the sense of Definition $\ref{ab}$. Then the skew-framed immersion
$(f,\Xi,\kappa_M)$ admits an abelian structure, for which the second component in the formula
$(\ref{abelcomp})$ is empty.
\[  \]

\subsubsection*{The justification of the example}
Let us define the map, determined by an
abelian structure by the formula $(\ref{redu})$.
The Conditions 1 and 2 in Definition $\ref{ab}$ implies the conditions 1 and 2 in Definition
$\ref{strucIdId}$. The abelian structure is well defined.
\[  \]

\subsubsection*{The fundamental class with local coefficients of a framed immersion}

Consider a  $\Z/2^{[2]}$--framed immersion 
$(g,\Psi,\eta_N)$, $g: N^{n-2k}_{b \times \bb} \looparrowright \R^n$, 
and assume that a mapping $(\ref{redu})$, which determines a reduction of the characteristic mapping,
is well-defined. Assume that the manifold
$N^{n-2k}_{b \times \bb}$ is connected. Assume that a marked point
$pt \in N^{n-2k-1}_{pt}$ is fixed. Assume that the image of the marked point by the mapping 
$(\ref{redu})$ is a point in the subspace $(\ref{inclK})$.
Let us prove that the image of the fundamental class
 $[N^{n-2k}_{b \times \bb};pt]$ by the mapping  $(\ref{redu})$ determines an element 
\begin{eqnarray}\label{etadd16}
\eta_{b \times \bb,\ast}([N^{n-2k}_{b \times \bb};pt]) \in H_{n-2k}(K(\I_{b \times \bb} \int_{\chi^{[2]}} \Z,1);\Z[\Z/2]).
\end{eqnarray}

Denote by 
\begin{eqnarray}\label{Npt}
N^{n-2k-1}_{pt} \subset N^{n-2k}_{b \times \bb}
\end{eqnarray}
the submanifold, which is defined as the regular preimage of the subspace
 $(\ref{incl})$. The restriction of the reduction mapping on the submanifold
$(\ref{Npt})$ determines a $\I_{b \times \bb}$--reduction of the restriction of the characteristic mapping 
$\eta_N$ to this  
manifold.

Consider the skeleton of the space
 $K(\I_{b \times \bb},1)$, which is realized as a  $O(2)/\I_{b \times \bb}$-bundle 
 over the Grassman manifold  $Gr_{O}(2,n)$ of 
$2$-planes in $n$--dimensional space. Denote this skeleton by
\begin{eqnarray}\label{KK}
KK(\I_{b \times \bb},1) \subset K(\I_{b \times \bb},1).
\end{eqnarray}

The following free involution
\begin{eqnarray}\label{chiK}
\chi^{[2]}: KK(\I_{b \times \bb},1) \to KK(\I_{b \times \bb},1),
\end{eqnarray}
acts on the space $(\ref{KK})$, this involution corresponds to the automorphism 
$(\ref{chidd})$ (and denote the same). Let us consider the skeleton of the pair of spaces
 $(\ref{incl})$ as the cylinder of the involution  $(\ref{KK})$, and denote this cylinder by
\begin{eqnarray}\label{KKint}
KK(\I_{b \times \bb} \int_{\chi^{[2]}} \Z,1) \subset K(\I_{b \times \bb} \int_{\chi^{[2]}} \Z,1).
\end{eqnarray}
The involution 
$(\ref{chiK})$ induces the involution 
\begin{eqnarray}\label{chiKint}
\chi^{[2]}: KK(\I_{b \times \bb}\int_{\chi^{[2]}} \Z ,1) \to KK(\I_{b \times \bb}\int_{\chi^{[2]}} \Z,1),
\end{eqnarray}
which is extended to the involution on the hole space
 $K(\I_{b \times \bb} \int_{\chi^{[2]}} \Z,1)$.

The  universal 
 $\I_{b \times \bb} \int_{\chi^{[2]}} \Z$--bundle over the skeleton $(\ref{KKint})$ is well-defined, 
 denote this bundle by 
$\tau_{b \times \bb, \int}$. Denote the restriction of the bundle  $\tau_{b \times \bb, \int}$ on
the subspace $(\ref{KK})$ by $\tau_{b \times \bb}$.

The mapping
$(\ref{redu})$ determines a $\I_{b \times \bb} \int _{\chi^{[2]}} \Z$--reduction of the characteristic mapping of the framing  $\Psi$ and the framing 
 $\Psi_{b \times \bb}: \nu_g \equiv k \eta_{b \times \bb}^{\ast}
(\tau_{b \times \bb, \int})$ is well-defined. 

Let us visualize this reduction by the following way. The normal bundle 
 $\nu_g$ of the immersion $g$ is decomposed into the Whitney sum of the $k$ copies of a 2-dimensional bundle.
Let us denote the first term of this decomposition by $\nu_{b \times \bb,\int} \subset \nu_g$. 
The bundle $\nu_{b \times \bb,\int}$ is a twisted Whitney sum of the two line bundles $\kappa_b \tilde{\oplus} \kappa_{\bb}$. 
%the restriction of these line bundles over the submanifold  $(\ref{NN})$ is well-defined. 
These two line bundles $\kappa_b$, $\kappa_{\bb}$ are permuted by the parallel transformation over a path
$l \subset N^{n-2k}_{b \times \bb}$, if the projection of the path by the mapping $(\ref{pbb})$ 
represents a generator in  $H_1(S^1)$.

Let us consider an arbitrary cell
$\alpha$ of a regular cell decomposition of the manifold $N^{n-2k}_{b \times \bb}$. Assume that a path  $\phi_{\alpha}$, which is joins the center of the cell $\alpha$ with the marked point $pt \in N^{n-2k}_{b \times \bb}$
is given. The restriction of the bundle  $\nu_{b \times \bb,\int}$ to the cell  $\alpha$ is classified by the mapping  $\eta_{b\times \bb}(\alpha, \phi_{\alpha}): \alpha \to KK(\I_{b \times \bb} \int_{\chi^{[2]}} \Z,1)$. 
A change of 
 $\phi_{\alpha}$ in a residue class of paths with the fixed ends in the group  $H_1(S^1)$ determines another mapping  $\eta_{b \times \bb}(\alpha, \phi_{\alpha})$, which is defined by means of the composition of the classified mapping mapping with the involution  $(\ref{chiKint})$. In the cell complex of the space  $(\ref{KKint})$ with $\Z[\Z/2]$--local coefficients a change of the path, which is attached to the cell $\alpha$, corresponds to the change
 of the base  of the corresponding generator by the multiplication on the element  $t$ and the change of the 
 parametrization mapping for the cell by the involution  $(\ref{chiKint})$. Therefore, the element
 $(\ref{etadd16})$ is well-defined. 

%Let us assume that the manifold $N^{n-2k}_{b \times \bb}$ is connected and is equipped with the marked point.
Let us prove that the image of the fundamental class of the target manifold by the mapping
$(\ref{Neta7})$ determines the element
\begin{eqnarray}\label{etadlok}
\eta^{loc}_{Ab,\ast}([\eta_N^{7k}]^{op}) \in D_{m_{\sigma}}(\I_{b \times \bb};\Z/2[\Z/2]),
\end{eqnarray}
which is mapped into the element
\begin{eqnarray}\label{etadd}
\eta_{Ab,\ast}([\eta_N^{7k}]^{op}) \in H_{m_{\sigma}}(K(\I_{b \times \bb},1)).
\end{eqnarray}
by means of the homomorphism $(\ref{deltaD2})$.

Consider the submanifold
$(\ref{NN})$ and consider the decomposition of this manifold into connected components:
\begin{eqnarray}\label{NNlok}
N_{\eta^{7k}}^{n-16k} \equiv \cup_i N_{i,\eta^{7k}}^{n-16k} \subset N^{n-2k}_{b \times \bb}.
\end{eqnarray}
On each connected component of the manifold
$(\ref{NNlok})$  let us take a marked point  $pt_i \in N_{i,\eta^{7k}}^{n-16k}$.
Take a path
 $\rho_i$ on the manifold  $N^{n-2k}_{b \times \bb}$ from the point  $pt_i$ to the point $pt$.
For an arbitrary  $i$ the isomorphism of the fiber $\kappa_b \oplus \kappa_{\bb}$ over the point $pt$ and the fiber over the point $pt_i$ is well-defined.

Therefore the following mapping
\begin{eqnarray}\label{Neta7lok}
\eta_{N_{i,\eta^{7k}}}(\rho_i): N_{i,\eta^{7k}}^{n-16k} \to K(\I_b \times \II_b,1).
\end{eqnarray}
is well-defined. The immersed manifold
 $N_{\eta^{7k}}^{n-16k}$ is equipped with a framing $\Psi_{Ab}$. This framing over each component
$N_{i,\eta^{7k}}^{n-16k}$ of  $N_{\eta^{7k}}^{n-16k}$ is totally determined by the choice of the prescribed coordinate system in the fiber over the marked point $pt_i \in N_{i,\eta^{7k}}^{n-16k}$.

At the opposite site, the choice of the coordinate system in the fiber over 
 $pt_i$ is changed by the transformation on the element $s(i) \in \Z/2$ from a residue class of the subgroup 
$\I_{b \times \bb} \subset \Z/2^{[2]}$.  For the non-trivial residue class the transformation
is given by the element   $ab \in \D \setminus \I_b$. Denote  the element
$\eta_{Ab,i,\ast}(\rho_i)(N_{i,\eta^{7k}}^{n-16k}) \in H_{m_{\sigma}}(\I_{b \times \bb})$, 
by $x_i(\rho_i)$. Let us define the element
\begin{eqnarray}\label{Xi}
 X_i \in D_{m_{\sigma}}(\I_{b \times \bb};\Z/2[\Z/2]) 
\end{eqnarray} 
is equal to $x_i(\rho_i) + 0t$, if the framing  $\Psi_{Ab}$ over the manifold
$N_{i,\eta^{7k}}^{n-16k}$ is agree with the framing, which is obtained by means of the parallel translation
of the framing
$\Psi_{b \times \bb}$ along the path $\rho_i$, and define
$Xi_i=0+ \chi^{[2]}_{\ast}(x_i)t$, if the considered framings are not agree. 
The element 
 $(\ref{etadlok})$ is well-defined.
By the construction, the element
$(\ref{etadlok})$ does not depended on a choice of a path  $\rho_i$ and a reduction $\Psi_{Ab}$.

The following lemma is proved by a straightforward calculation.

\begin{lemma}\label{dd}

Let us assume that the mapping $(\ref{redu})$ is well defined as in Definition $\ref{strucIdId}$. 
Then the following two properties are satisfied.

--1. The element $(\ref{etadlok})$, which is constructed by means of the mapping  $(\ref{Neta7})$, 
has the image with respect to the homomorphism
 $(\ref{deltaD2})$, such that the decomposition of this image over the standard base of the group
$H_{m_{\sigma}}(K(\I_b \times \II_b,1))$ contains not more then one non-trivial element, which is determined
by the coefficient of the monomial $t_{b,i} \otimes t_{\bb,i}$, see. $(\ref{tdtdd})$, $i=\frac{m_{\sigma}}{2}=\frac{n-16k}{2}$. 
This coefficient coincides with the characteristic number
$(\ref{44})$ for the $\Z/2^{[2]}$--framed immersion $(g,\Psi,\eta_N)$.

--2. The element $(\ref{etadd16})$ belongs to the subgroup  $(\ref{DZ2lok})$, $i=n-16k$. 
\end{lemma}

\subsubsection*{Proof of Lemma  $\ref{dd}$}
Let us prove the statement 1 of the lemma in the case
$m_{\sigma}=14$. The general case is analogous.
Consider the manifold
 $N_{\eta^{7k}}^{14}$, in the proof we denote this manifold by $N^{14}$ for short.
 The manifold $N^{14}$ is equipped with the mapping
  $(\ref{Neta7})$, in the proof we denote this mapping by
$\eta: N^{14} \to K(\I_{b \times \bb},1)$ for short.
Consider all characteristic numbers with $\Z/2$--coefficients for the mapping $\eta$, which are induced from
the universal characteristic classes $(\ref{kappab})$ by the mapping $(\ref{Neta7})$.
We will define the induced classes the same.

Because $N^{14}$ is oriented,  non-trivial numbers, possibly,  are following:
$\kappa_b\kappa^{13}_{\bb}$, $\kappa_b^3\kappa^{11}_{\bb}$, $\kappa_b^5\kappa^{9}_{\bb}$,  $\kappa_b^7\kappa^{7}_{\bb}$, $\kappa_b^9\kappa^{5}_{\bb}$,
$\kappa_b^3\kappa^{11}_{\bb}$, $\kappa_b\kappa^{13}_{\bb}$.

Let us prove that the characteristic number  $\kappa_b\kappa^{13}_{\bb}$
is trivial. Consider a submanifold $K^3 \subset N^{14}$,
which is dual to the characteristic class  $\kappa_b\kappa_{\bb}^{10}$.
The normal bundle  $\nu_N$ of the manifold $N^{14}$ is isomorphic to the bundle $k\kappa_b \oplus k\kappa_{\bb}$,
where $k \equiv 0 \pmod{8}$. The restriction of the bundle $\nu_N$ on the submanifold $K^3$ is trivial.
Therefore, the normal bundle  $\nu_K$ of the manifold $K^3$ is stably equivalent to the bundle
$\kappa_b \oplus 2\kappa_{\bb}$. The characteristic class  $w_2(K^3)$ is trivial, in particular,
$\langle \kappa^3_{\bb};[K^3]\rangle=0$. We have  $\langle \kappa^3_{\bb};[K^3]\rangle=
\langle \kappa_b\kappa^{13}_{\bb};[N^{14}]\rangle$. This proves that the characteristic number
$\kappa_b\kappa^{13}_{\dd}$ is trivial. Analogically, the following characteristic numbers
$\kappa_b^5\kappa^{9}_{\bb}$, $\kappa_b^9\kappa^{5}_{\bb}$, $\kappa_b\kappa^{13}_{\bb}$ are trivial.

Let us prove that the characteristic number
 $\kappa_b^3\kappa^{11}_{\bb}$ is trivial. Take a submanifold
 $K^7 \subset N^{14}$, which is dual to the characteristic class
$\kappa_b^2\kappa_{\bb}^5$. the normal bundle  $\nu_K$ of the manifold $K^6$ is stably equivalent to the bundle
$2\kappa_b \oplus 5\kappa_{\bb}$. Because $w_5(K^7)=0$, the characteristic class
 $\kappa_{\bb}^{5}$ is trivial. In particular,  $\langle \kappa_b\kappa_{\bb}^6;[K^7]\rangle=0$.
The following equation is satisfied:
$\langle \kappa_b\kappa_{\bb}^6;[K^7]\rangle = \langle \kappa_b^3\kappa_{\bb}^{11};[N^{14}]\rangle$.
This proves that the characteristic number
$\kappa_b^3\kappa^{11}_{\bb}$ is trivial. Analogously, the characteristic number
$\kappa_b^{11}\kappa^{3}_{\bb}$ is trivial.

Obviously, the characteristic number
 $\langle \kappa_b^7\kappa_{\bb}^7;[N^{14}] \rangle$ coincides with the
 characteristic number
 $(\ref{44})$. The statement 1 of the lemma is proved.
 
Let us prove the statement 2. Because the manifold  $N_{b \times \bb}^{n-16k}$ in the case $\sigma \ge 5$
is oriented. Consider the decomposition of the element  $(\ref{etadd16})$ over the base of the group $H_{n-16k}(K(\I_{b \times \bb} \int_{\chi^{[2]}} \Z,1);\Z/2[\Z/2])$, using Lemma $\ref{obstr}$. By Lemma
$\ref{imX}$, which is proved in Section  6, the element $R$ is not involved to the considered expansion.
Another case, the mapping
 $\eta_{b \times \bb}$ is not satisfied Condition 2 in Definition  $\ref{strucIdId}$. 
 The statement 2 is proved.
 Lemma $\ref{dd}$ is proved.
\[  \]

\begin{definition}\label{5}
Let $[(f,\Xi,\kappa_M)] \in Imm^{sf}(n-k,k)$,
 $f: M^{n-k} \looparrowright \R^n$, $\kappa_M \in H^1(M^{n-k};\Z/2)$
 be skew-framed by
$\Xi$. We say that the pair $(M^{n-k},\kappa_M)$ admits a
compression of order $q$ if the mapping  $\kappa_M : M^{n-k} \to
\RP^{\infty}$ can be represented as a composition $\kappa = I
\circ \kappa_M' : M^{n-k} \to \RP^{n-k-q-1} \subset \RP^{\infty}$,
where $I$ denotes the inclusion. We say that the element $[(f,
\Xi, \kappa_M)]$ admits a compression of order $q$ if this
cobordism class contains a triple $(f', \Xi', \kappa_{M'})$,  so
that the pair $(M'^{n-k},\kappa_{M'})$ admits a compression of
order $q$.
\end{definition}

\begin{theorem}$\label{th6}$
Let $m_{\sigma}=2^{\sigma}-2$, $\sigma \ge 5$, $n \ge 4m_{\sigma}+6$. Assume that the
element $\alpha \in Imm^{sf}(n-\frac{n-n_\sigma}{16},\frac{n-m_\sigma}{16})$
admits a compression of order $q=\frac{m_\sigma}{2}-1$ (in particular, in the case  $\sigma=5$, $q=16$). Then the
element $\alpha$ admits an $\I_{b \times \bb}$--structure.
\end{theorem}

Let
\begin{eqnarray}\label{d2}
d^{(2)}: \RP^{n-k} \times \RP^{n-k} \to \R^n \times \R^n
\end{eqnarray}
be an arbitrary $T_{\RP^{n-k}}$, $T_{\R^n}$--equivariant mapping, which is transversal
along the diagonal
$\R^n_{diag} \subset \R^n \times \R^n$ like in the formula (41)[A1]. (We will use only the case 
$k'=k+q+1$,  as in denotations of Theorem $\ref{th6}$.)

Denote
$(d^{(2)})^{-1}(\R^n_{diag})/T_{\RP^{n-k}}$ by  $\N=\N(d^{(2)})$. Let us call this polyhedron 
the polyhedron of (formal) intersection of the equivariant mapping
 $d^{(2)}$.

The polyhedron  $\N$, generally speaking, contains a nonempty boundary  $\partial \N$  
this boundary is represented by critical points of the mapping $d^{(2)}$.
Denote by
 $\N_{\circ}$ the open polyhedron 
$\N \setminus \partial \N$,
denote by $U(\partial \N)_{\circ}$ the deleted regular neighborhood of the boundary  $\partial \N$.

Let us assume the following equivariant mapping
 $(\ref{d2})$ and the following generic ${\rm{PL}}$--mapping  
\begin{eqnarray}\label{dmain}
 d: \RP^{n-k} \to \R^n
\end{eqnarray}  
 are given.

 \begin{definition}\label{hol} 
Let us say that a  formal (equivariant) mapping
$d^{(2)}$, given by  $(\ref{d2})$, is  holonomic if this mapping is the extension of a mapping 
\begin{eqnarray}\label{dmain}
 d: \RP^{n-k} \to \R^n.
\end{eqnarray}  
\end{definition}

\begin{definition}\label{abstru}
Let a  formal (equivariant) mapping
 $(\ref{d2})$ be given. 
 Let us say that 
$d^{(2)}$ admits an abelian structure, if the following condition is satisfied.

--  On the open polyhedron  $\N_{\circ}$
of formal self-intersection points of  $d^{(2)}$  the following mapping
\begin{eqnarray}\label{abels}
\eta_{b \times \bb\circ}: (\N_{\circ},U(\partial \N)_{\circ}) \to (K(\I_{b \times \bb} \int_{\chi^{[2]}} \Z,1), K(\I_{b \times \bb},1)),
\end{eqnarray}
is well-defined. This mapping determines a reduction of the structure mapping
 $$\eta_{\circ}: (\N_{\circ},U(\partial \N)_{\circ}) \to (K(\Z/2^{[2]},1), K(\I_{b \times \bb},1)),$$
which satisfies the boundary conditions (about the notion  "`structure mapping"' see 
 [A1, formula (48)]).
\end{definition}

The following lemma is proved in [A3].

\begin{lemma}\label{7}
Assume that the following dimensional restriction
\begin{eqnarray}\label{dimdim}
n-k' \equiv -1 \pmod{4}, \quad k'
\ge 7, \quad n \equiv 0 \pmod{2}
\end{eqnarray}
is satisfied.
Then there exists a formal (equivariant) mapping $d^{(2)}$, which admits an abelian structure in the sense of Definition $\ref{abstru}$.
\end{lemma}

\subsubsection*{Proof of Theorem $\ref{th6}$}
Put $k=\frac{n-n_\sigma}{16}$. Let the element $\alpha$ be represented by a skew-framed immersion
 $(f,\Xi,\kappa_M)$, $f: M^{n-k}
\looparrowright \R^n$.  By assumption there exists a compression
$\kappa_M': M^{n-k} \to  \RP^{n-k-q-1}$ such that the composition
$M^{n-k} \to \RP^{n-k-q-1} \subset K(\I_d,1)$ coincides with the
mapping $\kappa_M: M^{n-k} \to K(\I_d,1)$.

Let us define $k'$ as the maximal  integer, which is less or equal to  $k+q+1$, and such that  
$n-k' \equiv -1\pmod{4}$. In the case $\sigma=5$ we get $k'=k+q+1$, in the case $\sigma \ge 6$
we get $k'=k+q-1$.
 Because $q=\frac{m_{\sigma}}{2} \ge 14$, for
$k' \ge  7$ the both dimensional restrictions  $(\ref{dimdim})$ are satisfied.

By Lemma
$\ref{7}$ there a formal (equivariant) mapping $d^{(2)}$,  which admits an abelian structure
 in the sense of Definition  $\ref{abstru}$.

Let us construct a skew-framed immersion
$(f_1,\Xi_1,\kappa_1)$, for which the immersed $\Z/2^{[2]}$--framed  self-intersection manifold 
contains a closed component
 $N_{b \times \bb}^{n-2k}$, as it is required in the formula 
$(\ref{abelcomp})$.

 Define the immersion $f_1: M^{n-k} \looparrowright
\R^n$ using [Corollary 31, A1]  as a result of  a $C^0$--small regular deformation of the composition
$d \circ \kappa: M^{n-k} \to \RP^{n-k-q-1} \subset \RP^{n-k} \to \R^n$ in the prescribed regular homotopy class of the immersion $f: M^{n-k} \looparrowright \R^n$. Define a caliber of the deformation  $d \circ \kappa \mapsto f_a$ 
much less then the radius of the regular neighborhood  of the polyhedron of self-intersection points
of the mapping
$d$.

Denote by $N^{n-2k}_{b \times \bb}$
the self-intersection manifold of the immersion  $f_1$.
The following decomposition of the manifold  into the union of two 
manifolds with boundaries along the common boundary is well defined: 
\begin{eqnarray}\label{Nregb}
 N^{n-2k}_{b \times \bb} = N^{n-2k}_{b \times \bb,\N(d^{(2)})} \cup_{\partial}
N^{n-2k}_{reg}.
\end{eqnarray}
In this formula $N^{n-2k}_{b \times \bb,\N(d^{(2)})}$
is a manifold with boundary, which is immersed into a regular (immersed)
 neighborhood 
$U_{\N(d)}$ of the polyhedron with the boundary  $\N(d)$ of self-intersection points of the 
mapping  $d$. The manifold  $N^{n-2k}_{reg}$ with boundary is immersed
into a regular immersed neighborhood (denote this immersed neighborhood by $U_{reg}$) of self-intersection points of the mapping $d$ outside of critical points.
The common boundary of the manifolds  $N^{n-2k}_{b \times \bb,\N(d^{(2)})}$, $N^{n-2k}_{reg}$
is a closed manifold of the dimension  $n-2k-1$, this manifold is immersed into the boundary
 $\partial(U_{reg})$ of the immersed neighborhood $U_{reg}$.

Denote the $\Z/2^{[2]}$--framed immersion of the self-intersection  $(\ref{Nregb})$ by $g_{b \times \bb}$.
Let us prove that the $\Z/2^{[2]}$--framing
  over $(\ref{Nregb})$ is reduced to a framing
$(\eta_{b \times \bb},\Psi_{b \times \bb})$ with the structure group
 $\I_{b \times \bb} \int_{\chi^{[2]}} \Z$.

Define a mapping
 $N^{n-2k}_{b \times \bb}
\stackrel{\eta_{b \times \bb}}{\longrightarrow} K(\I_{b \times \bb} \int_{\chi^{[2]}} \Z,1)$,
which is determined the required reduction of the characteristic mapping of the framing. 
Define on the submanifold $N^{n-2k}_{reg} \subset N^{n-2k}_{b \times \bb}$ the mapping
$\kappa_{b,N_{reg}}: N^{n-2k}_{reg} \to K(\I_b,1)$ by the composition of the projection $N^{n-2k}_{reg} \to \RP^{n-k'}$ and the inclusion $\RP^{n-k'}
\subset \RP^{\infty} = K(\I_b,1)$.

The  cohomology class  $\kappa_{\bb,N_{reg}} \in
H^1(N^{n-2k}_{reg};\Z/2)$ is defined as the orientation class
of the line bundle $p \otimes \kappa_{b,N_{reg}}$, this bundle is the tensor product of the line bundle, associated with
the canonical double covering
 $p: \bar N^{n-2k}_{reg} \to N^{n-2k}_{reg}$ and the line bundle with the characteristic class  $\kappa_{d,N_{reg}}$, namely
$$ \kappa_{\bb,N_{reg}} = p \otimes \kappa_{b,N_{reg}}. $$

 The pair of the cohomology classes
 $\kappa_{\bb,N_{reg}}$, $\kappa_{b,N_{reg}}$ determines the required mapping
 \begin{eqnarray}\label{etabbbNreg}
\eta_{b \times \bb,N_{reg}}:
N^{n-2k}_{reg} \to K(\I_{b \times \bb},1) \subset K(\I_{b \times \bb} \int_{\chi^{[2]}} \Z,1).
\end{eqnarray}
The denotations are agree with the definition of universal characteristic classes 
 $(\ref{kappab})$.

Let us define on the first component
$N^{n-2k}_{b \times \bb,\N(d^{(2)})}$ in the formula $(\ref{Nregb})$ the mapping
\begin{eqnarray}\label{etadddNd2}
\eta_{b \times \bb,\N(d^{(2)})}:
N^{n-2k}_{b \times \bb,\N(d^{(2)})} \to K(\I_{b \times \bb} \int_{\chi^{[2]}} \Z),1)
\end{eqnarray}
as the composition of the projection
$N^{n-2k}_{b \times \bb,\N(d^{(2)})} \to \N(d)$ and the mapping $\eta_{b \times \bb}: \N(d^{(2)}) \to
K(\I_{b \times \bb} \int_{\chi^{[2]}} \Z),1)$, which is constructed in Lemma $\ref{7}$.

Restrictions of the mappings $\eta_{b \times \bb,\N(d^{(2)})}$, $\eta_{b \times \bb,N_{reg}}$
to the common boundaries $\partial N^{n-2k}_{b \times \bb,\N(d^{(2)})}$,  $\partial N^{n-2k}_{reg}$
are homotopic, because the mapping $\eta_{\N(d^{(2)})}: \N(d^{(2)}) \to K(\I_{b \times \bb} \int_{\chi^{[2]}} \Z),1)$
satisfies the boundary conditions on $\partial \N(d^{(2)})$. Therefore the mapping
\begin{eqnarray}\label{etadtd}
\eta_{b \times \bb}: N^{n-2k}_{b \times \bb}
\to K(\I_{b \times \bb} \int_{\chi^{[2]}} \Z),1)
\end{eqnarray}
 is well defined as the result of the gluing of the two mappings
 $\eta_{b \times \bb,\N(d^{(2)})}$ и $\eta_{b \times \bb,\N_{reg}}$.
 The mapping
 $(\ref{etadtd})$ determines a reduction of the characteristic mapping of the $\Z/2^{[2]}$--framing of the immersion $g_{b \times \bb}$.

Let us comes back to the immersion 
$f_{a}$ of the manifold $(\ref{Nreg})$. Consider the submanifold
 $M^{n-8k}_{2,a} \subset M^{n-2k}$, which is defined as the inverse image  $M^{n-8k}_{2,a}=
(\kappa_M')^{-1}(\RP^{n-8k-q-1})$ of the projective subspace 
$\RP^{n-8k-q-1} \subset \RP^{n-k-q-1}$ of the codimension $7k$.
Define the immersion
$f_{2,a}: M^{n-8k}_{2,a} \looparrowright \R^n$ by the restriction
of the immersion $f_a$ on the submanifold $M^{n-8k}_{2,a}$.

By construction of
$f_a$ the image of  $f_{2,a}$ is inside a regular neighborhood of the image of the submanifold $\RP^{n-8k-q-1} \subset \RP^{n-7k-k'}$ by the mapping $d$. Because
$n-8k-q-1 = n- \frac{n}{2} + \frac{m_{\sigma}}{2} - \frac{m_{\sigma}}{2} -1= \frac{n}{2}-1$,
we get $d(\RP^{n-8k-q-1})$ is an embedded submanifold.

Denote by
$N_{2,a}^{n-16k}$ the self-intersection manifold of the immersion   $f_{2,a}$. We have $\dim(N_{2,a})=n-16k=m_{\sigma}$.
$n-16k=m_{\sigma}$. The following inclusion is well-defined:
\begin{eqnarray}\label{incl}
N_{2,a}^{n-16k} \subset N_{a,reg}^{n-2k},
\end{eqnarray}
where the manifold
 $N_{a,reg}^{n-2k}$ is determined by the formula  $(\ref{Nreg})$. 

In particular, a reduction of the classifying mapping
$\eta_{2,a}: N^{n-16k}_{2,a} \to K(\Z/2^{[2]},1)$ into the subspace $K(\I_{b \times \bb},1) \subset K(\Z/2^{[2]},1)$
is well-defined:
\begin{eqnarray}\label{red2a}
\eta_{Ab,a}: N^{n-16k}_{2,a} \to K(\I_{b \times \bb},1). 
\end{eqnarray}

Therefore in the case 
$\sigma \ge 5$ without loss of a generality we may assume that
\begin{eqnarray}\label{emptyset}
N^{n-16k}_{2} \cap N^{n-2k}_{b \times \bb,\N(d)} = \emptyset. 
\end{eqnarray}

The triple $(f_1,\Xi_1,\kappa_1)$ defines the required skew-framed immersion in the cobordism class
$x$. Let us define the second component in the formula 
 $(\ref{abelcomp})$ as the empty component.  The reduction mapping  $(\ref{red2a})$ (recall that
$N^{n-16k}_{2,reg}$ is re-denoted by  $N_{\eta^{7k}}^{n-16k}$) coincides with the required reduction,
which is given by the formula $(\ref{Neta7})$/ This proves Property 1 from Definition  $\ref{strucIdId}$. 
Property 2 is analogous. 

Арф-инвариант, определеный при помощи $\Z/2^{[2]}$--оснащенного многообразия
$(\ref{Nregb})$ совпадает с Арф-инвариантом исходного  элементa $x$.
 The Arf-invariant, which is defined by the   $\Z/2^{[2]}$--framed immersion of the manifold $(\ref{Nregb})$
coincides with
the Arf-invariant of the  element  $x$.

Let us present a sketch of the proof without the assumption that the mapping 
 $d^{(2)}$ is holonomic. Analogously to [Theorem  23, A1], let us generalize the construction above and let us define 
a $\Z/2^{[2]}$--framed immersion  $(g_{b \times \bb},\Psi,\eta)$ of the manifold  $(\ref{Nregb})$, which satisfies Conditions 1 and 2,
but, probably, which is not an immersion of a self-intersection of a framed immersion 
Evidently, there exists a skew-framed immersion
 $(f,\Xi,\kappa)$, which represents the element $x$, for which an arbitrary  $\Z/2^{[2]}$--framed immersion, in particular, the $\Z/2^{[2]}$--framed immersion
$(g_{b \times \bb},\Psi,\eta)$, is a closed (probabely, non-connected) component of the self-intersection manifold of this skew-framed immersion $(f,\Xi,\kappa)$.  The
second component in the formula  $(\ref{abelcomp})$ is defined as the last component of the self-intersection manifold of the skew-framed immersion $(f,\Xi,\kappa)$. Properties 1 and 2 are evident. 
Theorem $\ref{th6}$ is proved.

\section{$\H_{b \times \bb}$--structure  of a $\Z/2^{[2]}$--framed immersion and $\J_a \times \JJ_a$-structure
(bicyclic structure) of a  $\Z/2^{[3]}$--framed immersion}

We define the group
 $\I_a$ as the cyclic subgroup of order  $4$ in the dihedral group
 $\I_a \subset \Z/2^{[2]}$, see section 2 of [A1]. We shall now define an analogous subgroup
\begin{eqnarray} \label{iaa}
i_{\J_a \times \JJ_a}: \J_a \times \JJ_a \subset
 \Z/2^{[4]},
\end{eqnarray}
which is isomorphic to the direct product of the two cyclic groups of the order $4$.

Recall that the group $\Z/2^{[4]}$ is defined in terms of a basis
$(\e_1, \dots, \e_8)$ of the Euclidean space $\R^8$.
We denote the generators of the factors of the group
 $\J_a \times
\JJ_a$ by $a$ and $\aa$
respectively. We shall describe transformations in
$\Z/2^{[4]}$, corresponding to
each generator. We introduce a new basis $\{\f_1, \dots, \f_8\}$, by the formulas
$\f_{2i-1}=\frac{\e_{2i-1}+ \e_{2i}}{\sqrt{2}}$, $\f_{2i}= \frac{\e_{2i-1}- \e_{2i}}{\sqrt{2}}$,
$i=1, \dots, 4$.

We shall show that the group of transformations
 $\J_a \times \JJ_a$ has invariant orthogonal
 $(2,2,2,2)$--dimensional subspaces, which we denote by $\R^2_{a,+}$, $\R^2_{a,-}$,
 $\R^2_{\aa,+}$, $\R^2_{\aa,-}$.

The subspace
 $\R^2_{a,+}=Lin(\f_1+\f_5,\f_3+\f_7)$ is generated by the pair of vectors
$(\f_1+\f_5,\f_3+\f_7)$. The subspace $\R^2_{a,-}=Lin(\f_1-\f_5,\f_3-\f_7)$
is generated by the pair of vectors $(\f_1-\f_5,\f_3-\f_7)$. The
subspace
$\R^2_{\aa,+}=Lin(\f_2+\f_4,\f_6+\f_8)$ is generated by the pair of vectors  $(\f_2+\f_4,\f_6+\f_8)$.
The subspace $\R^2_{\aa,-}=Lin(\f_2-\f_4,\f_6-\f_8)$ is generated by the pair of vectors
$(\f_2-\f_4,\f_6-\f_8)$.

It is convenient to pass to a new basis
\begin{eqnarray}\label{h}
\frac{\f_1+\f_5}{\sqrt{2}}=\h_{1,+},
\frac{\f_1-\f_5}{\sqrt{2}}=\h_{1,-},
\frac{\f_3+\f_7}{\sqrt{2}}=\h_{2,+},
\frac{\f_3-\f_7}{\sqrt{2}}=\h_{2,-},
\end{eqnarray}
\begin{eqnarray}\label{hh}
\frac{\f_2+\f_4}{\sqrt{2}}=\hh_{1,+},
\frac{\f_2-\f_4}{\sqrt{2}}=\hh_{1,-},
\frac{\f_6+\f_8}{\sqrt{2}}=\hh_{2,+},
\frac{\f_6-\f_8}{\sqrt{2}}=\hh_{2,-}.
\end{eqnarray}
  
  The pairs of vectors $(\h_{1,+},\h_{2,+})$, $(\h_{1,-}, \h_{2,-})$ are bases for the subspaces $\R^2_{a,+}=Lin(\h_{1,+},\h_{2,+})$,  $\R^2_{a,-}=Lin(\h_{1,-}, \h_{2,-})$ respectively. In addition, the pairs of vectors
$(\hh_{1,+},\hh_{2,+})$,
$(\hh_{1,-}, \hh_{2,-})$ are
bases for the subspaces
$\R^2_{\aa,+}=Lin(\hh_{1,+},\hh_{2,+})$,  $\R^2_{\aa,-}=Lin(\hh_{1,-}, \hh_{2,-})$ respectively.

The generator $a$ of order $4$ is represented by a rotation
through angle $\frac{\pi}{2}$ in each of the planes $\R^2_{a,+}$,
$\R^2_{a,-}$ and by the central symmetry in the plane
$\R^2_{\aa,-}$, (Evidently, the image of the generator $a$
commutes with the presentation of the generator $\aa$ see below).
The generator $\aa$ is represented by a rotation through angle
$\frac{\pi}{2}$ in the planes $\R^2_{\aa,+}$ $\R^2_{\aa,-}$ and by
the central symmetry in the plane $\R^2_{a,-}$. So the subgroup
($\ref{iaa}$) is well defined.

Define the subgroup 
$i_{\H_{b \times \bb}, \J_a \times \JJ_a}: \H_{b \times \bb} \subset \J_a \times \JJ_a$, as
the direct product of the diagonal subgroup
$\I_a \subset \J_a \times \JJ_a$ with the elementary subgroup $\JJ_d \subset \JJ_a$ of the
second factor. 
The subgroup $\H_{b \times \bb} \subset \J_a \times \JJ_a$ coincides with the preimage of the subgroup
$\Z/2 \subset \Z/4$ with respect to the homomorphism
\begin{eqnarray}\label{omega4}
\omega^{[4]}: \I_a \times \II_a \to \Z/4,
\end{eqnarray}
defined by the formula $(x \times y) \mapsto xy$.

\subsubsection*{Remark}
The group $\H_{b \times \bb}$ is isomorphic to the group $\E_{b \times \bb}$, which is defined
in [A1, Section 2]. The corresponding subgroups in $\Z/2^{[3]}$ are distinguished. 
\[  \]

Define the subgroup
$i_{\I_{b \times \bb}, \H_{b \times \bb}}: \I_{b \times \bb} \subset \H_{b \times \bb}$,
as the kernel of the homomorphism 
\begin{eqnarray}\label{omega3}
\omega^{[3]}: \H_{b \times \bb} \to \Z/2,
\end{eqnarray}
which is determined by the formula
$(x \times y) \mapsto x$ in terms of the generators of the ambiance group. 

Consider the diagonal subgroup $\Z/2^{[3]} \subset \Z/2^{[4]}$, this subgroup is generated by transformations
in the direct sum of the subspaces
 $diag(\R^2_{a,+},\R^2_{\aa,+})$,
$diag(\R^2_{a,-},\R^2_{\aa,-})$.
This subgroup is the transformation group of the unit vectors, which are collinear to  the vectors
$\h_{1,+} + \hh_{1,+},\h_{2,+} + \hh_{2,+},\h_{1,-} + \hh_{2,-}, \h_{2,-} + \hh_{2,-}$.
This collection of vectors gives the standard base in the subspace
 $diag(\R^2_{a,+},\R^2_{\aa,+}) \oplus diag(\R^2_{a,-},\R^2_{\aa,-})$.
 The complement to this subspace is given by the formula
$antidiag(\R^2_{a,+},\R^2_{\aa,+}) \oplus antidiag(\R^2_{a,-},\R^2_{\aa,-})$.
In this complement the standard base is given analogously.

In the notations above the diagonal subgroup
 $\Z/2^{[2]} \subset \Z/2^{[3]}$ is the transformation group of
 the unit vectors, which are collinear to the vectors
 $\h_{1,+} + \hh_{1,+} + \h_{2,+} + \hh_{2,+}$, $\h_{1,-} + \hh_{2,-} + \h_{2,-} + \hh_{2,-}$.
 It is easy to check that the first vector is collinear to the vector
 $\e_1 + \e_3 + \e_5 + \e_7$ and the second vector is collinear to the vector
$\e_2 + \e_4 + \e_5 + \e_6$.

There is an inclusion $i_{\H_{b \times \bb}}: \H_{b \times \bb} \subset
 \Z/2^{[3]}$, 
 which is compatible
with the inclusion ($\ref{iaa}$).
 Moreover, the following diagram is commutative:
 
\begin{eqnarray}
\label{a,aa}
\begin{array}{ccc}
\qquad \I_{b \times \bb} & \stackrel {i_{b \times \bb}} {\longrightarrow}& \qquad \Z/2^{[2]} \\
i_{b \times \bb,\H_{b \times \bb}} \downarrow &  & i^{[3]} \downarrow \\
 \qquad \H_{b \times \bb} &  \stackrel {i_{\H_{b \times \bb}}}
{\longrightarrow}& \qquad \Z/2^{[3] } \\
i_{\H_{b \times \bb},\J_a \times \JJ_a} \downarrow &  &  i^{[4]} \downarrow \\
\qquad \J_a \times \JJ_a & \stackrel{i_{\J_a \times \JJ_a}}{\longrightarrow} &  \qquad \Z/2^{[4]}.\\
\end{array}
\end{eqnarray}

Let us define automorphisms
\begin{eqnarray}\label{chi3A}
\chi^{[3]}: \H_{b \times \bb} \to \H_{b \times \bb},
\end{eqnarray}
\begin{eqnarray}\label{chi3A}
\chi^{[4]}: \J_a \times \JJ_a \to \J_a \times \JJ_a.
\end{eqnarray}
of order 2.

Let us also define automorphisms
\begin{eqnarray}\label{chi3ZA}
\chi^{[3]}: \Z/2^{[3]} \to \Z/2^{[3]},
\end{eqnarray}
\begin{eqnarray}\label{chi4ZA}
\chi^{[4]}: \Z/2^{[4]} \to \Z/2^{[4]},
\end{eqnarray}
of order 2,
which are denoted the same.
The automorphism $(\ref{chi3ZA})$
is defined by the permutation of the corresponding basis vectors
of the subspace $diag(\R^2_{a,+},\R^2_{\aa,+}) \oplus diag(\R^2_{a,-},\R^2_{\aa,-})$ with the indexes $a$ and $\aa$.

Define the automorphism  $(\ref{chi3A})$  such that its restriction to the diagonal subgroup 
$diag(\J_a,\JJ_a)=\I_a \subset \H_{b \times \bb}$
 coincides with the identity, and the restriction of this automorphism to the subgroup  
$\I_{b \times \bb} \subset \H_{b \times \bb}$ 
coincides with the automorphism
 $\chi^{[2]}$. Evidently, the automorphism  $(\ref{chi3A})$ is uniquely well defined.

The automorphism  $(\ref{chi4ZA})$
 is defined by means of the standard bases of the subspaces $diag(\R^2_{a,+},\R^2_{\aa,+}) \oplus diag(\R^2_{a,-},\R^2_{\aa,-})$, $antidiag(\R^2_{a,+},\R^2_{\aa,+}) \oplus antidiag(\R^2_{a,-},\R^2_{\aa,-})$
 by the permutation of the basis vectors with the indexes $a$ and $\aa$. 
 %The automorphism $(\ref{chi4A})$ define by the formula  $(x \times y) \mapsto (y \times x)$.
% Evidently, the automorphisms $(\ref{chi3A})$-$(\ref{chi4ZA})$ are uniquely well defined.

Recall that the group $\I_{b \times \bb} \int_{\chi^{[2]}} \Z$ is defined by the formula $(\ref{chi2Z})$.
Define the analogous subgroups
\begin{eqnarray}\label{chi3Z}
\H_{b \times \bb} \int_{\chi^{[3]}} \Z,
\end{eqnarray}
\begin{eqnarray}\label{chi4Z}
(\J_a \times \JJ_a) \int_{\chi^{[4]}} \Z,
\end{eqnarray}
as the semi-direct product of the corresponding groups with involutions with the group  $\Z$.

Namely, define the group  $(\ref{chi3Z})$ as the factorgroup of the group
 $\H_{b \times \bb} \ast \Z$ by the relations   $zxz^{-1}=\chi^{[3]}(x)$, where $z \in \Z$ is the generator,
$x \in \H_{b \times \bb}$ is an arbitrary element. The classifying space  $K(\H_{b \times \bb} \int_{\chi^{[3]}} \Z,1)$
is the semi-direct product of the standard circle  $S^1$ and the space $K(\H_{b \times \bb},1)$. 
The shift map  $K(\H_{b \times \bb},1) \to K(\H_{b \times \bb},1)$ of the cyclic cover over
 $K(\H_{b \times \bb} \int_{\chi^{[3]}} \Z,1)$ is induced by the automorphism  $\chi^{[3]}$.
The definition of the group $(\ref{chi4Z})$ is analogous.

Let us consider homology groups $H_{i}(K(\H_{b \times \bb} \int_{\chi^{[3]}} \Z,1))$,
$H_{i}(K((\J_a \times \JJ_a) \int_{\chi^{[4]}} \Z,1))$.
 In particular, for an odd $\ast=i$, the second group contains the elements, which is represented by the
 fundamental classes of the submanifolds
  $S^{i}/\i \times pt \subset S^i/\i \times S^i/\i \subset K(\J_a,1) \times K(\JJ_a,1)$, $pt \times S^{i}/\i  \subset S^i/\i \times S^i/\i \subset K(\J_a,1) \times K(\JJ_a,1)$, denote the elements by 
\begin{eqnarray}\label{tataa}
t_{a,i},  \in H_{i}(K(\J_a,1));\quad t_{\aa,i} \in H_{i}(K(\JJ_a,1)).
\end{eqnarray}

The homology groups with local coefficients
$\Z[\Z/2]$ are defined analogously to
 $(\ref{Z2lok})$, $(\ref{Zlok})$, $(\ref{DZlok})$. For example: 
\begin{eqnarray}\label{DZlokaa}
D_{i}(\J_a \times \JJ_a;\Z[\Z/2]) \subset H_{i}(K((\J_a \times \JJ_a) \int_{\chi^{[4]}} \Z,1);\Z[\Z/2]),
\end{eqnarray}
\begin{eqnarray}\label{DZlokad}
D_{i}(\H_{b \times \bb};\Z[\Z/2]) \subset H_{i}(K(\H_{b \times \bb} \int_{\chi^{[3]}} \Z,1);\Z[\Z/2]).
\end{eqnarray}
Analogous groups are defined with local $\Z/2[\Z/2]$--coefficients:
\begin{eqnarray}\label{DZ2lokaa}
D_{i}(\J_a \times \JJ_a;\Z/2[\Z/2]) \subset H_{i}(K((\J_a \times \JJ_a) \int_{\chi^{[4]}} \Z,1);\Z/2[\Z/2]),
\end{eqnarray}
\begin{eqnarray}\label{DZ2lokad}
D_{i}(\H_{b \times \bb};\Z/2[\Z/2]) \subset H_{i}(K(\H_{b \times \bb} \int_{\chi^{[3]}} \Z,1);\Z/2[\Z/2]).
\end{eqnarray}

Analogously to $(\ref{deltaD})$, $(\ref{deltaD2})$ the following homomorphism (isomorphisms) are well-defined:
\begin{eqnarray}\label{deltaDaa}
\Delta^{[4]}: D_{i}(\J_a \times \JJ_a;\Z/2[\Z/2]) \to H_{i}(K(\J_a \times \JJ_a,1)),
\end{eqnarray}
\begin{eqnarray}\label{deltaDad}
\Delta^{[3]}: D_{i}(\H_{b \times \bb};\Z/2[\Z/2]) \to H_{i}(K(\H_{b \times \bb},1)).
\end{eqnarray}

The following lemma is analogous to Lemma 
$\ref{obstr}$.

\begin{lemma}\label{oaa}

--1. The group $H_{2s}(K((\J_a \times \JJ_a)\int_{\chi^{[4]}}\Z,1);\Z[\Z/2])$ 
(correspondingly, $H_{2s}(K(\H_{b \times \bb} \int_{\chi^{[3]}}\Z,1);\Z[\Z/2])$)
contains the direct factor 
$ D_{2s}(\J_a \times \JJ_a;\Z[\Z/2])$
(correspondingly, contains the direct factor $ D_{2s}(\H_{b \times \bb};\Z[\Z/2])$)

--2. The base of  the subgroup 
$\bigoplus_{i_1+i_2=2s} H_{i_1}(K(\J_a,1);\Z) \otimes H_{i_2}(K(\JJ_a,1);\Z) \subset D_{2s}(\J_a \times \JJ_a;\Z[\Z/2])$ determines the base 
of the subgroup $\mathrm{Im}(B) \cap D_{2s}(\I_{b \times \bb};\Z/2[Z/2])$, where
$$B: H_{2s}(K((\J_a \times \JJ_a)\int_{\chi^{[4]}}\Z,1);\Z[\Z/2]) \to H_{2s}(K((\J_a \times \JJ_a)\int_{\chi^{[4]}}\Z,1);\Z/2[\Z/2])$$
is the reduction homomorphism of modulo 2.
(Correspondingly, the base of the subgroup
$\bigoplus_{i_1+i_2=2s} H_{i_1}(K(\I_a,1);\Z) \otimes H_{i_2}(K(\Z/2,1);\Z) \subset D_{2s}(\H_{b \times \bb};\Z[\Z/2])$ 
(see the homomorphism $(\ref{omega3})$) determines the base 
of the subgroup %и элемент $R$ 
 $\mathrm{Im}(B) \cap D_{2s}(\I_{b \times \bb};\Z/2[\Z/2])$,
$$B: H_{2s}(K(\H_{b \times \bb} \int_{\chi^{[3]}}\Z,1);\Z[\Z/2]) \to H_{2s}(K(\H_{b \times \bb}\int_{\chi^{[3]}}\Z,1);\Z/2[\Z/2]).$$
%индуцированного приведением коэффициентов по модулю $2$. 
\end{lemma}
\[  \]

Define the epimorphism
\begin{eqnarray}\label{om4}
\omega^{[4]}: (\J_a \times \JJ_a) \int_{\chi^{[4]}} \Z \to \Z/4.
\end{eqnarray}
This epimorphism is defined by the extension of the homomorphism  $(\ref{omega4})$ 
from the subgroup
$\J_a \times \JJ_a$ to the hole group, the generator of the factor   $\Z$ 
is in the kernel of the epimorphism  $(\ref{om4})$.

Analogous  epimorphism
\begin{eqnarray}\label{om3}
\omega^{[3]}: \H_{b \times \bb} \int_{\chi^{[3]}} \Z \to \Z/4.
\end{eqnarray}
is well defined.

The representation $\Phi^{[2]}$ by the formula $(\ref{Phi2})$ is generalized as follows:
\begin{eqnarray}\label{Phi3}
\Phi^{[3]}: \H_{b \times \bb} \int_{\chi^{[3]}} \Z \to \Z/2^{[3]},
\end{eqnarray}
\begin{eqnarray}\label{Phi4}
\Phi^{[4]}: (\J_a \times \JJ_a) \int_{\chi^{[4]}} \Z \to \Z/2^{[4]},
\end{eqnarray}
where the image of the generator of the cyclic group
 $\Z$ in $\Z/2^{[3]}$ (correspondingly, in $\Z/2^{[4]}$) is represented by the automorphism $(\ref{chi3ZA})$ ($(\ref{chi4ZA})$).

The automorphisms  $\chi^{[i]}$, $i=2,3,4,$ in the images and in the preimages of the diagram $(\ref{a,aa})$
correspond with respect to horizontal homomorphisms. Therefore the following diagram of groups is well defined:
\begin{eqnarray}\label{a,aa,chi}
\begin{array}{ccc}
\qquad \I_{b \times \bb} \int_{\chi^{[2]}} \Z & \stackrel {\Phi^{[2]}} {\longrightarrow}& \qquad \Z/2^{[2]} \\
\\
i_{b \times \ bb,\H_{b \times \bb}} \downarrow &  & i^{[3]} \downarrow \\
\\
 \qquad \H_{b \times \bb} \int_{\chi^{[3]}} \Z &  \stackrel {\Phi^{[3]}}
{\longrightarrow}& \qquad \Z/2^{[3] } \\
\\
i_{\H_{b \times \bb},\J_a \times \JJ_a} \downarrow &  &  i^{[4]} \downarrow \\
\\
\qquad (\J_a \times \JJ_a) \int_{\chi^{[4]}} \Z & \stackrel{\Phi^{[4]}}{\longrightarrow} &  \qquad \Z/2^{[4]}.\\
\end{array}
\end{eqnarray}

%An inclusion $i_{d \times \dd,\_a \times \II_d}= i_{d \times
%\dd,\I_a} \times p_{d \times \dd,\II_d}: \I_d \times \II_d \subset
%\I_a \times \II_{d}$ is also defined. The homomorphism $i_{d
%\times \dd,\I_a}: \I_d \times \II_d \to \I_a$ is defined as the
%composition of the projection homomorphism $\I_d \times \II_d \to
%\I_d$ and the inclusion $\I_d \subset \I_a$. Moreover, the
%following diagram is commutative:

%\begin{eqnarray} \label{a,dd}
%\begin{array}{ccccc}
%\I_d \times \II_d & \stackrel  {i_{d \times \dd,\I_a \times
%\II_d}}{\longrightarrow} &
%\Z/2^{[2]} &  \stackrel {i_{diag,\Z/2^{[2]}}}{\longrightarrow}& \Z/2^{[2]} \times \Z/2^{[2]}  \\
%i_{d \times \dd,\I_a \times \II_d} \bigcap \qquad & & &  & \bar i_{[3]} \bigcap \\
%\I_a \times \II_d & & \stackrel{i_{\I_a \times \II_d}}{\longrightarrow} & & \Z/2^{[3]}.\\
%\end{array}
%\end{eqnarray}

\begin{definition}\label{[3,4]}
Let a  $\Z/2^{[3]}$--framed ($\Z/2^{[4]}$--framed) immersion
$(h,\Lambda,\zeta_L)$, $h: L^{n-4k} \looparrowright \R^n$
 ($h: L^{n-8k} \looparrowright \R^n$) represent an element
$z \in Imm^{\Z/2^{[3]}}(n-4k,4k)$ ($z \in
Imm^{\Z/2^{[4]}}(n-8k,8k)$). We say that this $\Z/2^{[3]}$--framed
($\Z/2^{[4]}$--framed) immersion is an $\H_{b \times \bb}$--framed
($\J_a \times \JJ_a$--framed) immersion if the
following two conditions are satisfied:
\[  \]
 
1. The structure mapping
$\zeta_L: L^{n-4k} \to K(\Z/2^{[3]},1)$ (correspondingly, $\zeta_L: L^{n-8k} \to K(\Z/2^{[4]},1)$) admits a reduction: this mapping is the composition of the mapping
\begin{eqnarray}\label{redad}
\zeta_{\H_{b \times \bb}}: L^{n-4k} \to K(\H_{b \times \bb} \int_{\chi^{[3]}} \Z,1)
\end{eqnarray}
(correspondingly,
\begin{eqnarray}\label{redaa}
\zeta_{\J_a \times \JJ_a}: L^{n-8k} \to K((\J_a \times \JJ_a) \int_{\chi^{[4]}} \Z,1))
\end{eqnarray}
 and the mapping $\Phi^{[3]}: K(\H_{b \times \bb} \int_{\chi^{[3]}} \Z,1) \to K(\Z/2^{[3]},1)$ $\quad$ 
 (correspondingly, the mapping $\Phi^{[4]}:
K((\J_a \times \JJ_a) \int_{\chi^{[4]}} \Z,1) \to K(\Z/2^{[4]},1)$).
\[  \]
 
 --2. The mapping $\bar{\zeta}_{\H_{b \times \bb}}: \bar{L}^{n-4k} \to K(\I_{b \times \bb} \int_{\chi^{[2]}} \Z,1)$ 
(correspondingly, the mapping $\bar \zeta_{\J_a \times \JJ_a}: \bar{L}^{n-8k} \to K(\I_{b \times \bb} \int_{\chi^{[2]}} \Z,1)$),
which is defined by means of the 2-sheeted (correspondingly, by means of the 4-sheeted) covering over the mapping
 $\zeta_{\H_{b \times \bb}}: L^{n-4k} \to K(\H_{b \times \bb} \int_{\chi^{[3]}} \Z,1)$ (correspondingly, over the mapping  $\zeta_{\J_a \times \JJ_a}: L^{n-8k} \to
K((\J_a \times \JJ_a) \int_{\chi^{[4]}} \Z,1)$)%, which (by Condition 1) satisfies Condition 1
%of Definition
%$\ref{ab}$, 
satisfies Condition 2 of Definition  $\ref{ab}$.
\end{definition}

Let us clarifies Condition 2 in the definition above.
Consider a submanifold
\begin{eqnarray}\label{LLad}
\bar L_{\eta^{6k}}^{n-16k} \subset \bar{L}^{n-4k},
\end{eqnarray}
which represents the homology Euler class
$[\eta_{\bar{L}_{\eta^{6k}}}]^{op} \in H_{n-16k}(\bar{L}^{n-4k};\Z/2)$
 of the bundle $3k\eta$
over $\bar L^{n-4k}$, where $2\eta = p^{\ast}(\zeta_L)$, $p: \bar L^{n-4k} \to L^{n-4k}$
 is the canonical 2-sheeted covering.

(Correspondingly, a submanifold
\begin{eqnarray}\label{LLaa}
\bar {L}_{\eta^{4k}}^{n-16k} \subset \bar{L}^{n-8k},
\end{eqnarray}
which represents the Euler class
$[\eta_{\bar{L}_{\eta^{4k}}}]^{op} \in H_{n-16k}(\bar{L}^{n-8k};\Z/2)$
of the bundle $4k\eta$ over $\bar L^{n-8k}$,where
$4\eta_{\bar{L}} = p^{\ast}(\zeta_L)$, $p: \bar L^{n-8k} \to L^{n-8k}$
 is the canonical $4$-sheeted covering.)

The restriction of the mapping
 $(\ref{redad})$ to the submanifold  $(\ref{LLad})$ in the regular $\Z/2^{[2]}$--framed cobordism class
 has to be represented by the composition
of the mapping
\begin{eqnarray}\label{Neta3}
\eta_{\bar{L}_{\eta^{6k}}}: \bar{L}_{\eta^{6k}}^{n-16k} \to K(\I_{b \times \bb},1)
\end{eqnarray}
and the standard inclusion
$K(\I_{b \times \bb},1)
\subset K(\I_{b \times
\bb} \int_{\chi^{[2]}} \Z,1)$.

(Correspondingly, the restriction of the mapping $(\ref{redaa})$ on the submanifold $(\ref{LLaa})$
in the regular $\Z/2^{[3]}$-framed cobordism class
 has to be presented by the composition
of the mapping
\begin{eqnarray}\label{Neta}
\eta_{\bar{L}_{\eta^{4k}}}: \bar{L}_{\eta^{4k}}^{n-16k} \to K(\H_{b \times \bb},1)
\end{eqnarray}
and the standard inclusion
 $K(\H_{b \times \bb},1)
\subset K(\H_{b \times
\bb} \int_{\chi^{[3]}} \Z,1)$.)
\[  \]

Let us investigate the characteristic class
 $\ref{thetaz2^3}$ for the cobordism group of $\H_{b \times \bb}$--framed (correspondingly, $\J_a \times \JJ_a$--framed)
immersions.

The cohomology group $H^4(K(\H_{b \times \bb} \int_{\chi^{[3]}},1);\Z/2)$ (correspondingly, $H^8(K((\J_a
\times \JJ_a) \int_{\chi^{[4]}},1);\Z/2)$ 
  contains an element $\tau_{\H_{b \times \bb}}$ (correspondingly, $\tau_{\J_a \times \JJ_a}$),
   which is defined in
by the equation ($\ref{i3}$) (correspondingly, ($\ref{i4}$)) below.

Consider the mapping $\Phi^{[3]}: K(\H_{b \times \bb}\int_{\chi^{[3]}},1) \to K(\Z/2^{[3]},1)$ 
(correspondingly, $\Phi^{[4]}: K((\J_a
\times \JJ_a) \int_{\chi^{[4]}},1) \to K(\Z/2^{[4]},1)$)
 and consider the pull-back
$(\Phi^{[3]})^{\ast}(\tau_{[3]})$  (correspondingly, $(\Phi^{[4]})^{\ast}(\tau_{[4]})$) 
of the characteristic Euler class
$\tau_{[3]} \in H^4(K(\Z/2^{[3]},1);\Z/2)$ (correspondingly, $\tau_{[4]} \in
H^8(K(\Z/2^{[4]},1);\Z/2)$) of the universal bundle. 

Define
\begin{eqnarray}\label{i3}
 (\Phi^{[3]})^{\ast} (\tau_{[3]}) = \tau_{\H_{b \times \bb}} \in H^4(K(\H_{b \times \bb}\int_{\chi^{[3]}}\Z,1);\Z/2),
\end{eqnarray}
(correspondingly,
\begin{eqnarray}\label{i4}
 \Phi^{[4]})^{\ast} (\tau_{[4]}) = \tau_{\J_a \times \JJ_a} \in H^8(K((\J_a \times \JJ_a)\int_{\chi^{[4]}}\Z,1);\Z/2).
\end{eqnarray}

In section 1, for a
 $\Z/2^{[s+1]}$---framed immersion $(h,\Lambda,\zeta_L)$, together with a
$2^{s}$--dimensional characteristic class
 $\zeta_L \in H^{2^{s}}(L^{n-k2^{s}};\Z/2)$,
 we also considered
a $2$-dimensional characteristic class $\bar \zeta_{[2],L} \in
H^2(\bar L_{[2]}^{n-k2^{s}k};\Z/2)$.

For a mapping $\zeta_{\H_{b \times \bb}}:
L^{n-4k} \to K(\H_{b \times \bb}\int_{\chi^{[3]}} \Z,1)$
(correspondingly, for a mapping $\zeta_{\J_a
\times \JJ_a}: L^{n-8k} \to K((\J_a \times \JJ_a)\int_{\chi^{[4]}} \Z,1)$) as an analog
of the characteristic class $\bar
\zeta_{[2],L}^{\ast}(\tau_{[2]})$ there serves the characteristic
class $\bar \zeta_{b \times \bb}^{\ast}(\tau_{b \times \bb}) \in
H^2(\bar L^{n-4k};\Z/2)$, for $s=3$ (correspondingly, $\bar \zeta_{a
\times \aa}^{\ast}(\tau_{a \times \aa}) \in H^2(\bar L^{n-8k};\Z/2)$, for $s=4$, in this formula the  covering
$\bar L^{n-8k} \to L^{n-8k}$ is a 4-sheeted covering). 
%We now define this 2-dimensional
%characteristic class.

Define the mapping $\bar \zeta_{\H_{b
\times \bb}}$ as 2-sheeted covering over the mapping 
$\zeta_{\H_{b \times \bb}}$ with respect to the subgroup $i_{b \times \bb,\H_{b \times
\bb}}: \I_{b \times \bb}\int_{\chi^{[2]}} \Z \subset \H_{b \times \bb}\int_{\chi^{[3]}} \Z $.
Define the mapping $\bar \zeta_{\J_a \times \JJ_a}$
as 2-sheeted covering over the mapping 
$\zeta_{\J_a \times \JJ_a}$ with respect to the subgroup $i_{b \times \bb,a \times \aa}: 
\I_{b \times \bb}\int_{\chi^{[2]}} \Z \subset (\J_{a} \times J_{\aa})\int_{\chi^{[4]}} \Z $.

The characteristic class $\bar \zeta_{\H_{b \times \bb}}^{\ast}(\tau_{b \times \bb})$
(correspondingly, the class $\bar \zeta_{\J_a \times \JJ_a}^{\ast}(\tau_{b \times \bb})$)
  is induced from the universal
class $\tau_{b \times \bb} \in
H^2(K(\I_{b \times \bb} \int_{\chi^{[2]}} \Z),1);\Z/2)$
by the mapping $\bar \zeta_{\H_{b
\times \bb}}: \bar L^{n-4k} \to K(\I_{b \times
\bb}\int_{\chi^{[2]}} \Z),1)$ (correspondingly, by the mapping $\bar \zeta_{\J_a \times \JJ_a}: 
\bar L^{n-8k} \to K(\I_{b \times \bb}\int_{\chi^{[2]}} \Z,1)$, in this formula the mapping $\bar \zeta_{\J_a \times \JJ_a}$
is the 4-sheeted covering over the  mapping $\zeta_{\J_a \times \JJ_a}$ ).
We need to define dual homology classes and the analogous formulas to $(\ref{etadlok})$.

Let us consider the immersion
 $h: L^{n-4k}_{\H_{b \times \bb}} \looparrowright \R^n$ (correspondingly, the immersion 
$h: L^{n-8k}_{\J_a \times \JJ_a} \looparrowright \R^n$) as in Definition  $\ref{[3,4]}$.
Let the mapping 
\begin{eqnarray}\label{1zetaad}
\zeta_{\H_{b \times \bb}}: L_{\H_{b \times \bb}}^{n-4k} \to K(\H_{b \times \bb}\int_{\chi^{[3]}} \Z,1)
\end{eqnarray}
be given (correspondingly, the mapping
\begin{eqnarray}\label{1zetaaa}
\zeta_{\J_a \times \JJ_a}: L^{n-8k}_{\J_a \times \JJ_a} \to K((\J_a \times \JJ_a)\int_{\chi^{[4]}} \Z,1))
\end{eqnarray}
be given).

Over the target space of the mapping
 $(\ref{1zetaad})$ the universal $4$-dimensional $\H_{b \times \bb}\int_{\chi^{[3]}} \Z$-bundle 
 is well defined. 
 (Correspondingly, Over the target space of the mapping $(\ref{1zetaaa})$ the universal $8$-dimensional 
 $(\J_a \times \JJ_a)\int_{\chi^{[3]}} \Z$-bundle 
 is well defined.) 
 Assume that the manifold
 $L_{\H_{b \times \bb}}^{n-4k}$ (correspondingly, the manifold
$L^{n-8k}_{\J_a \times \JJ_a}$)  is connected, and the mapping
$(\ref{1zetaad})$ (correspondingly, the mapping $(\ref{1zetaaa})$) is punctured.

Let us define the element
\begin{eqnarray}\label{etaadlok}
\zeta^{loc}_{\H_{b \times \bb},\ast}([L_{\H_{b \times \bb}}]) \in H_{n-4k}(K(\H_{b \times \bb} \int_{\chi^{[3]}} \Z ,1);\Z/2[\Z/2])
\end{eqnarray}
(correspondingly, the element 
\begin{eqnarray}\label{etaaalok}
\zeta^{loc}_{\J_a \times \JJ_a,\ast}([L_{\J_a \times \JJ_a}]) \in H_{n-8k}(K((\J_a \times \JJ_a) \int_{\chi^{[4]}} \Z ,1);\Z/2[\Z/2]))
\end{eqnarray}
analogously to the formula
$(\ref{Xi})$ for the element  $(\ref{etadlok})$. Below this element is defined even 
for a non-punctured reduction mapping, assuming that the target manifold is the self-intersection manifold
of a connected $\Z/2^{[2]}$-framed (correspondingly, $\Z/2^{[3]}$-framed) immersed manifold. 
%for which a $\I_{b \times \bb} \int_{\chi^{[2]}}\Z$--reduction (correspondingly, a $\H_{b \times \bb} %\int_{\chi^{[3]}}\Z$--reduction) of the classifying mapping of the framing
%is well-defined.

Let us assume that a 
$\Z/2^{[3]}$--framed immersion  $h: L^{n-4k} \looparrowright \R^n$ (correspondingly, a 
$\Z/2^{[4]}$--framed immersion  $h: L^{n-8k} \looparrowright \R^n$)
is the immersion of self-intersection manifold of a 
$\Z/2^{[2]}$--framed immersion  $(g,\eta_N,\Psi)$, $g: N^{n-2k} \looparrowright \R^n$
(correspondingly, of a $\Z/2^{[3]}$--framed immersion $(g,\eta_N,\Psi)$, $g: N^{n-4k} \looparrowright \R^n$).

Assume that in the manifold
 $N^{n-2k}$ a closed connected component  $N_{b \times \bb}^{n-2k} \subset N^{n-2k}$
 is marked, comp. with the formula 
$(\ref{abelcomp})$ (correspondingly, in the manifold 
$N^{n-4k}$ a closed connected component  $N_{b \times \bb}^{n-4k} \subset N^{n-4k}$ is marked). 
Moreover, a punctured mapping
\begin{eqnarray}\label{1etadd}
\eta_{b \times \bb}: N^{n-2k}_{b \times \bb} \to K(\I_{b \times \bb} \int_{\chi^{[2]}} \Z,1)
\end{eqnarray}
(correspondingly, a punctured mapping 
\begin{eqnarray}\label{1etaad}
\eta_{b \times \bb}: N^{n-4k}_{b \times \bb} \to K(\H_{b \times \bb} \int_{\chi^{[3]}} \Z,1)),
\end{eqnarray}
which determines a reduction of the characteristic mapping of the $\Z/2^{[2]}$--framing
(correspondingly, of the $\Z/2^{[3]}$--framing) over the marked component is given.

Denote by  $L_{b \times \bb}^{n-4k} \subset L^{n-4k}$ (correspondingly, by $L_{b \times \bb}^{n-8k} \subset L^{n-8k}$) 
the component of the self-intersection manifold of the immersion $g$, restricted to the marked component
$N_{b \times \bb}^{n-2k}$ (correspondingly, to the marked component $N_{b \times \bb}^{n-4k}$).
Assume that the manifold  $L_{b \times \bb}^{n-4k} \subset L^{n-4k}$ is decomposed into 2 subcomponents as in the following formula: 
\begin{eqnarray}\label{compL3a}
L_{b \times \bb}^{n-4k} = L_{\H_{b \times \bb}}^{n-4k} \cup L_{b \times \bb,[3]}^{n-4k}
\end{eqnarray}
(correspondingly, the manifold $L_{b \times \bb}^{n-8k} \subset L^{n-8k}$ is decomposed into 2 subcomponents as in the following formula: 
\begin{eqnarray}\label{compL4a}
L_{b \times \bb}^{n-8k} = L_{\J_a \times \JJ_a}^{n-8k} \cup L_{b \times \bb,[4]}^{n-8k}).
\end{eqnarray}

Let us assume that the first component $L^{n-4k}_{\H_{b \times \bb}}$ in the formula $(\ref{compL3a})$, 
(correspondingly, to the first component $L^{n-8k}_{\J_a \times \JJ_a}$ in the formula $(\ref{compL4a})$)
which
generally speaking, is non-connected and a reduction mapping $(\ref{1zetaad})$ 
(correspondingly, a reduction mapping $(\ref{1zetaaa})$)
of the characteristic mapping is given,
we do not assume that the reduction mapping is punctured. 
 
 Let us consider the immersion of the canonical 2-sheeted covering
\begin{eqnarray}\label{barL4k}
\bar{L}_{\H_{b \times \bb}}^{n-4k} \looparrowright N_{b \times \bb}^{n-2k}
\end{eqnarray}
(correspondingly, the immersion of the canonical 2-sheeted covering
\begin{eqnarray}\label{barL8k}
\bar{L}_{\J_a \times \JJ_a}^{n-8k} \looparrowright N_{b \times \bb}^{n-4k}
\end{eqnarray}
over the self-intersection manifold $L_{\H_{b \times \bb}}^{n-4k}$ (correspondingly, over the self-intersection manifold $L^{n-8k}_{\J_a \times \JJ_a}$).
 
Denote by  
\begin{eqnarray}\label{2etadd}
\eta_{\H_{b \times \bb}}: \bar L^{n-4k}_{\H_{b \times \bb}} \to K(\I_{b \times \bb} \int_{\chi^{[2]}} \Z,1)
\end{eqnarray}
(correspondingly, by 
\begin{eqnarray}\label{2etaad}
\eta_{\J_a \times \JJ_a}: \bar L^{n-8k}_{\J_a \times \JJ_a} \to K(\H_{b \times \bb} \int_{\chi^{[3]}} \Z,1))
\end{eqnarray}
the restriction of the mapping 
$(\ref{1etadd})$ (correspondingly, the restriction of the mapping $(\ref{1etaad})$)
on the 2-sheeted covering.

Assume, that the mapping
 $(\ref{2etadd})$ is homotopic to the corresponding 2-sheeted covering over the  mapping$(\ref{1zetaad})$.
(Correspondingly,  assume, that the mapping
 $(\ref{2etaad})$ is homotopic to the corresponding 2-sheeted covering over the mapping $(\ref{1zetaaa})$.) 

The characteristic mapping
$(\ref{1etadd})$ (correspondingly, the mapping  $(\ref{2etaad})$) determines the following homology class 
\begin{eqnarray}\label{baretadlok}
\eta^{loc}_{\H_{b \times \bb},\ast}([\bar{L}_{\H_{b \times \bb}}]) \in H_{n-4k}(K(\I_{b \times \bb} \int_{\chi^{[2]}} \Z ,1);\Z/2[\Z/2]),
\end{eqnarray}
which coincides with the transfer to the corresponding 2-sheeted covering of the element $(\ref{etaadlok})$.
(Correspondingly the homology class
\begin{eqnarray}\label{baretaalok}
\eta^{loc}_{\J_a \times \JJ_a,\ast}([\bar{L}_{\J_a \times \JJ_a}]) \in H_{n-8k}(K(\H_{b \times \bb} \int_{\chi^{[3]}} \Z ,1);\Z/2[\Z/2])
\end{eqnarray}
is well-defined. This homology class
coincides with the transfer to the corresponding 2-sheeted covering of the element
 $(\ref{etaaalok})$.)

\begin{definition}\label{ad-struct}
Let $(g,\Psi,\eta_N)$ be a
$\Z/2^{[2]}$--framed immersion, $g:
N^{n-2k} \looparrowright \R^n$, which represents an element  $y \in
Imm^{\Z/2^{[2]}}(n-2k,2k)$, assuming  $n > 16k$. Let $(h,\Lambda,\zeta_L)$
be a
$\Z/2^{[3]}$--framed immersion, $h:
L^{n-4k} \looparrowright \R^n$, is an immersion of self-intersection points of the immersion  $g$,
which represents an element
$z=\delta^{\Z/2^{[3]}}_{k}(y) \in Imm^{ \Z/2^{[3]}}(n-4k,4k)$.

Assume that a closed component $N^{n-2k}_{b \times \bb}$ of the manifold
$N^{n-2k}$ is punctured (comp. with the formula $(\ref{abelcomp})$). 
Assume that the self-intersection manifold
$L_{b \times \bb}^{n-4k}$ of the immersion 
$g \vert_{N^{n-2k}_{b \times \bb}}$ is decomposed into two components like in the formula
$(\ref{compL3a})$.

We say that  $\Z/2^{[2]}$--framed immersion
$(g,\Psi,\eta_N)$ admits an $\H_{b \times \bb}$--structure, if on the component  $N^{n-2k}_{b \times \bb}$
a punctured mapping $(\ref{1etadd})$ is given, and this mapping determines a reduction of the restriction of the characteristic mapping  $\eta_N$; on the component $L^{n-4k}_{\H_{b \times \bb}}$  a mapping  $(\ref{1zetaad})$
is given (we do not assume that this mapping is punctured), which determines a reduction of the restriction of the classifying mapping  $\zeta_L$ and the canonical double covering over this reduction mapping is homotopic to the mapping $(\ref{1etadd})$. Additionally, the following 3 conditions are satisfy. 

--1. The $\Z/2^{[2]}$--framed immersion    $(g,\Psi,\eta_N)$ satisfies Condition 2 from Definition
$\ref{strucIdId}$ (recall, that this condition means that the Arf-invariant 
is calculated using the component  $N^{n-2k}_{b \times \bb}$).

--2. The punctured mapping 
$(\ref{1etadd})$ satisfies Property 2 of Lemma   $\ref{dd}$ (recall, that this condition means that 
the element 
$\eta_{b \times \bb,\ast}([N^{n-2k}_{b \times \bb},pt]) \in H_{n-2k}(K(\I_{b \times \bb} \int_{\chi^{[2]}} \Z,1);\Z[\Z/2])$, which is constructed by the formula $(\ref{etadd16})$ for the mapping  $(\ref{1etadd})$,
belongs to the subgroup $D_{n-2k}(\I_{b \times \bb};\Z/2[\Z/2])$.

--3. The element  
\begin{eqnarray}\label{2etaddbarL}
\eta^{loc}_{b \times \bb,\ast}([\bar{L}_{b \times \bb}]) \in H_{n-4k}(K(\I_{b \times \bb} \int_{\chi^{[2]}} \Z ,1);\Z/2[\Z/2]),
\end{eqnarray}
which is defined analogously to the element 
 $(\ref{baretadlok})$, but using the both components from the formula $(\ref{compL3a})$,%$(\ref{abelcomp})$ 
and the element 
$(\ref{baretadlok})$, which is constructed only for the first component 
in the formula  $(\ref{compL3a})$ %$(\ref{abelcomp})$ 
are related in the group
$H_{n-4k}(K(\I_{b \times \bb} \int_{\chi^{[2]}} \Z ,1);\Z/2[\Z/2])$ by the following equation: 
\begin{eqnarray}\label{dd!}
\eta^{loc}_{b \times \bb,\ast}([\bar{L}_{b \times \bb}]) = r^!(\zeta^{loc}_{\H_{b \times \bb},\ast}([L_{\H_{b \times \bb}}])),
\end{eqnarray}
where
\begin{eqnarray}\label{r!}
r^!: H_{n-4k}(K(\H_{b \times \bb} \int_{\chi^{[3]}} \Z) ,1);\Z/2[\Z/2]) \to
\end{eqnarray}
$$H_{n-4k}(K(\I_{b \times \bb} \int_{\chi^{[2]}} \Z) ,1);\Z/2[\Z/2])$$
is the transfer homomorphism, which is associated with the right upper inclusion of the subgroup $2$
in the diagram 
$(\ref{a,aa,chi})$, this inclusion is re-denoted by
\begin{eqnarray}\label{rast}
r: \I_{b \times \bb} \int_{\chi^{[2]}} \Z  \subset \H_{b \times \bb} \int_{\chi^{[3]}} \Z 
\end{eqnarray}
for short.
\end{definition}

Let us express in the framework of Definition 
 $\ref{ad-struct}$ the homology class  $(\ref{etadlok})$ from the homology class
$(\ref{baretadlok})$  and the normal Euler class of the immersion.
%Наряду с элементом $(\ref{baretadlok})$ (соответственно, $(\ref{baretaalok})$) рассматривать
Define the element 
\begin{eqnarray}\label{3baretaddlok}
\eta^{loc}_{b \times \bb,\ast}([\bar{\zeta}_{L}^{3k}]^{op}) \in H_{m_{\sigma}}(K(\I_{b \times \bb} \int_{\chi^{[2]}} \Z) ,1);\Z/2[\Z/2]),
\end{eqnarray}
where the homology class  $([\bar{\zeta}_{L}^{3k}]^{op}) \in H_{m_{\sigma}}(\bar{L}_{\H_{b \times \bb}}^{n-4k};\Z/2[\Z/2])$ is defined as the result of the intersection of the fundamental class
$[\bar{L}_{\H_{b \times \bb}}] \in H_{n-4k}(\bar{L}_{\H_{b \times \bb}}^{n-4k};\Z/2[\Z/2])$ 
with the  Euler class
$\bar{\zeta}_{L}^{3k} \in H^{12k}(\bar {L}_{\H_{b \times \bb}}^{n-4k};\Z/2[\Z/2])$ %(соответственно, $\bar{\zeta}_{\bar 
of the bundle $3k\bar{\zeta}_{L}$. 

\begin{lemma}\label{ad}
Assume that a
 $\Z/2^{[2]}$--framed immersion $(g,\Psi,\eta_N)$ is given, and this immersion allows $\H_{b \times \bb}$--structure in the sense of Definition $\ref{ad-struct}$. Then the following properties are satisfied:

--1. The element $(\ref{3baretaddlok})$ belongs to the subgroup  $(\ref{DZ2lok})$, $i=m_{\sigma}$ 
and this element is lifted to the subgroup $(\ref{DZlok})$. The image of this  element  in the group
$H_{m_{\sigma}}(K(\I_{b \times \bb},1))$ by means of the homomorphism 
$(\ref{deltaD2})$ satisfies Property 1 fro Lemma $\ref{dd}$, analogously to the element $(\ref{etadlok})$.

--2. The element  $\zeta^{loc}_{L_{\H_{b \times \bb}},\ast}([L_{\H_{b \times \bb}}])$,
which is defined by the formula
$(\ref{etaadlok})$, belongs to the subgroup  $D_{n-4k}(\H_{b \times \bb};\Z/2[\Z/2])$, 
which is defined by the formula $(\ref{DZ2lokad})$ for $i=n-4k$, and this elements is lifted to the subgroup 
$(\ref{DZlokad})$.

--3. The Arf-invariant of the $\Z/2^{[2]}$--framed immersion  $(g,\Psi,\eta_N)$ coincides with the Arf-invariant which is calculated by the formula $(\ref{thetaz2^3})$ only for  $\Z/2^{[3]}$--framed immersed component  $L^{n-4k}_{\H_{b \times \bb}}$.
\end{lemma}

\subsubsection*{Proof of Lemma  $\ref{ad}$}
Prove Statement 1.
Consider the immersion 
$g_{b \times \bb}: N_{b \times \bb}^{n-2k} \looparrowright \R^n$ of the marked component and apply
to this immersion the Herbert Theorem with the local coefficient system $\Z/2[\Z/2]$. 
(The statement and a proof of this version of the Herbert Theorem 
is analogous to the statement and the proof of the Herbert Theorem with $\Z/2$-coefficients, see [A1], Proposition 8, and the reference there.) We get that homology class
\begin{eqnarray}\label{barL_b}
[\bar{L}_{b \times \bb}] \in %\mathrm{Im}(H_{n-4k}(\bar{L_b}^{n-4k};\Z/2[\Z/2]) \to
H_{n-4k}(N_{b \times \bb}^{n-2k};\Z/2[\Z/2]))
\end{eqnarray} 
coincides with the homology class %(гомоморфизм групп гомологий индуцирован погружением $(\ref{barL4k})$)
$[\eta_{N_{b \times \bb}}^{k}]^{op} \in H_{n-4k}(N_{b \times \bb}^{n-2k};\Z/2[\Z/2])$.
By assumption, the equation $(\ref{dd!})$ is satisfied, therefore we have: 
%образы классов гомологий $(\ref{barL_b})$ и $[\bar{L}_{b \times \bb}]$ в группе
%$H_{n-4k}(K(\I_{b \times \bb} \int_{\chi^{[2]}} \Z) ,1);\Z/2[\Z/2])$
%равны. Поэтому возможна замена в левой части. 
\begin{eqnarray}\label{eqvbarL_b}
\eta_{b \times \bb,\ast}^{loc}([\eta_{N_{b \times \bb}}^{k}]^{op}) = \eta_{b \times \bb,\ast}^{loc} ([\bar{L}_{\H_{b \times \bb}}]).
\end{eqnarray}
In this formula we re-express the right side using the transfer homomorphism:
\begin{eqnarray}\label{barL_b}
\eta_{b \times \bb,\ast}^{loc}([\eta_{N_{b \times \bb}}^{k}]^{op}) = r^!(\zeta_{\H_{b \times \bb},\ast}^{loc}([L_{\H_{b \times \bb}}])).
\end{eqnarray}
%поскольку выполнено равенство $r^!(\zeta_{\H_{b \times \bb},\ast}^{loc}([L_{\H_{b \times \bb}}])=\eta_{b \times %\bb,\ast}^{loc}([\bar{L}_{b \times \bb}])$, где гомоморфизм $r$ определен по формуле $(\ref{rast})$.
%$[\bar{L_b}]) \in H_{n-4k}(\bar{L}_{\bar \zeta};\Z/2[\Z/2])$ (это класс гомологий
%равен классу гомологий $r^!_{\zeta}([L])$, соответствующему  правой части равенства $(\ref{dd!})$ и
%$[\eta_N^{k}]^{op} \in H_{n-4k}(N^{n-2k};\Z/2[\Z/2])$ (этот класс гомологий
%равен классу гомологий $[\bar{L}]$, соответствующему  левой части равенства $(\ref{dd!})$.

Consider the product of the homology classes in the both sides of the equation 
with the cohomology class 
$\eta^{\ast}_{b \times \bb} \circ r^{\ast}(\zeta^{3k}) \in H^{12k}(N^{n-2k}_{b \times \bb};\Z/2[\Z/2])$.
In the right side of the equation we have the homology class  $r^{!}(\zeta^{loc}_{\H_{b \times \bb},\ast}([L_{\H_{b \times \bb}}])) \cap r^{\ast}((\zeta^{[3]})^{3k})$, which coincides with the homology class 
$r^!(\zeta^{loc}_{\H_{b \times \bb},\ast}([L_{\H_{b \times \bb}}]) \cap (\zeta^{[3]})^{3k})$,
namely, with the homology class 
$(\ref{3baretaddlok})$.
In the right side of the equation we have the homology class
$\eta_{b \times \bb,\ast}^{loc}([\eta_{N_{b \times \bb}}^{7k}]^{op})$. 
Let us apply the homomorphism $(\ref{deltaD2})$ to this homology classes, we get that the class coincides with
  $(\ref{etadlok})$. For this class the required property is proved in Lemma $\ref{dd}$.
Statement 1 is proved.

Prove Statement 2.
  In the case $\sigma \ge 5$ the codimension of  $\Z/2^{[3]}$--framed immersed manifold 
   $L^{n-4k}_{\H_{b \times \bb}}$ is even, therefore the considered manifold is oriented and its fundamental class
  belongs to the homology group with integer coefficients. 

The transfer homomorphism 
$$ r^!: H_{n-4k}(\H_{b \times \bb} \int_{\chi^{[3]}} \Z;\Z/2[\Z/2]) \to H_{n-4k}(\I_{b \times \bb} \int_{\chi^{[2]}} \Z;\Z/2[\Z/2])$$
with $\Z/2[\Z/2]$--coefficients is a monomorphism. This follows from the following fact: the element $R$,
which is described in Lemma $\ref{obstr}$, belongs to the image of the transfer homomorphism.
Therefore it is sufficient  to prove that the element 
$(\ref{baretadlok})$ belongs to the subgroup  $(\ref{DZ2lok})$, $i=n-4k$. Using the equation
$(\ref{dd!})$ it is sufficient to prove that the element  $(\ref{2etaddbarL})$ belongs to the considered group.
The homology class  $(\ref{2etaddbarL})$ is expressed from the homology class $(\ref{etadd16})$. 
By Property 1 from Definition $(\ref{ad-struct})$ the homology class $(\ref{etadd16})$
belongs to the subgroup 
$(\ref{DZ2lok})$, $i=n-2k$. Statement 2 is proved.

Statement 3 follows from Statement 2 above,  and the property of the homology class
 $(\ref{etadd16})$,which is formulated in Statement 1 of Lemma $\ref{dd}$.
 Lemma  $\ref{ad}$ is proved.
\[  \]

\begin{definition}\label{aa-struct}
Let $(g,\Psi,\eta_N)$ be a
$\Z/2^{[3]}$--framed immersion, $g:
N^{n-4k} \looparrowright \R^n$, which represents an element  $y \in
Imm^{\Z/2^{[3]}}(n-4k,4k)$, assuming  $n > 16k$. Let $(h,\Lambda,\zeta_L)$
be a
$\Z/2^{[4]}$--framed immersion, $h:
L^{n-8k} \looparrowright \R^n$, is an immersion of self-intersection points of the immersion  $g$,
which represents an element
$z=\delta^{\Z/2^{[4]}}_{k}(y) \in Imm^{ \Z/2^{[4]}}(n-8k,8k)$.

Assume that a component $N^{n-4k}_{b \times \bb}$ of the manifold  $N^{n-4k}$ is marked and punctured. 
Assume that the self-intersection manifold
$L_{b \times \bb}^{n-8k}$ of the immersion  $g \vert_{N^{n-4k}_{b \times \bb}}$ is decomposed into two components, like in the formula  $(\ref{compL4a})$.

We will say that the
$\Z/2^{[3]}$--framed immersion 
$(g,\Psi,\eta_N)$ admits an $\J_a \times \JJ_a$--structure, if  
on the component  $N^{n-4k}_{b \times \bb}$ a punctured mapping
$(\ref{1etaad})$ is well-defined, and this mapping determines a reduction of the restriction of the characteristic mapping $\eta_N$;  on the component $L^{n-8k}_{\J_a \times \JJ_a}$  a mapping  $(\ref{1zetaaa})$
is given (we do not assume that this mapping is punctured), which determines a reduction of the restriction of the classifying mapping  $\zeta_L$ and the canonical double covering over this reduction mapping is homotopic to the mapping $(\ref{1etaad})$. Additionally, the following 3 conditions are satisfy. 

--1. The mapping $(\ref{2etadd})$  satisfies Condition 3 from Lemma $\ref{ad}$ 
(recall, that this condition means that the Arf-invariant 
is calculated using the component  $N^{n-4k}_{b \times \bb}$).

--2. The mapping $(\ref{2etadd})$ satisfies Condition  2 from Lemma $\ref{ad}$.

--3. The element  
\begin{eqnarray}\label{3etaddbarL}
\eta^{loc}_{b \times \bb,\ast}([\bar{L}_{b \times \bb}]) \in H_{n-8k}(K(\H_{b \times \bb} \int_{\chi^{[3]}} \Z ,1);\Z/2[\Z/2]),
\end{eqnarray}
which is defined analogously to the element 
 $(\ref{baretaalok})$, but using the both components from the formula  $(\ref{compL4a})$,%$(\ref{abelcomp})$ 
and the element $(\ref{baretaalok})$, which is constructed only for the first component 
in the formula  $(\ref{compL4a})$ %$(\ref{abelcomp})$ 
are related in the group $H_{n-8k}(K(\H_{b \times \bb} \int_{\chi^{[3]}} \Z ,1);\Z/2[\Z/2])$
by the following relation:
\begin{eqnarray}\label{dda!}
\eta^{loc}_{b \times \bb,\ast}([\bar{L_{b \times \bb}}]) = r^!(\zeta^{loc}_{\J_a \times \JJ_a,\ast}([L_{\J_a \times \JJ_a}])),
\end{eqnarray}
where
\begin{eqnarray}\label{r!}
r^!: H_{n-8k}(K((\J_a \times \JJ_a) \int_{\chi^{[4]}} \Z) ,1);\Z/2[\Z/2]) \to
\end{eqnarray}
$$H_{n-8k}(K(\H_{b \times \II_b} \int_{\chi^{[3]}} \Z) ,1);\Z/2[\Z/2])$$
is the transfer homomorphism, which is associated with the right middle inclusion of the subgroup $2$
in the Diagram
$(\ref{a,aa,chi})$, which is re-denoted by 
\begin{eqnarray}\label{rast3}
r: \H_{b \times \bb} \int_{\chi^{[3]}} \Z  \subset (\J_a \times \JJ_a) \int_{\chi^{[4]}} \Z
\end{eqnarray}
for short.
\end{definition}

\begin{lemma}\label{aa}
Assume a
 $\Z/2^{[3]}$--framed immersion  $(g,\Psi,\eta_N)$ is given, which admits an $\J_a \times \JJ_a$--structure
in the sense of Definition  $\ref{aa-struct}$. The following two  properties are satisfied:

--1. In the group $H_{m_{\sigma}}(K(\H_{b \times \bb} \int_{\chi^{[3]}} \Z,1);\Z/2[\Z/2])$
the element $(\ref{baretaalok})$ is equal to the element $\eta_{b \times \bb,\ast}^{loc}([\eta_{N_{b \times \bb}}^{3k}]^{op})$. In particular, the images of these elements by the homomorphism 
$(\ref{deltaDad})$ for $i=m_{\sigma}$ are equal in the group 
 $H_{m_{\sigma}}(K(\H_{b \times \bb},1))$.

--2. The Arf-invariant of the $\Z/2^{[3]}$--framed immersion  $(g,\Psi,\eta_N)$ coincides with the Arf-invariant, which is calculated by the formula  $(\ref{thetaz2^3})$, using only the component  $L^{n-8k}_{\J_a \times \JJ_a}$.
\end{lemma}

\subsubsection*{Proof of Lemma $\ref{aa}$}
The proof is analogous to the proof of Lemma 
$\ref{ad}$.
\[  \]

\begin{example}\label{IaIa}
Let the $\Z/2^{[2]}$--framed (correspondingly, $\Z/2^{[3]}$--framed) immersion
$(g,\eta_N,\Psi,)$, $g: N^{n-2k} \looparrowright \R^n$ (correspondingly, $g:
N^{n-4k} \looparrowright \R^n$) be represented an element $y \in
Imm^{\Z/2^{[2]}}(n-2k,2k)$ (correspondingly, an element $y \in Imm^{\Z/2^{[3]}}(n-4k,4k)$), $n
> 16k$.
Assume that the considered immersion is a $\I_{b \times \bb}$--framed
immersion (correspondingly, is a $\H_{b \times \bb}$--framed immersion) in the sense of Definition $\ref{ab}$
 (correspondingly, in the sense of Definition $\ref{[3,4]}$.
Let a $\Z/2^{[3]}$--framed (correspondingly, a $\Z/2^{[4]}$--framed) immersion
$(h,\zeta_L,\Lambda)$, $h: L^{n-4k} \looparrowright \R^n$ (correspondingly, $h:
L^{n-8k} \looparrowright \R^n$) be represented the element $z =
\delta^{\Z/2^{[3]}}_k(y) \in Imm^{\Z/2^{[3]}}(n-4k,4k)$ (correspondingly, $z =
\delta^{\Z/2^{[4]}}_k(y) \in Imm^{\Z/2^{[4]}}(n-8k,8k)$) and be the
immersion of the self-intersection manifold of the immersion $(g,\eta_N,\Psi)$.
Assume that the immersion $(h,\zeta_L,\Lambda)$ is a 
$\H_{b \times \bb}$--framed immersion  (correspondingly, is a $\J_a \times
\JJ_a$--framed immersion) in the sense of Definition  $\ref{[3,4]}$. Then the $\Z/2^{[2]}$--framed 
(correspondingly, the  $\Z/2^{[3]}$--framed) immersion $(g,\Psi,\eta_N)$ admits a 
$\H_{b \times \bb}$--structure  (correspondingly, a  $\J_a \times \JJ_a$--structure), which is defined
by the reduction mapping  $\zeta_{\H_{b \times \bb}}$ (correspondingly, by  $\zeta_{\J_a
\times \JJ_a}$) of the characteristic mapping  $\zeta_L$ (comp. with Example 19 in [A1]).
\end{example}

\subsubsection*{Justification of the example $\ref{IaIa}$}
The example is obvious.
\[  \]

The following theorems are analogs of Theorem $\ref{th6}$.

\begin{lemma}\label{Hbb}

Assume that the
 $\Z/2^{[2]}$--framed immersion
$(g,\Psi,\eta_N)$ represents an element $y \in
Imm^{\Z/2^{[2]}}(n-\frac{n-n_\sigma}{8}, \frac{n-n_\sigma}{8})$, $n \ge 254$,
see $(\ref{check})$, and a reduction of the characteristic mapping $\eta_N$ by 
the mapping $(\ref{1etadd})$ is given, such that Conditions 1,2, from Definition  $\ref{ad-struct}$
are satisfied.
Then the element $\delta^{\Z/2^{[3]}}_{k}(y)$ in the group
$Imm^{\Z/2^{[3]}}(n-\frac{n-m_{\sigma}}{4},\frac{n-m_{\sigma}}{4})$ is represented by a
$\Z/2^{[3]}$--framed immersion  $(h,\Lambda,\zeta_L)$, which admits a $\H_{b \times \bb}$--structure
 of the $\Z/2^{[2]}$--framed immersion  $(g,\Psi,\eta_N)$ in the sense of Definition $\ref{ad-struct}$.
\end{lemma}

\begin{lemma}\label{Jaa}
Assume that the
 $\Z/2^{[3]}$--framed immersion
$(g,\Psi,\eta_N)$ represents an element $z \in
Imm^{\Z/2^{[3]}}(n-\frac{n-n_\sigma}{4}, \frac{n-n_\sigma}{4})$, $n \ge 254$,
see $(\ref{check})$, and a reduction of the characteristic mapping $\eta_N$ by 
the mapping $(\ref{1etaad})$ is given, such that Conditions 1,2, from Definition  $\ref{aa-struct}$
are satisfied.

Then the element $\delta^{\Z/2^{[4]}}_{k}(z)$ in the group
$Imm^{\Z/2^{[4]}}(n-\frac{n-m_{\sigma}}{2},\frac{n-m_{\sigma}}{2})$ is represented by a
$\Z/2^{[4]}$--framed immersion  $(h,\Lambda,\zeta_L)$, which admits a $\J_a \times \JJ_a$--structure
 of the $\Z/2^{[3]}$--framed immersion  $(g,\Psi,\eta_N)$  in the sense of Definition $\ref{aa-struct}$.
\end{lemma}

\begin{corollary}\label{cor13}
Assume that the assumptions and the dimensional restriction of Theorem $\ref{th6}$ is satisfied.
Then the element $\delta_k^{\Z/2^{[4]}} \circ \delta_k^{\Z/2^{[3]}} \circ
\delta_k^{\Z/2^{[2]}}(\alpha)$, which is defined by means of the composition 
of the homomorphisms $(\ref{6})$, $k=\frac{n-m_{\sigma}}{16}$, is represented by
a $\Z/2^{[4]}$--framed immersion  $(h,\Lambda,
\zeta_L)$, which admits 
 a bicyclic structure ($\J_a \times \JJ_a$--structure) in the sense of Definition $\ref{aa-struct}$. 
 
Moreover, the projection  of the element
\begin{eqnarray}\label{Lzetalok}
\zeta^{loc}_{\J_a \times \JJ_a,\ast}([\zeta_{L_{\J_a \times \JJ_a}}^{k}]^{op}) \in H_{m_{\sigma}}(K(\J_a \times \JJ_a) \int_{\chi^{[4]}} \Z ,1);\Z/2[\Z/2]))
\end{eqnarray}
into the direct factor  $(\ref{DZ2lokaa})$, $i = m_{\sigma}$, after the expansion over the standard base,
involves not more then the only basic element 
 $t_{a,i} \otimes t_{\aa,i}$, see. $(\ref{tataa})$, $i=\frac{m_{\sigma}}{2}=\frac{n-16k}{2}$. 
 The coefficient at this basic element coincides with the Arf-invariant
 $(\ref{44})$, which is calculated for the  $\Z/2^{[2]}$--framed immersion $(g,\eta_N,\Psi)$.
\end{corollary}

\subsubsection*{Proof of Corollary  $\ref{cor13}$}
By Theorem
 $\ref{th6}$ without loss of a generality, we may assume that the element  $y=\delta_k^{\Z/2^{[2]}}(x) \in
Imm^{\Z/2^{[2]}}(n-\frac{n-m_{\sigma}}{8},
\frac{n-m_{\sigma}}{8})$ is represented by a  $\Z/2^{[2]}$--framed immersion 
$(g,\eta_N,\Psi)$, such that the self-intersection manifold of this immersion contains a closed marked component
$N_{b \times \bb}^{n-8k}$, like in the formula $(\ref{abelcomp})$, and the mapping 
$(\ref{redu})$ on this marked component is well-defined. Then the both conditions in Lemma $\ref{dd}$
are satisfied, therefore Conditions 
1 and 2 from Definition  $\ref{ad-struct}$ are satisfied. 

By Lemma $\ref{Hbb}$ we may assume that the element $y$ is represented by a
$\Z/2^{[2]}$--framed immersion, which admits a  $\H_{b \times \bb}$--structure.
An immersion, which represents the element  $z = \delta_k^{\Z/2^{[3]}}(y)$ contains a marked component 
$N^{n-4k}_{b \times \bb}$. By Lemma 
$\ref{Jaa}$ a  $\J_a \times \JJ_a$--structure of an immersion, which is represent the element  $z$ is well-defined.  
Properties  1 and 2 in Definition  $\ref{aa-struct}$ follow from Lemma  $\ref{ad}$. 
 
Let us consider the image of the element
 $(\ref{Lzetalok})$ from the group  $H_{m_{\sigma}}(K(\J_a \times \JJ_a) \int_{\chi^{[4]}} \Z ,1);\Z/2[\Z/2])$
by the composition of the following two transfers: 
$$ H_{m_{\sigma}}(K(\J_a \times \JJ_a) \int_{\chi^{[4]}} \Z ,1);\Z/2[\Z/2]))  \to $$
$$  H_{m_{\sigma}}(K(\I_{b \times \bb} \int_{\chi^{[2]}} \Z ,1);\Z/2[\Z/2])),$$
which are induced by the inclusions  $(\ref{rast3})$, $(\ref{rast})$.
Properties 1 from Lemmas 
$\ref{ad}$, $\ref{aa}$ imply that the image of the element $(\ref{Lzetalok})$ coincides to the element
$(\ref{etadlok})$ from the subgroup 
$$D_{m_{\sigma}}(\I_{b \times \bb};\Z/2[\Z/2]) \subset H_{m_{\sigma}}(K(\I_{b \times \bb} \int_{\chi^{[2]}} \Z) ,1);\Z/2[\Z/2])).$$ 
The image by the composition of the two transfers of the projection of the element  $(\ref{Lzetalok})$ onto the subgroup  $(\ref{DZ2lokaa})$, $i = m_{\sigma}$, is also coincides with
 $(\ref{etadlok})$. 

Properties of the expansion of the element 
$(\ref{etadlok})$ are given by Lemma $\ref{dd}$, Statement 1.  We get the required properties
for the element 
 $(\ref{Lzetalok})$. 
Corollary
 $\ref{cor13}$ is proved.

\section{$\Q \times \Z/4$--structure (quaternionic-cyclic structure) on $\Z/2^{[4]}$--framed immersion}

Let us recall the definition of the quaternion subgroup $\Q \subset
\Z/2^{[3]}$, which contains the subgroup $\I_a \subset \Q$, see
[A1], section 2, formulas (22),(23),(24).

Let us define subgroups:

\begin{eqnarray} \label{QZ4}
i_{\J_a \times \JJ_a, \Q \times \Z/4}: \J_a \times \JJ_a \subset \Q \times \Z/4,
\end{eqnarray}

\begin{eqnarray} \label{iQaa}
i_{\Q \times \Z/4}: \Q \times \Z/4 \subset
 \Z/2^{[5]},
\end{eqnarray}

\begin{eqnarray} \label{iaaa}
i_{\J_a \times \JJ_a \times \Z/2}: \J_a \times \JJ_a \times \Z/2
\subset
 \Z/2^{[5]}.
\end{eqnarray}

Define the subgroup  $(\ref{QZ4})$. Define the epimorphism on the subgroup $\J_a \times \JJ_a \to \Z/4 \subset \Q$
by the formula
$(x \times y) \mapsto xy$. The kernel of this epimorphism coincides with the antidiagonal subgroup
$\II_a = antidiag(\J_a \times \JJ_a) \subset \J_a \times \JJ_a$,  this subgroup a the direct factor.
This factor is mapped onto the subgroup
 $\Z/4$ by the formula $(x \times x^{-1}) \mapsto x$.
The complement of this factor is the subgroup $\J_a \subset \J_a \times \JJ_a$. The subgroup $(\ref{QZ4})$ is well defined.

To define the subgroups  $(\ref{iQaa})$, $(\ref{iaaa})$ consider the basis
$(\h_{1,+}, \h_{2,+}, \h_{1,-},
\h_{2,-}, \hh_{1,+}, \hh_{2,+}, \hh_{1,-}, \hh_{2,-})$
in the space $\R^8$, determined by $(\ref{h})$, $(\ref{hh})$.

We define an analogous basis  of    $\R^{16}$. This basis consists
of $16$ vectors, this set of the basis vectors is divided two
subsets:

\begin{eqnarray}\label{Qh}
\h_{1,\ast, \ast \ast},\h_{2,\ast, \ast \ast};
\end{eqnarray}

\begin{eqnarray}\label{Qhh}
\hh_{1,\ast, \ast \ast}, \hh_{2,\ast, \ast \ast}.
\end{eqnarray}
where the symbols $\ast$, $\ast \ast$ independently takes the
values  $+,-$.

Let us define the subgroup $(\ref{iQaa})$. The representation
$i_{\Q \times \Z/4}$ is given such that the generator $\j$ of the factor  $\Q
\subset \Q \times \Z/4$ acts  in each
4-dimensional subspace

\begin{eqnarray}\label{11}
diag(Lin(\h_{1,\ast,\ast \ast},\h_{2,\ast, \ast \ast},\h_{1,\ast,-\ast \ast},\h_{2,\ast, -\ast \ast}),
\end{eqnarray}
$$
Lin(\hh_{1,\ast,\ast \ast},\hh_{2,\ast, \ast \ast},\hh_{1,\ast,-\ast \ast},\hh_{2,\ast, -\ast \ast})),
$$

\begin{eqnarray}\label{12}
diag(Lin(\h_{1,-\ast,\ast \ast},\h_{2,-\ast, \ast \ast},\h_{1,-\ast,-\ast \ast},\h_{2,-\ast, -\ast \ast}),
\end{eqnarray}
$$
Lin(\hh_{1,-\ast,\ast \ast},\hh_{2,-\ast, \ast \ast},\hh_{1,-\ast,-\ast \ast},\hh_{2,-\ast, -\ast \ast})),
$$

\begin{eqnarray}\label{13}
antidiag(Lin(\h_{1,\ast,\ast \ast},\h_{2,\ast, \ast \ast},\h_{1,\ast,-\ast \ast},\h_{2,\ast, -\ast \ast}),
\end{eqnarray}
$$
Lin(\hh_{1,\ast,\ast \ast},\hh_{2,\ast, \ast \ast},\hh_{1,\ast,-\ast \ast},\hh_{2,\ast, -\ast \ast})),
$$

\begin{eqnarray}\label{14}
antidiag(Lin(\h_{1,-\ast,\ast \ast},\h_{2,-\ast, \ast \ast},\h_{1,-\ast,-\ast \ast},\h_{2,-\ast, -\ast \ast}),
\end{eqnarray}
$$
Lin(\hh_{1,-\ast,\ast \ast},\hh_{2,-\ast, \ast \ast},\hh_{1,-\ast,-\ast \ast},\hh_{2,-\ast, -\ast \ast}))
$$
by the standard transformations, given by the matrix
(23),  [A1].

Let us note that each 4-dimensional space, described above, corresponds to one of the two subspaces
$\R^2_{a,\ast}$, or to one of the two subspaces $\R^2_{\aa,\ast}$, the definition of this subspaces
is given below the formulas
$(\ref{h})$, $(\ref{hh})$. The generator
 $\i \in \Q$ is represented in the direct sum of the two copies of the corresponding spaces, according to the representation of the  generator of the group
  $\J_a$ given by the matrix (23) [A1]. The generator of the second factor $\Z/4 \subset \Q \times \Z/4$
 $\i \in \Q$ is represented in the direct sum of the two copies of the corresponding spaces, according to the representation of the
 generator of the subgroup
 $\II_a \equiv antidiag(\I_a \times \II_a) \subset \I_a \times \II_a$. The representation
$(\ref{iQaa})$ is well defined.

Let us define the representation $(\ref{iaaa})$ as follows. The
factor  $\J_a \times \JJ_a \subset \J_a \times \JJ_a \times
\Z/2$
is represented in each 4-dimensional subspace $(\ref{11})$-$(\ref{14})$
by the formula $(\ref{iaa})$, this formula is applied in each
subspace with the prescribed basis. The factor  $\Z/2
\subset \J_a \times \JJ_a \times \Z/2$ is represented:

--in 8-dimensional subspace, the direct sum of the subspaces $(\ref{11})$, $(\ref{13})$ by the identity;

--in 8-dimensional subspace, the direct sum of the subspaces $(\ref{12})$, $(\ref{14})$  by the central symmetry;

 the representation
$(\ref{iaaa})$ is well-defined.

We define an order 4 automorphism $\chi^{[5]}$ of the subgroup $\Q \times \Z/4$.
The restriction of this automorphism to the subgroup  $(\ref{QZ4})$
coincides with the automorphism  $\chi^{[4]}$. The extension of the  automorphism $\chi^{[4]}$ on the subgroup to the automorphism $\chi^{[5]}$ on the group is defined by
the identity on the generator
$\j$. It is easy to verify that the  automorphism described above is uniquely well defined.

Consider the homomorphism
\begin{eqnarray}\label{prQ}
p_{\Q}: \Q \times \Z/4 \to \Q,
\end{eqnarray} 
which is the projection on the first factor. The kernel of the homomorphism
$p_{\Q}$ coincides with the image of the antidiagonal subgroup $\II_a \subset \J_a \times \JJ_a$
by the inclusion $(\ref{QZ4})$.  Evidently, the following equation is satisfied: 
\begin{eqnarray}\label{chi5com}
\chi^{[5]} \circ p_{\Q}  = p_{\Q}.
\end{eqnarray} 

Analogously, define the  automorphism (involution) $\chi^{[5]}$ of the group 
$\J_a \times \JJ_a \times \Z/2$ (we denote this new automorphism the same). 
Define the homomorphism 
\begin{eqnarray}\label{prZ4Z2}
p_{\Z/4 \times \Z/2}: \J_a \times \JJ_a \times \Z/2 \to \Z/4 \times \Z/2
\end{eqnarray} 
with the kernel
$\II_a \subset \J_a \times \JJ_a \times \Z/2$. Evidently, the following equation is satisfied:
$\chi^{[5]} \circ p_{\Z/4 \times \Z/2}  = p_{\Z/4 \times \Z/2}$.

Define the following groups
\begin{eqnarray}\label{chi5Z}
(\Q \times \Z/4) \int_{\chi^{[5]}} \Z,
\end{eqnarray}
\begin{eqnarray}\label{chi55Z}
(\J_a \times \JJ_a \times \Z/2) \int_{\chi^{[5]}} \Z,
\end{eqnarray}
as semi-direct products of the corresponding groups with automorphisms and the cyclic group $\Z$.
(see the analogous definitions $(\ref{chi2Z})$, $(\ref{chi3Z})$, $(\ref{chi4Z})$).

Let us define the epimorphism
\begin{eqnarray}\label{pchi5Z}
\omega^{[5]}: (\Q \times \Z/4) \int_{\chi^{[5]}} \Z \to \Q,
\end{eqnarray}
the restriction of this epimorphism on the subgroup
 $(\ref{QZ4})$ coincides with the epimorphism $(\ref{prQ})$. 
It is sufficient to use the formula  $(\ref{chi5com})$, define  $z \in Ker(p_{\Q})$, 
where $z \in \Z$ is the generator.  

Evidently, the following epimorphism
\begin{eqnarray}\label{pchi55Z}
\omega^{[5]}: (\J_a \times \JJ_a \times \Z/2) \int_{\chi^{[5]}} \Z \to \Z/4 \times \Z/2,
\end{eqnarray}
is well-defined, we denote this epimorphism the same way like the epimorphism
$(\ref{pchi5Z})$.

Define the automorphism (involution) of the group  $\Z/2^{[5]}$, which is also denoted by  $\chi^{[5]}$.
In the standard basis of the subspaces
 $(\ref{11})$-$(\ref{14})$ the automorphism $\chi^{[5]}$ is given by the same formulas as the automorphism  $\chi^{[4]}$,
each of this space is an invariant space of
 $\chi^{[5]}$. From the definition it is easy to verify that
 $\chi^{[5]}$ commutes with $(\ref{iQaa})$, $(\ref{iaaa})$.

 Moreover, the following homomorphisms are well defined:
\begin{eqnarray}\label{chi5ZQ}
\Phi^{[5]}: (\Q \times \Z/4) \int_{\chi^{[5]}} \Z \to \Z/2^{[5]},
\end{eqnarray}
\begin{eqnarray}\label{chi5ZZ}
\Phi^{[5]}: (\J_a \times \JJ_a \times \Z/2) \int_{\chi^{[5]}} \Z \to \Z/2^{[5]}.
\end{eqnarray}
These homomorphisms are analogous to the homomorphism
$(\ref{Phi4})$ and are included into the following commutative diagrams  $(\ref{a,Qaa})$, $(\ref{a,aaa})$,
(see the analogous diagram
 $(\ref{a,aa})$):

\begin{eqnarray}\label{a,Qaa}
\begin{array}{ccc}
\qquad (\J_a \times \JJ_a) \int_{\chi^{[4]}} \Z & \stackrel{\Phi^{[4]} \times \Phi^{[4]})}
{\longrightarrow} & \qquad \Z/2^{[4]} \times \Z/2^{[4]} \\
i_{\J_a \times \JJ_a, \Q \times \Z/4} \downarrow \qquad & &  i_{[5]} \downarrow \\
\qquad (\Q \times \Z/4) \int_{\chi^{[5]}} \Z &  \stackrel {\Phi^{[5]}}{\longrightarrow}& \qquad \Z/2^{[5]}, \\
\end{array}
\end{eqnarray}

In this diagram the left vertical homomorphism
$$i_{\J_a \times \JJ_a, \Q \times \Z/4}:
(\J_a \times \JJ_a) \int_{\chi^{[4]}} \Z \to (\Q \times \Z/4) \int_{\chi^{[5]}} \Z$$
is induced by the homomorphism
 $(\ref{QZ4})$, the right vertical homomorphism
$$i_{[5]}: \Z/2^{[4]} \times \Z/2^{[4]} \subset \Z/2^{[5]}$$
is defined by the formula $(\ref{9})$.

\begin{eqnarray}\label{a,aaa}
\begin{array}{ccc}
\qquad (\J_a \times \JJ_a) \int_{\chi^{[4]}} \Z & \stackrel{\Phi^{[4]} \times \Phi^{[4]}}
{\longrightarrow} & \qquad \Z/2^{[4]} \times \Z/2^{[4]} \\
i_{\J_a \times \JJ_a, \J_a \times \JJ_a \times \Z/2} \downarrow \qquad & &  i_{[5]} \downarrow \\
\qquad (\J_a \times \JJ_a \times \Z/2) \int_{\chi^{[5]}} \Z &  \stackrel {\Phi^{[5]}}{\longrightarrow}& \qquad \Z/2^{[5]}, \\
\end{array}
\end{eqnarray}
In this diagram the left vertical homomorphism
$$i_{\J_a \times \JJ_a, \J_a \times \JJ_a \times \Z/2}: (\J_a \times \JJ_a) \int_{\chi^{[4]}} \Z \to
(\J_a \times \JJ_a \times \Z/2) \int_{\chi^{[5]}} \Z$$
is induced by the inclusion of the factor.

 The
following definition is analogous to Definition $\ref{[3,4]}$
(cf Definition 15 of [A1]).

\begin{definition}\label{[5]}
Let a  $\Z/2^{[5]}$--framed  immersion $(h,\Lambda,\zeta_L)$, $h:
L^{n-16k} \looparrowright \R^n$
 represent an element
$z \in Imm^{\Z/2^{[5]}}(n-16k,16k)$. We say that this
$\Z/2^{[5]}$--framed  immersion is an $\Q \times \Z/4$--framed
 immersion if 
 the following two conditions are satisfied.

 --1.  The structure mapping
$\zeta_L: L^{n-16k} \to K(\Z/2^{[5]},1)$ is represented as a
composition of a mapping $\zeta_{\Q
\times \Z/4}: L^{n-16k} \to K((\Q \times \Z/4) \int_{\chi^{[5]}} \Z),1)$
 and the mapping $\Phi^{[5]}: K((\Q \times
\Z/4) \int_{\chi^{[5]}} \Z),1) \to K(\Z/2^{[5]},1)$,
this mapping is induced by the homomorphism $(\ref{chi5ZQ})$.

--2. The mapping  $\bar{\zeta}_{\Q \times \Z/4}: \bar{L}^{n-16k} \to K(\I_{b \times \bb} \int_{\chi^{[2]}} \Z,1)$,
which is defined as the $8$-sheeted covering over the mapping
$\zeta_L: L^{n-4k} \to K(\Z/2^{[5]},1)$, which (by Condition 1) satisfies Condition 1
of Definition 
 $\ref{ab}$, satisfies also Condition 2 of Definition  $\ref{ab}$.
\end{definition}

The cohomology group
 $H^{16}(K((\Q \times \Z/4) \int_{\chi^{[5]}} \Z ,1);\Z/2)$ contains an element $\tau_{\Q \times \Z/4}$,
this element is defined by the following equation
($\ref{ii5}$). Consider the mapping
$\Phi^{[5]}: K((\Q \times \Z/4) \int_{\chi^{[5]}} \Z),1) \to K(\Z/2^{[5]},1)$
and take the pull-back
 $(\Phi^{[5]})^{\ast}(\tau_{[5]})$ of the Euler class
 $\tau_{[5]} \in
H^{16}(K(\Z/2^{[5]}),1);\Z/2)$  of the universal bundle. 
Define
\begin{eqnarray}\label{ii5}
(\Phi^{[5]})^{\ast} (\tau_{[5]}) = \tau_{\Q \times \Z/4} \in H^{16}(K((\Q \times
 \Z/4) \int_{\chi^{[5]}} \Z),1);\Z/2).
\end{eqnarray}

Let us assume that the manifold $L^{n-16k}$ is the
self-intersection manifold of a  $\Z/2^{[4]}$--framed immersion
$(h,\Lambda,\zeta_L)$, and the immersion of this manifold into
$\R^n$ is a $\Z/2^{[5]}$--framed immersion which is  a $\Q \times
\Z/4$--framed immersion.

The mapping $\bar \zeta_{b \times \bb}$ is defined as the
8-sheeted covering over the mapping $\zeta_{L}$ with respect to
the subgroup $i_{\I_{b \times \bb},\Q \times \JJ_a}: \I_{b \times \bb} \int_{\chi^{[2]}} \Z \subset
(\Q \times \JJ_a) \int_{\chi^{[5]}} \Z $. Over the  preimage of this mapping, i.e.
over the manifold  $L^{n-16k}$, the corresponding 8-sheeted
covering  $p_{b
\times \bb,\Q \times \JJ_a}$ is well-defined.
Let us re-denote the characteristic class  $\bar \zeta_{[2],L}$ by
 $\bar \zeta_{b \times \bb} \in H^2(\bar
L_{b \times \bb}^{n-16k};\Z/2)$.

\subsubsection*{Quaternionic-cyclic structure}

Let a $\Z/2^{[4]}$--framed immersion $(g,\Psi,\eta_N)$, $g:
N^{n-8k} \looparrowright \R^n$ represent an element $y \in
Imm^{\Z/2^{[4]}}(n-8k,8k)$, assuming $n > 16k$.
Additionally, let us assume that there exists a mapping
\begin{eqnarray}\label{n-8k}
\eta_{a \times \aa}: N^{n-8k}_{a \times \aa} \to K((\J_a \times \JJ_a) \int_{\chi^{[4]}} \Z),1).
\end{eqnarray}

Let a $\Z/2^{[5]}$--framed immersion $(h,\Lambda,\zeta_L)$, $h:
L^{n-16k} \looparrowright \R^n$, is the self-intersection manifold of the restriction of the immersion
$g$ on the component  $N^{n-8k}_{a \times \aa}$ 
(denote this restriction by  $g_{a \times \aa}$. 
Note that the
following equation $n-16k=m_{\sigma}$, where $m_{\sigma}$ is
defined by the equation $(\ref{check})$, is satisfied.

Let us assume that the manifold $L^{n-16k}_{a \times \aa}$ is represented by the
following disjoint union of two components:
\begin{eqnarray}\label{compon}
L^{n-16k}_{a \times \aa} = L^{n-16k}_{\Q \times \Z/4} \cup L^{n-16k}_{\I_a
\times \II_a \times \Z/2}.
\end{eqnarray}

Let as assume that the following mapping
\begin{eqnarray}\label{LQII}
\lambda = \zeta_{\Q \times \Z/4} \cup \zeta_{\J_a \times \JJ_a
\times \Z/2} : L^{n-16k}_{\Q \times \Z/4} \cup L^{n-16k}_{\J_a
\times \JJ_a \times \Z/2} \to
\end{eqnarray}
$$ K((\Q \times \Z/4) \int_{\chi^{[5]}} \Z),1) \cup K((\J_a \times \JJ_a \times \Z/2) \int_{\chi^{[5]}} \Z),1).$$
is well-defined.

The following definition is analogous to Definition 20 from [A1].

\begin{definition}\label{Q-struct}
Let a $\Z/2^{[4]}$--framed immersion $(g,\Psi,\eta_N)$, equipped with a punctured mapping
$(\ref{n-8k})$ of a marked component $N^{n-8k}_{a \times \aa} \subset N^{n-8k}$, which is defined the reduction of the characteristic mapping $\eta_{[4],N}$ be given. Assume that the Arf-invariant for  $(g,\Psi,\eta_N)$ 
is totally determined by means of the marked component.

 Let us say that a $\Z/2^{[4]}$--framed immersion
$(g,\Psi,\eta_N)$,  admits a $\Q \times \Z/4$--structure  if, respectively to the formula  $(\ref{compon})$,
 the mapping  $(\ref{LQII})$ is well-defined, moreover, the pair of mappings $(\eta_{a \times \aa},\lambda)$
 satisfies the following conditions:

--1. The pair of mappings $(\eta_{a \times \aa,N},\zeta_{\Q
\times \Z/4})$  are related by the following
commutative diagram:

\begin{eqnarray}\label{etaIaIa1}
\begin{array}{ccccc}
\bar{L}^{m_{\sigma}}_{a \times \aa, \Q \times \Z/4} &
\looparrowright & N^{n-8k} &
 \stackrel{\eta_{a \times \aa}}{\longrightarrow} & K((\J_a \times \JJ_a) \int_{\chi^{[4]}} \Z),1)\\
 &&&&\\
 \pi_{a \times \aa, \Q \times \Z/4}  \downarrow  & & & & i_{\J_a \times \JJ_a,\Q \times
\Z/4} \downarrow  \\
&&&&\\
 L^{m_{\sigma}}_{\Q \times \Z/4} & & &
\stackrel{\zeta_{\Q \times \Z/4}}{\longrightarrow} & K((\Q \times
\Z/4) \int_{\chi^{[5]}} \Z),1),
\end{array}
\end{eqnarray}
where the right vertical mapping $\pi_{a \times \aa, \Q
\times \Z/4}: \bar{L}^{m_{\sigma}}_{a \times \aa,\Q \times \Z/4} \to
L^{n-16k}_{\Q \times \Z/4}$ is the canonical 2-sheeted
covering over the component of self-intersection manifold of the
immersion $g_{a \times \aa}$ in the formula $(\ref{compon})$, the left vertical
mapping is induced by the corresponding subgroup.

--2. The pair of mappings $(\eta_{a \times \aa},\zeta_{L_{\J_a
\times \JJ_a \times \Z/2}})$  are related by
the following commutative diagram:

\begin{eqnarray}\label{etaIaIa2}
\begin{array}{ccccc}
\bar{L}^{m_{\sigma}}_{a \times \aa,\J_a \times \JJ_a \times
\Z/2} & \looparrowright & N^{n-8k} & \stackrel{\eta_{\J_a \times
\JJ_a,N}}{\longrightarrow} & K((\J_a \times \JJ_a) \int_{\chi^{[4]}} \Z),1)
\\
&&&&\\
 \pi_{a \times \aa,\J_a \times \JJ_a \times \Z/2}  \downarrow & & & & i_{\J_a \times \JJ_a,\J_a \times
\JJ_a \times \Z/2} \downarrow \\
&&&&\\
L^{m_{\sigma}}_{\J_a \times \JJ_a \times \Z/2} & & &
\stackrel{\zeta_{\J_a \times \JJ_a \times \Z/2}}{\longrightarrow}
& K((\J_a \times \JJ_a \times \Z/2) \int_{\chi^{[5]}} \Z),1),
\end{array}
\end{eqnarray}
where the left vertical mapping $\pi_{a \times \aa,
\J_a \times \JJ_a \times \Z/2 }: \bar{L}^{n-16k}_{a \times \aa,\J_a \times
\JJ_a \times \Z/2} \to L^{n-16k}_{\J_a \times \JJ_a \times \Z/2}$  is the canonical 2-sheeted covering over the
component of the self-intersection manifold of the immersion $g_{a \times \aa}$
in the formula $(\ref{compon})$, the left vertical mapping is
induced by the corresponding subgroup.
\end{definition}

\begin{example}\label{QIIa}
Let a $\Z/2^{[4]}$--framed immersion $(g,\Psi,\eta_N)$, $g:
N^{n-8k} \looparrowright \R^n$ be represented an element $y \in
Imm^{\Z/2^{[4]}}(n-8k,8k)$ and be an $\J_a \times \JJ_a$--framed
 immersion, (see Definition $\ref{[3,4]}$), where $n
> 16k$.

Let a $\Z/2^{[5]}$--framed  immersion $(h,\Lambda,\zeta_L)$, $h:
L^{n-16k} \looparrowright \R^n$, which is the immersion of the self-intersection manifold of 
$(g,\Psi,\eta_N)$ be represented the element $z =
\delta^{\Z/2^{[5]},k}(y) \in Imm^{\Z/2^{[5]}}(n-16k,16k)$  and be
an $\Q \times \Z/4$--framed  immersion (see Definition
$\ref{[5]}$). Assume that the restriction of the reduction mapping of $\zeta_L$ and the the restriction
of the mapping $\eta_{a \times \aa}$ to the canonical 2-sheeted covering $\bar{L}^{n-16k}$ of $L^{n-16k}$
are homotopic.

 Then the $\Z/2^{[4]}$--framed  immersion
$(g,\eta_N,\Psi)$, which is equipped with the mapping 
$\eta_{a \times \aa}$,
 admits an $\Q \times \Z/4$--structure, given by the
reduction $\zeta_{\Q \times \Z/4}$ of the classifying mapping
$\zeta_{[5],L}$. The manifold in the decomposition
$(\ref{compon})$ contains the empty second component (cf
Example 21 of [A1].
\end{example}
\[  \]

\subsubsection*{Justification of the example $\ref{QIIa}$}
 Consider the manifold   $L^{m_{\sigma}}$ of self-intersection points of the immersion
 $g$, given by the formula
 $(\ref{compon})$, and define the decomposition of this manifold such that $L^{m_{\sigma}}$
coincides  with the first component $L^{m_{\sigma}}_{\Q \times
\Z/4}$, i.e. the second component $L^{m_{\sigma}}_{\J_a \times
\JJ_a \times \Z/2}$ is empty.
The commutativity of the diagram
$(\ref{etaIaIa1})$ follows from the diagram  $(\ref{a,Qaa})$ and
from the definition $\ref{[5]}$, the diagram $(\ref{etaIaIa2})$ is represented
by the empty manifold.

\bigskip \bigskip
The following theorem is analogous to  Lemmas $\ref{Hbb}$, $\ref{Jaa}$.

\begin{theorem}$\label{th13}$
Assume that the dimensional restriction is $n=2^{\ell}-2$, $\ell \ge 11$, $\sigma \ge 5$ (in Theorem $\ref{th6}$
 is required  $\ell \ge 8$). Assume that a 
 $\Z/2^{[4]}$--framed immersion
$(g,\Psi,\eta_N)$ represents an element $y \in
Imm^{\Z/2^{[4]}}(n-\frac{n-n_\sigma}{2}, \frac{n-n_\sigma}{2})$,
 and a mapping $(\ref{1zetaaa})$, which determines a $(\J_a \times \JJ_a) \int_{\chi^{[4]}} \Z$--reduction
 of the characteristic mapping $\eta_N$ is given.

 Then the element
$\delta^{\Z/2^{[5]}}_{k}(y)$ in the group 
$Imm^{\Z/2^{[5]}}(n-\frac{n-m_{\sigma}}{2},\frac{n-m_{\sigma}}{2})$ is represented by a
$\Z/2^{[5]}$--framed immersion  $(h,\zeta_L,\Lambda)$ 
and using this immersion a
 $\Q \times \Z/4$--structure of the $\Z/2^{[4]}$--framed immersion
$(g,\eta_N,\Psi)$ is defined.
\end{theorem}

The following corollary is analogous to the Corollary
 $\ref{cor13}$.

\begin{corollary}\label{cor19}

Assume that the hypotheses of Theorem $\ref{th6}$ hold under the
stronger dimensional restrictions from Theorem $\ref{th13}$. Then
for an arbitrary $x$ the element
\begin{eqnarray} \label{ddddx}
\delta_k^{\Z/2^{[5]}} \circ \delta_k^{\Z/2^{[4]}} \circ \delta_k^{\Z/2^{[3]}} \circ
\delta_k^{\Z/2^{[2]}}(x),
\end{eqnarray}
defined by the composition of homomorphisms $(\ref{6})$,
$k=\frac{n-n_\sigma}{16}$, is represented by a $\Z/2^{[5]}$--framed
immersion $(h,\Lambda,
\zeta_L)$  and using this immersion   
a $\Q \times \Z/4$--structure of the element $\delta_k^{\Z/2^{[4]}} \circ \delta_k^{\Z/2^{[3]}} \circ
\delta_k^{\Z/2^{[2]}}(x)$ is defined.

With respect to  the decomposition
 $(\ref{compon})$ the following homology classes are well-defined: 
\begin{eqnarray}\label{LzetaQlok}
\bar{\zeta}^{loc}_{\ast}([\bar{L}_{a \times \aa,\Q \times \Z/4}]) \in H_{m_{\sigma}}(K((\J_a \times \JJ_a) \int_{\chi^{[4]}} \Z) ,1);\Z/2[\Z/2])),
\end{eqnarray}
\begin{eqnarray}\label{LzetaQQlok}
\bar{\zeta}^{loc}_{\ast}([\bar{L}_{a \times\aa,\J_a \times \JJ_a \times \Z/2}]) \in H_{m_{\sigma}}(K((\J_a \times \JJ_a) \int_{\chi^{[4]}} \Z) ,1);\Z/2[\Z/2])).
\end{eqnarray}
  The projection of the element
 $(\ref{LzetaQlok})$ to the factor $(\ref{DZlokaa})$ by the homomorphism 
$(\ref{deltaDaa})$
is expended 
over the standard basis. This expansion contains not more then one nontrivial element, which is defined by the coefficient of the monomial
$t_{a,i} \otimes t_{\aa,i}$, see. $(\ref{tataa})$, $i=\frac{m_{\sigma}}{2}=\frac{n-16k}{2}$. 
This coefficient coincides with the characteristic number
$(\ref{44})$, which is calculated for the  $\Z/2^{[2]}$--framed immersion  $(g,\Psi,\eta_N)$.
The element  $(\ref{LzetaQQlok})$ is trivial.
\end{corollary}

\subsubsection*{Proof of Corollary $\ref{cor19}$}
Let us consider the immersion 
 $g_{a \times \aa}$, which is defined as the restriction of the immersion $g$ on the marked component
$N^{n-8k}_{a \times \aa} \subset N^{n-8k}$ (Recall, that for the marked component a reduction 
of the characteristic mapping is given, and the Arf-invariant of the cobordism class $x$ is determined
by the cobordism class of the marked component.)
Let us apply to the immersion 
 $g_{a \times \aa}$ the Herbert's theorem with  $\Z/2[\Z/2]$--local coefficients system, see the analogous construction in the proof of Lemma $\ref{ad}$.

From this theorem we get that the sum of the homology classes
$(\ref{LzetaQlok})$, $(\ref{LzetaQQlok})$ is equal to the homology class $(\ref{Lzetalok})$. In Corollary $\ref{cor13}$
the required property for the class $(\ref{Lzetalok})$ is proved. 
 It is sufficiently to prove that the homology class
 $(\ref{LzetaQQlok})$ is trivial. Consider the projection $\J_a \times \JJ_a \times \Z/2 \to \J_a \times \JJ_a$,
which is induced the projection $(\J_a \times \JJ_a \times \Z/2) \int_{\chi^{[5]}} \Z \to (\J_a \times \JJ_a) \int_{\chi^{[4]}} \Z $. Consider the corresponding mapping 
$$p: K((\J_a \times \JJ_a \times \Z/2) \int_{\chi^{[5]}} \Z,1) \to
K((\J_a \times \JJ_a) \int_{\chi^{[4]}} \Z,1). $$

The element
 $(\ref{LzetaQQlok})$ is obtained from the element 
\begin{eqnarray}\label{pLzetaQQlok}
p_{\ast} \circ \zeta^{loc}_{\ast}([L_{\J_a \times \JJ_a \times \Z/2}]) \in H_{m_{\sigma}}(K((\J_a \times \JJ_a) \int_{\chi^{[4]}} \Z) ,1);\Z/2[\Z/2])).
\end{eqnarray}
by means of the composition with the canonical 2-sheeted covering. Therefore the element
$(\ref{LzetaQQlok})$ is trivial. Corollary $\ref{cor19}$ is proved.

\section{Solution of the Kervaire invariant problem}

In this section we prove the following result.

\subsubsection*{Main Theorem}

There exists a natural number $l_0$ such that for an arbitrary
integer $\ell \ge l_0$ the Kervaire
invariant defined by formula ($\ref{1}$) is trivial. 
(Recall, $n=2^{\ell}-2$,  $\ell \ge 12$, $m_{\sigma}=2^{\sigma}-2$, 
where $\sigma$ is defined by the
equation $(\ref{sigma})$. We may put $m_{\sigma}=30$, $\sigma=5$) 
Assume that an element
 $x \in
Imm^{sf}(n-1,1)$ admits a compression of the order $m_{\sigma}+2$, then $\Theta^{sf}(x)=0$.
\[  \]

\subsubsection*{Remark}
The Main Theorem  could be proved  under weaker dimensional
assumption then $n \ge 4094$ by the considered approach. See a remark in the
Introduction of the part $III$ [A3] of the paper. 
\[  \]

\begin{theorem}{\textbf{Compression Theorem}}\label{comp}
For an arbitrary positive integer $d$ there exists a positive
integer $\ell=\ell(d)$ such that for an arbitrary element in the $2$-component of the
cobordism group  $Imm^{sf}(2^{l'}-3,1)$, assuming $l' \ge \ell$,
admits a compression of the order  $(d-1)$ $($ see Definition
$\ref{5}$ $)$.
\end{theorem}

\subsubsection*{Remark}
The proof of Theorem $\ref{comp}$ is presented in Section 7.
By the Pontryagin-Thom  construction (in the form by Welles) the cobordism group $Imm^{sf}(2^{l}-3,1)$
is isomorphic to the stable homotopy group  $\Pi_{2^{l}-2}(K(\Z/2,1))$. The space  $Q(K(\Z/2,1))$
is 2-primary. This implies that the cobordism group $Imm^{sf}(2^{l}-3,1)$ has no odd torsions.
From Theorem  $\ref{comp}$ an explicit sub-exponential estimation
of the dimension $2^{\ell}-2$, $\ell \ge
l_0=l_0(d)$, for which an arbitrary element in $Imm^{sf}(2^{l}-3,1)$ admits a $d-1$-compression, could be possible.  
Prof. D. Ravenel in [R]
gave an explicit formula for  $l_0(d)$
for small $d$. 
\[  \]

To prove the Main Theorem is sufficiently to prove that the residue class in the cobordism group $Imm^{sf}(2^{l'}-3,1)$
which is determined by the non-trivial Arf-invariant, is the empty, or, contains an element 
 $x$, which admits a compression of the order  $16$. The condition of a compression of the order  $q$, $3q < 2^{l'}-2$ of an element $x$ is equivalent to the following condition: the adjoined element 
$y \in \Pi_{2^{l'-2}}$ to the element $x$ admits a $q+1$--desuspension in the unstable domain, i.e. the stable homotopy class $y$ is represented with sphere-of-origin
$2^{l'}-2-q$. Accordingly to the result  [R-S], if the group $\Pi_{126}$ contains an element with Kervaire-invariant 1,
it will have a representative with sphere-of-origin $116$.
\[  \]

\subsubsection*{Proof of Main Theorem from Corollary $\ref{cor19}$}
The proof is analogous to the proof of Proposition  29 [A1].
Compute a positive integer $k$ from the equation $n-16k =
m_{\sigma}$, $k \ge 7$, $k \equiv 0 \pmod{2}$, this is possible if
the condition  $\ell \ge 9$ is satisfied. By Theorem $(\ref{th13})$ we have $\ell \ge 12$. 
Consider a triple $(f:
M^{n-k} \looparrowright \R^n,\kappa,\Xi)$, representing the given
element  $x \in Imm^{sf}(n-k,k)$.

Consider the element  $\delta_k^{\Z/2^{[4]}} \circ
\delta_k^{\Z/2^{[3]}} \circ \delta_k^{\Z/2^{[2]}}(x) \in
Imm^{\Z/2^{[4]}}(n-8k,8k)$, see the formula $(\ref{ddddx})$,
represented by a $\Z/2^{[4]}$--framed immersion $(g,\eta_N,\Psi)$.

Denote by $L^{n-16k}$ the self-intersection manifold of the
immersion  $g$. By Corollary
  $\ref{cor19}$ we may assume that the triple
$(g,\eta_N,\Psi)$ admits a $\Q \times \Z/4$--structure (see
Definition $\ref{Q-struct}$).

Let us assume that the classifying mapping  $\eta_N$ satisfies the
condition of Example  $\ref{QIIa}$. This means that the following
equalities are satisfied:
$$ \eta_N = i_{\J_a \times \JJ_a, \Z/2^{[4]}} \circ \eta_{a \times \aa}, $$
$$ \zeta_L = i_{\Q \times \Z/4, \Z/2^{[5]}} \circ \zeta_{\Q \times \Z/4}. $$
Let us prove the Theorem under this assumption.

Let us denote by $\tilde N^{n-8k-2} \subset N^{n-8k}$ the
submanifold, representing the Euler class of the vector bundle
 $\eta_{a \times \aa}^{\ast}(\psi_{+})$, where by $\psi_{+}$ denotes
 a 2-dimensional vector bundle over the classifying space
$K((\J_a \times \JJ_a) \int_{\chi^{[4]}} \Z,1)$,
 given by the formula
 \begin{eqnarray}\label{psi+}
\psi_{+} = (\omega^{[4]})^{\ast}(\psi_{\Z/4}), 
\end{eqnarray}
 where
$\psi_{\Z/4}$ is the universal 2-bundle over $K(\Z/4,1)$, the mapping
 $\omega^{[4]}$ is induced by the homomorphism $(\ref{om4})$. Because
the classifying map $\eta$ admits a $\J_a \times \JJ_a$--reduction, the submanifold $\tilde
N^{n-2k-2} \subset N^{n-2k}$ is co-oriented
(we do not use that  $k$ is even).

Let us denote by
 $\tilde g: \tilde N^{n-8k-2} \looparrowright \R^n$ the
 restriction of the immersion $g$ to the submanifold  $\tilde N^{n-8k-2} \subset N^{n-8k}$,
assuming that the immersion $\tilde g$ is generic. The normal
bundle of the immersion
 $\tilde g$ is the Whitney sum $\nu_g \oplus \tilde{\nu}_{\J_a}$, where $\nu_g$  denotes
the normal bundle of the immersion $g$, restricted to  $\tilde
N^{n-8k-2} \subset N^{n-8k}$ (this bundle has the structure
group  $(\J_a \times \JJ_a) \int_{\chi^{[4]}} \Z$,
 this group determines a reduction the structure group of $\Z/2^{[4]}$--framing
$\Psi$), by $\tilde{\nu}_{\J_a}$ is denoted
the normal bundle of the submanifold  $\tilde N^{n-8k-2} \subset
N^{n-8k}$ (this bundle is an $(\omega^{[4]})^{\ast}(\Z/4)$-bundle).

Let us denote by  $\tilde L^{n-16k-4}$ the self-intersection
manifold of the immersion  $\tilde g$. The manifold  $\tilde
L^{n-16k-4}$ is a submanifold of the manifold  $L^{n-16k}$,
$\tilde L^{n-16k-4} \subset L^{n-16k}$. The 
immersion
 $\tilde h h \vert_{\tilde L}: \tilde L^{n-16k-4}
\looparrowright \R^n$ is well-defined.

The normal bundle of this immersion  $\tilde h$ is isomorphic to the
Whitney sum $\nu_h \oplus \tilde \nu_{\Q}$, where by $\nu_h$ is
denoted the restriction of the normal bundle of the immersion
 $h$ over the submanifold $\tilde L^{n-16k-4} \subset L^{n-16k}$,
 by  $\tilde \nu_{\Q}$ denotes the normal bundle of the submanifold $\tilde L^{n-16k-4} \subset
L^{n-16k}$. 

Let us repeat the proof of Lemma 7 [A1] (the
commutativity of the left square of the diagram   (8) of [A1]). 
This arguments proves that the submanifold $\tilde L^{n-16k-4} \subset
L^{n-16k}$ represents the Euler class of the bundle
$(\zeta^{\ast}_{\Q \times \Z/4})^{\ast}(\psi_{\Q \times \Z/4})$,
where  $\psi_{\Q \times \Z/4}$  is  the $SO(4)$-bundle over $K((\Q \times \Z/4) \int_{\chi^{[5]}} \Z,1)$, given by
the representation  $(23)-(25)$ of [A1]) as follows:
\begin{eqnarray}\label{psiQZ4}
\psi_{\Q \times \Z/4} = (\omega^{[5]})^{\ast}(\psi_{\Q}),
\end{eqnarray}
the mapping
$$\omega^{[5]}: K((\Q \times \Z/4) \int_{\chi^{[5]}} \Z,1) \to K(\Q,1)$$ 
corresponds to the epimorphism  $(\ref{pchi5Z})$.

 Consider 2-sheeted covering $p_{a \times \aa, \Q \times
\JJ_a}: K((\J_a \times \JJ_a) \int_{\chi^{[4]}} \Z,1) \to K((\Q \times \JJ_a)\int_{\chi^{[5]}} \Z,1)$
  over the universal space
and denote the bundle 
$p_{a \times \aa,\Q \times \Z/4}^{\ast}(\psi_{\Q \times \Z/4})$
by  $ \psi_{\Q \times \Z/4}^!$, where the bundle  $\psi_{\Q \times \Z/4}$
is defined by the formula $(\ref{psiQZ4})$. 

For the universal bundle   $ \psi_{\Q \times \Z/4}^!$  the following formula is
satisfied: 
$$ \psi_{\Q \times \Z/4}^! = \psi_+ \oplus \psi_-, $$
where the bundle $\psi_+$ (this bundle is given by the formula $(\ref{psi+})$)
admits a lift $\psi_{+,U}$ to a complex
$U(1)$--bundle, the bundle $\psi_-$ is an $SO(2)$--bundle, obtained
from $\psi_{+,U}$ by means of the complex conjugation and
forgetting the complex structure. The bundles $\psi_+$, $\psi_-$
satisfy the equation: $e(\psi_+)=-e(\psi_-)$.

Let us denote by $m \in H^{16k}( N^{n-8k};\Z[\Z/2])$
 the cohomology class, with local coefficient system (this cohomology class is defined analogously to the homology class $(\ref{etadd16})$), which is dual to
 the fundamental class of the oriented submanifold
$\bar L^{n-16k} \subset N^{n-8k}$ in the oriented manifold
 $N^{n-8k}$. Let us denote by  $e_g \in H^{16k}(N^{n-8k};\Z[\Z/2])$
the Euler class of the normal bundle $\nu_g$ of the immersion $g$.
By the analog of the Herbert theorem for the immersion $g: N^{n-8k}
\looparrowright \R^n$ with the self-intersection manifold
$L^{n-16k}$ (an analogous theorem was formulated    in Theorem 1.1 of [E-G], the case  $r=1$, but only  with integer coefficients)
 the following formula is satisfied:
\begin{eqnarray}\label{1.501}
e_g+m=0,
\end{eqnarray}
where  the both cohomology classes in this formula are defined by means
of the $\Z[\Z/2]$-local coefficients system.

Let us denote by $\tilde m \in H^{16k-4}( N^{n-8k};\Z[\Z/2])$
 the cohomology class, dual to
 the fundamental class of the oriented submanifold
$\bar{\tilde L}^{n-16k-4} \subset \tilde N^{n-8k-2} \subset
N^{n-8k}$ in the oriented manifold
 $N^{n-8k}$. Let us denote by  $ e_{\tilde g} \in H^{16k-4}(N^{n-8k};\Z[\Z/2])$
the cohomology class, dual to the Euler class of the normal bundle
$\nu_{\tilde g}$ of the immersion $\tilde g$ on the submanifold
$\tilde N^{n-8k} \subset N^{n-8k}$. 

By the analog of the Herbert theorem for
the immersion $\tilde g: \tilde N^{n-8k} \looparrowright \R^n$
with the self-intersection manifold $\tilde L^{n-16k-4}$ 
 the following formula is satisfied:
\begin{eqnarray}\label{1.502}
e_{\tilde g}+ \tilde m=0,
\end{eqnarray}
where the both the cohomology classes in this formula are also defined by means
of the $\Z[\Z/2]$-local coefficients system.

 Because $\bar \lambda = \eta_{a \times \aa}$,
 we may use the equation: $\bar \lambda^{\ast}(\psi_{\Q \times \Z/4}^!) = \eta_{a \times
\aa}^{\ast}(\psi_+) \oplus \eta_{a \times
\aa}^{\ast}(\psi_-)$.
  The following
equation is satisfied: 
$$\tilde m = m \cdot e(\eta_{a \times \aa}^{\ast}(\psi_+))\cdot e(\eta_{a \times \aa}^{\ast}(\psi_-)),$$ 
where the right side is the product of
the three cohomology classes: $m$ and the two Euler classes of the
corresponding bundles. The following equation is satisfied:
$e_{\tilde g} =
e_g \cdot e^2(\eta_{a \times \aa}^{\ast}(\psi_+))$.

The equation $(\ref{1.502})$ can be rewritten in the following
form:
\begin{eqnarray}\label{1.503}
e_g \cdot e^2(\eta_{a \times \aa}^{\ast}(\psi_+)) +
e_g \cdot e(\eta_{a \times \aa}^{\ast}(\psi_+)) \cdot e(\eta_{a \times
\aa}^{\ast}(\psi_-))) =0.
\end{eqnarray}

Then we may take into account $(\ref{1.501})$ and the equation
$e(\eta_{\I_a \times \II_a}^{\ast}(\psi_-))=-e(\eta_{\I_a \times
\II_a}^{\ast}(\psi_+))$. 
%in this equation homology classes are defined with $\Z[\Z/2]$ coefficients. 
Let us rewrite the previous formula as
follows:
\begin{eqnarray}\label{1.504}
2e_g \cdot e^2(\eta_{a \times \aa}^{\ast}(\psi_+))=0.
\end{eqnarray}

Let us prove that the equation
$(\ref{1.504})$ implies that the expansion of the expression $(\ref{LzetaQlok})$
with respect to the standard base does not involves the generator
$t_{a,i} \otimes t_{\aa,i}$, see $(\ref{tataa})$, $i=\frac{m_{\sigma}}{2}=\frac{n-16k}{2}$.
$(\ref{tdtdd})$. The homology class
$[N] \cap e_g$ is equal to the homology class $[\tilde L]$, therefore from the expression $(\ref{1.504})$
follows that the homology class $2\eta_{a \times \aa, \ast}([\tilde L]) $
in the group $H_{m_{\sigma-4}}(K(\J_a \times \JJ_a) \int_{\chi^{[4]}} \Z,1);\Z[\Z/2])$
 is trivial.

The homology class $\eta_{a \times \aa,\ast}([\tilde L])$
is calculated by the formula 
\begin{eqnarray}\label{etiL}
\eta_{a \times \aa,\ast}([\tilde L]) = \eta_{a \times \aa,\ast}([L]) \cap (\omega^{[4],\ast}[\tau])^2,
\end{eqnarray}
where
$\omega^{[4]}$ is defined by the formula  $(\ref{om4})$, and where  $\tau \in H^2(K(\Z/4,1);\Z)$ is the generator.

Let us consider the expansion of the homology class $(\ref{etiL})$ 
%$\eta_{a \times \aa,\ast}([L]) \cap \omega^{[4],\frac{n-16k}{2}}$ 
over the standard basis of the group  $H_{m_{\sigma}}(K((\J_a \times \JJ_a) \int_{\chi^{[4]}} \Z);\Z[\Z/2])$, using  Lemma  $\ref{oaa}$. All basis elements of the subgroup $D(\J_a \times \JJ_a;\Z[\Z/2]) \subset H_{m_{\sigma}}(K((\J_a \times \JJ_a) \int_{\chi^{[4]}} \Z);\Z[\Z/2])$ 
(elements of this subgroups are detected by means of the monomomorphism $(\ref{deltaDaa})$) except, probably, the elements 
$t_{a,\frac{m_{\sigma}}{2}+4} \otimes t_{\aa,\frac{m_{\sigma}}{2}-4}$ and   $t_{a,\frac{m_{\sigma}}{2}-4} \otimes t_{\aa,\frac{m_{\sigma}}{2}+4}$ 
are involved with even coefficients.

In the expansion of the element $(\ref{etiL})$
%$\eta_{a \times \aa,\ast}([L]) \cap \omega^{[4],\frac{n-16k}{2}}$ 
over the standard basis 
the elements  $t_{a,\frac{m_{\sigma}}{2}} \otimes t_{\aa,\frac{m_{\sigma}}{2}-4}$,   $t_{a,\frac{m_{\sigma}}{2}-4} \otimes t_{\aa,\frac{m_{\sigma}}{2}}$
are involved with  coefficients of the same parity as the parity of the corresponding coefficients in the expansion of the homology class
 $(\ref{LzetaQlok})$. By the formula $(\ref{1.504})$ we get that all this two coefficients are even. 

By Corollary $\ref{cor19}$ the characteristic number 
$(\ref{44})$ is trivial.
The theorem  in the particular case
is proved.

Let us prove the theorem in a general case. Let us consider the
pair of mappings  $(\eta_{a \times \aa},
 \lambda)$, where
 $\eta_{a \times \aa}: N^{n-8k}_{a \times \aa} \to K((\J_a \times \JJ_a) \int_{\chi^{[4]}} \Z,1)$,
$\lambda = \zeta_{\Q \times \Z/4} \cup \zeta_{\J_a \times \JJ_a
\times \Z/2} : L^{n-16k}_{\Q \times \Z/4} \cup L^{n-16k}_{\J_a
\times \JJ_a} \to K((\Q \times \Z/4) \int_{\chi^{[5]}} \Z,1) \cup K((\J_a \times \JJ_a
\times \Z/2) \int_{\chi^{[5]}} \Z,1)$ and $L^{n-16k} = L^{n-16k}_{\Q \times \Z/4} \cup
L^{n-16k}_{\J_a \times \JJ_a \times \Z/2}$,
 where these two
mappings are determined by the quaternionic structure of the
$\Z/4^{[4]}$--framed immersion  $(g,\eta_N,\Psi)$.

 Let us consider the manifold $L^{n-16k}_{\Q \times \Z/4} \cup
L^{n-16k}_{\J_a \times \JJ_a \times \Z/2}$, defined by the
formula $(\ref{compon})$, and the manifolds $\bar
L^{n-16k}_{\Q \times \Z/4} \cup \bar L^{n-4k}_{\J_a \times \JJ_a
\times \Z/2}$, which  are the canonical 2-sheeted covering
manifolds over the components of $L^{n-16k}_{a \times \aa}$.

The formula  $(\ref{1.501})$ is satisfied, where the cohomology class
$m$ (this class is dual to the fundamental class  $[\bar L_{a \times \aa}]$ of
the submanifold $\bar L^{n-16k}_{a \times \aa} \subset N^{n-8k}_{a \times \aa}$)  decomposes
into the following sum:
\begin{eqnarray}\label{mm}
m=m_{\Q \times \Z/4} + m_{\J_a \times \JJ_a \times \Z/2},
\end{eqnarray}
correspondingly to the type of the components  of the
self-intersection manifold.

Let us consider the submanifold  $\tilde N^{n-8k-2} \subset
N^{n-8k}$, representing the Euler class of the bundle $\eta_{\J_a
\times \JJ_a}^{\ast}(\psi^+)$. The following immersion $\tilde g:
\tilde N^{n-8k-2} \looparrowright \R^n$ is well defined by the
restriction of the immersion  $g_{a \times \aa}$ to the submanifold $\tilde
N^{n-8k-2} \subset N^{n-8k}_{a \times \aa}$.  Let us denote by  $\tilde
L^{n-16k-4}$ the self-intersection manifold of the immersion
$\tilde g$.

The following inclusion $\tilde L^{n-16k-4} \subset L^{n-16k}_{a \times \aa}$ is
well-defined. In particular, the manifold $\tilde L^{n-16k-4}$ is
representing by the union of the following two components: $\tilde
L^{n-16k-4} = \tilde L^{n-16k-4}_{\Q \times \Z/4} \cup \tilde
L^{n-16k-4}_{\I_a \times \II_a \times \Z/2}$. From Lemma 34 of [A1]
we will prove the following. The submanifold $\tilde
L^{n-16k-4}_{\Q \times \Z/4} \subset L^{n-16k}_{\Q \times \Z/4}$
represents the Euler class of the bundle  $\zeta_{\Q \times \Z/4}^{\ast}(\psi_{\Q \times \Z/4})$.

The submanifold
$\tilde L^{n-2k-4}_{\I_a \times \II_a \times Z/2} \subset
L^{n-2k}_{\I_a \times \II_a \times \Z/2}$ represents the Euler
class of the bundle $\zeta_{\I_a \times \II_a \times
\Z/2}^{\ast}(\psi_{\I_a \times \II_a \times i\Z/2})$, where
$\psi_{\I_a \times \II_a \times \Z/2}$ is the universal 4-bundle
over the space $K((\I_a
\times \II_a \times \Z/2) \int_{\chi^{[5]}} \Z,1)$.
 (comp. with an analogous 
definition in the proof of Theorem 12 [A1]).

The cohomology class $\tilde{m}$ is well-defined
as in the formula $(\ref{1.502})$, moreover the following
formula is satisfied:
\begin{eqnarray}\label{tildemm}
\tilde m=\tilde m_{\Q \times \Z/4} +   \tilde m_{\J_a \times \JJ_a \times \Z/2},
\end{eqnarray}
where the terms in the right side of the formula are defined as
the homology classes, dual to the fundamental classes
$[\bar{\tilde L}_{\Q \times \Z/4}]$, $[\bar{\tilde L}_{\I_a
\times \II_a \times \Z/2}]$  of the canonical coverings over the
corresponding component.

The following formula determines a relation  between the cohomology classes $m_{\Q \times \Z/4}$ and $\tilde m_{\Q
\times \Z/4}$: 
$$\tilde m_{\Q \times \Z/4} = m_{\Q
\times \Z/4} \cdot e(\eta_{\J_a \times \JJ_a}^{\ast}(\psi_+))
\cdot e(\eta_{\J_a \times \JJ_a}^{\ast}(\psi_-)).$$  
The following formula determines a relation  between the cohomology classes
$m_{\J_a \times \JJ_a \times \Z/2}$  and $\tilde m_{\J_a \times
\JJ_a \times \Z/2}$: 
$$\tilde m_{\J_a \times \JJ_a
\times \Z/2} = m_{\J_a \times \JJ_a \times \Z/2} \cdot e^2(\eta_{\J_a
\times \JJ_a}^{\ast}(\psi_+)).$$

To prove the last formula we use the following formula: 
$$\bar{\zeta}_{\J_a \times \JJ_a \times \Z/2}(\psi_{\J_a \times \JJ_a \times \Z/2}^!)=
\eta_{a \times\aa}^{\ast}(\psi_+ \oplus \psi_+).$$
In this formula the mapping $\eta_{a \times \aa}$ is the restriction of the corresponding mapping to the immersed submanifold
$\bar{\tilde L^{n-16k-4}} \subset \tilde{N}^{n-8k-2} \subset N^{n-8k}_{a \times \aa}$,
$\psi_{\J_a \times \JJ_a \times \Z/2}^!$ is the pull-back of the universal $SO(4)$--bundle
over $K((\J_a \times \JJ_a \times \Z/2) \int_{\chi^{[5]}} \Z,1)$ to the 2-sheeted covering
$K((\J_a \times \JJ_a) \int_{\chi^{[4]}} \Z,1)$. Obviously, the bundle
$\psi_{\J_a \times \JJ_a \times \Z/2}^!$
is the Whitney sum of two copies of the  pull-back of the universal
$SO(2)$-bundles over $K(\Z/4,1)$ by the mapping
$K((\J_a \times \JJ_a) \int_{\chi^{[4]}} \Z,1) \to K(\Z/4 \times \Z/2,1) \to K(\Z/4,1)$.

The analog of the formula $(\ref{1.503})$ is the following:
\begin{eqnarray}\label{1.5033}
e_g \cdot e^2(\eta_{a \times \aa}^{\ast}(\psi_+)) - m_{\Q \times
\JJ_a} \cdot e^2(\eta_{a \times \aa}^{\ast}(\psi_+)) +
\end{eqnarray}
$$ m_{\J_a
\times \JJ_a \times \Z/2} \cdot e^2(\eta_{a \times
\aa}^{\ast}(\psi_+)) =0.$$

Let us multiply both sides of the formula $(\ref{tildemm})$ by
the cohomology class $e^2(\eta_{\I_a \times
\II_a}^{\ast}(\psi_+))$ and take the sum with the opposite sign
with $(\ref{1.5033})$, we get:
\begin{eqnarray}\label{1.5044}
2m_{\Q \times \Z/4} \cdot e^2(\eta_{\I_a \times
\II_a}^{\ast}(\psi_+))=0.
\end{eqnarray}
This is an analog of the formula $(\ref{1.504})$.

Let us prove that the Kervaire invariant of  the
$\Z/2^{[4]}$--framed immersion
 $(g,\Psi,\eta)$ is trivial. 
 The expansion of the element $(\ref{LzetaQlok})$ over the standard basis does not involved the monomial
 $t_{a,i} \otimes t_{\aa,i}$, см. $(\ref{tataa})$, $i=\frac{m_{\sigma}}{2}=\frac{n-16k}{2}$.
 The proof is analogous to the previous case using Corollary $\ref{cor19}$ and the formula 
 $(\ref{1.5044})$ instead of the formula $(\ref{1.504})$.
 
The Main Theorem is proved.

\section{Proof of Lemmas $\ref{Hbb}$, $\ref{Jaa}$ and Theorem $\ref{th13}$}

 We shall prove first Lemmas $\ref{Hbb}$ and $\ref{Jaa}$.   Let
us define the positive integer $m_\sigma$ by the formula $(\ref{check})$ 
(in this section to simplify the calculation of dimensions we assume that $\sigma \ge 6$
All the constructions are well-defined in the case $\sigma =5$). 

Let us denote by $ZZ_{\J_a \times \JJ_a}$ the
direct product of the two standard lens spaces  $\pmod{4}$, namely
\begin{eqnarray}\label{ZZIaIa}
  ZZ_{\J_a \times \JJ_a}= S^{n-\frac{n-m_{\sigma}}{8}+9}/\i \times S^{n-\frac{n-m_{\sigma}}{8}+9}/\i.
\end{eqnarray}
Obviously, $$2n > \dim(ZZ_{\J_a \times \JJ_a}) = \frac{7n+m_{\sigma}}{4}+18 > n.$$

On the space $ZZ_{\J_a \times \JJ_a}$ the standard involution $\chi^{[4]}: ZZ_{\J_a \times \JJ_a} \to ZZ_{\J_a \times \JJ_a}$
is given by the formula $\chi^{[4]}(x \times y)=(y \times x)$.

Let us define a polyhedron (a submanifold with singularities)
\begin{eqnarray}\label{Xaa}
 X_{a \times \aa} \subset ZZ_{\J_a \times \JJ_a}. 
 \end{eqnarray}
 Let us consider  the following
family:
$\{X_j, \quad j=0, \dots, j_{max}\}$ of submaniflods $ZZ_{\J_a \times \JJ_a}$:
$$X_0 = S^{n-\frac{n-m_{\sigma}}{8}+9}/\i  \times S^7/\i,$$
$$X_1= S^{n-\frac{n-m_{\sigma}}{8}+1}/\i \times S^{15}/\i, \quad \dots$$
$$ X_j = S^{n-\frac{n-m_{\sigma}}{8}+9-8j}/\i \times S^{8j+7}/\i, \quad$$
$$X_{j_{max}} = S^7/\i \times S^{n-\frac{n-m_{\sigma}}{8}+9}/\i,$$
where
\begin{eqnarray}\label{jmax}
j_{max}=\frac{7n+m_{\sigma}+16}{64}.
\end{eqnarray}

The dimension of each submanifold in the family is equal to $n-\frac{n-m_{\sigma}}{8}+16$
and the codimension in $ZZ_{\J_a \times \JJ_a}$ is equal to
$n-\frac{n-m_{\sigma}}{8}+2$. Let us define an embedding $$X_j \subset ZZ_{\J_a \times \JJ_a}$$
as the product of the two standard inclusions. 
Obviously, we get $\chi^{[4]}(X_{2j}) = X_{j_{max}-j}$.

The subpolyhedron $X_{a \times \aa} = \bigcup_{j=0} ^{j_{max}} X_{j}
\subset ZZ_{\J_a \times \JJ_a}$ is well defined. This subpolyhedron is invariant with respect to
the involution
 $\chi^{[4]}$. 
 %The polyhedron $X_{a \times \aa}$  is a manifold with singularities in 
 %the codimension $2$. 
 The restriction of the involution 
 $\chi^{[4]}$ on the subpolyhedron 
$X_{a \times \aa}$ denote the same.

Let us consider the following sequence of the index 2 subgroups, as in the Diagram $(\ref{a,aa})$:
\begin{eqnarray}\label{zep}
\I_{b \times \bb} \stackrel{{\rm i}_{b \times \bb,\H_{b \times \bb}}}{\longrightarrow} \H_{b \times \bb} 
\stackrel {{\rm i}_{b \times \bb, \J_a \times \JJ_a}}{\longrightarrow} \J_a
\times \JJ_a.
\end{eqnarray}
Let us define the following tower of 2-sheeted covers, associated
with the  sequence $(\ref{zep})$:
%\begin{eqnarray} \label{zz}
%ZZ_{b \times \bb}             \stackrel{p_{ZZ_{b \times
%\bb},ZZ_{\H_{b \times \bb}}}}{\longrightarrow} ZZ_{\H_{b \times \bb}} \stackrel{p_{ZZ_{\H_{b \times \bb}},ZZ_{\J_a %\times
%\JJ_a}}}{\longrightarrow} ZZ_{\J_a \times \JJ_a}.
%\end{eqnarray}
\begin{eqnarray} \label{zz}
ZZ_{\I_{b \times \bb}}   \longrightarrow ZZ_{\H_{b \times \bb}} \longrightarrow ZZ_{\J_a \times \JJ_a}.
\end{eqnarray}

The bottom space of the tower $(\ref{zz})$ is the standard
 skeleton
of the Eilenberg--Mac Lane space: $ZZ_{\J_a \times \JJ_a}
\subset K(\J_a,1) \times K(\JJ_a,1)$.
 The tower of 2-sheeted coverings
 $$K(\I_{b \times \bb},1) \to K(\H_{b \times \bb},1)  \to K(\J_a \times \JJ_a,1),$$
 associated with the sequence $(\ref{zep})$ is
well-defined. This tower determines the tower $(\ref{zz})$ by
means of the inclusion 
$ZZ_{\J_a \times \JJ_a} \subset K(\J_a,1)
\times K(\JJ_a,1) = K(\J_a \times \JJ_a,1)$.

Let us define the following tower of 2-sheeted covers:
\begin{eqnarray}\label{z}
X_{b \times \bb}   \longrightarrow X'_{b \times \bb}  \longrightarrow X_{a \times \aa}.
\end{eqnarray}
The bottom space of the tower $(\ref{z})$ is the subspace of the
bottom space of the tower $(\ref{zz})$, the inclusion is given by the inclusion 
$X_{a \times \aa} \subset ZZ_{\J_a \times
\JJ_a}$, see $(\ref{Xaa})$. The tower
$(\ref{z})$ is defined as the restriction of the tower
$(\ref{zz})$ over this subspace of the bottom.

Let us describe the top polyhedron of the tower $(\ref{z})$, which is a 
subpolyhedron $X_{b \times \bb} \subset ZZ_{\I_{b \times \bb}}$
explicitly. Let us define the
family $\{X_0,X_2,  \dots,  X_{j_{max}} \}$ of the standard submanifolds in 
$ZZ_{\I_{b \times \bb}} =
\RP^{n-\frac{n-m_{\sigma}}{8}+9} \times \RP^{n-\frac{n-m_{\sigma}}{8}+9}$
 by the following formulas:
\begin{eqnarray}\label{Xj}
X_0 = \RP^{n-\frac{n-m_{\sigma}}{8}+9} \times \RP^{7} \subset
\RP^{n-\frac{n-m_{\sigma}}{8}+9} \times \RP^{n-\frac{n-m_{\sigma}}{8}+9}, \dots
\end{eqnarray}
$$
X_j = \RP^{n-\frac{n-m_{\sigma}}{8}+9-8j} \times \RP^{8j+7} \subset
\RP^{n-\frac{n-m_{\sigma}}{8}+9} \times \RP^{n-\frac{n-m_{\sigma}}{8}+9}, \dots
$$
$$
X_{j_{max}} = \RP^{7} \times \RP^{n-\frac{n-m_{\sigma}}{8}+9} \subset
\RP^{n-\frac{n-m_{\sigma}}{8}+9} \times \RP^{n-\frac{n-m_{\sigma}}{8}+9}.
$$
 In this formula $j_{max}$ is defined by the formula
$(\ref{jmax})$. The subpolyhedron $X_{b \times \bb} \subset ZZ_{\I_{b
\times \bb}}$ is defined as the union of the standard submanifolds
in this family. The description of the subpolyhedron 
$X'_{b \times \bb} \subset ZZ_{\H_{b \times \bb}}$ 
 is obvious and
omitted.

On the space $ZZ_{\H_{b \times \bb}}$ (correspondingly, $ZZ_{\I_{b \times \bb}}$) the standard involution is well-defined
\begin{eqnarray}\label{ZZ1}
\chi^{[3]}: ZZ_{\H_{b \times \bb}} \to ZZ_{\H_{b \times \bb}} 
\end{eqnarray}
(correspondingly,
\begin{eqnarray}\label{ZZ2}
\chi^{[2]}: ZZ_{\I_{b \times \bb}} \to ZZ_{\I_{b \times \bb}})
\end{eqnarray}
by the formula $\chi^{[3]}(x \times y)=(y \times x)$ (correspondingly, $\chi^{[2]}(x \times y)=(y \times x)$).
The subpolyhedron  $X_{b \times \bb} \subset 
ZZ_{\I_{b \times \bb}}$ is invariant with respect 
to the involution    $\chi^{[2]}$.  The restriction of the involution 
$\chi^{[2]}$ on $X_{b \times \bb}$) is denoted  the same.

Let us define the submanifold $YY_{\H_{b \times \bb}} \subset
ZZ_{\H_{b \times \bb}}$  by the following formula:
$$ YY_{\H_{b \times \bb}} = (\RP^{n-\frac{n-m_{\sigma}}{4}+9} \times
\RP^{n-\frac{n-m_{\sigma}}{4}+9})/\i \subset $$
$$\RP^{n-\frac{n-m_{\sigma}}{8}+9} \times
\RP^{n-\frac{n-m_{\sigma}}{8}+9})/\i \subset ZZ_{\H_{b \times \bb}}.$$

Let us define the family of the standard submanifolds in 
$YY_{\H_{b \times \bb}} = (\RP^{n-\frac{n-m_{\sigma}}{4}+9} \times
\RP^{n-\frac{n-m_{\sigma}}{4}+9})/\i$.
\begin{eqnarray}\label{Yj}
Y_0 = (\RP^{n-\frac{n-m_{\sigma}}{4}+9} \times \RP^{7})/\i \subset
(\RP^{n-\frac{n-m_{\sigma}}{4}+9}) \times \RP^{n-\frac{n-m_{\sigma}}{4}+9})/\i, \dots
\end{eqnarray}
$$
Y_j = (\RP^{n-\frac{n-m_{\sigma}}{4}+9-8j} \times \RP^{8j+7})/\i \subset
(\RP^{n-\frac{n-m_{\sigma}}{4}+9} \times \RP^{n-\frac{n-m_{\sigma}}{4}+9})/\i, \dots
$$
$$
Y_{j'_{max}} = (\RP^{7} \times \RP^{n-\frac{n-m_{\sigma}}{4}+9})/\i \subset
(\RP^{n-\frac{n-m_{\sigma}}{4}+9} \times \RP^{n-\frac{n-m_{\sigma}}{4}+9})\i,
$$
where
\begin{eqnarray}\label{j'max}
j'_{max} =\frac{3n+m_{\sigma}+8}{32}.
\end{eqnarray}%j_{max}=\frac{7n+m_{\sigma}+16}{64}.
In the formula $(\ref{Yj})$ the action of the generator  $\i \in \Z/4$ 
is defined as the diagonal action, which is the standard on the factors.

Define the polyhedron 
\begin{eqnarray}\label{Y}
 Y_{b \times \bb} \subset YY_{\H_{b \times \bb}}
\end{eqnarray}
as the union of the submanifold of the family  $\{Y_j\}$. 

The polyhedron $Y_{b \times \bb}$
 is equipped with the involutions  $\chi^{[3]}$. 
 The definition of the involutions is standard.

The maps $(\ref{etaX})$, $(\ref{etaY})$ correspond to the
inclusions of the subgroups in ($\ref{zep}$) and are commuted with
the corresponded restrictions of maps in the tower of cover ($\ref{z}$).

  Let us denote the cylinder of the involution  $\chi^{[2]}$ on $X_{b \times \bb}$ by 
$X_{b \times \bb} \int_{\chi^{[2]}} S^1$.

The following natural inclusion is well defined: 
\begin{eqnarray}\label{etaX}
\eta_X: X_{b \times \bb} \int_{\chi^{[2]}} S^1 \subset
K(\I_{b \times \bb} \int_{\chi^{[2]}} \Z),1).
\end{eqnarray}

The space $Y_{b \times \bb} \int_{\chi^{[3]}} S^1$
 and the standard inclusion 
 \begin{eqnarray}\label{etaY}
\eta_Y: Y_{b \times \bb} \int_{\chi^{[3]}} S^1 \subset
K(\H_{b \times \bb} \int_{\chi^{[3]}} \Z,1).
\end{eqnarray}
are defined analogously.

Let us define a polyhedron  $J_X$.
For an arbitrary $j= 0,2, \dots, j_{max}$, where $j_{max}$ is given
by the formula $(\ref{jmax})$,  let us define the polyhedron $J_{j} =
S^{n-\frac{n-m_{\sigma}}{8}-8j+7} \times S^{8j+7}$
 The spheres
(the factors of this Cartesian product) $S^{n-\frac{n-m_{\sigma}}{8}-8j+9}$, $S^{8j+7}$
 is re-denoted
by $J_{j,1}$, $J_{j,2}$ correspondingly. Using this notation we
have:
 $$J_{j} = J_{j,1} \times J_{j,2}.$$

The standard inclusion
$i_{J_{j}}: J_{j,1} \times J_{j,2} \subset S^{\frac{n-m_{\sigma}}{8}+9} \times S^{\frac{n-m_{\sigma}}{8}+9}$
 is well-defined, each factor embeds in the sphere-image as the
standard subsphere. The union $\bigcup_{j=0}^{j_{max}}
Im(i_{J_j})$  of the images of these inclusions will be denoted by
\begin{eqnarray}\label{JX}
J_X  \subset
S^{\frac{n-m_{\sigma}}{8}+9} \times S^{\frac{n-m_{\sigma}}{8}+9}.
\end{eqnarray}
The polyhedron $J_X$ is well-defined.

Let us define a polyhedron (a manifold with singularities)  $J_Y$,
this polyhedron is a subpolyhedron on the polyhedron  $J_X$:
\begin{eqnarray}\label{JYJX}
{\rm{i}}_{J_Y,J_X}: J_Y \subset J_X.
\end{eqnarray}
For an arbitrary $j= 0, \dots, j'_{max}$, where $j'_{max}$ is given
by the formula $(\ref{j'max})$,  let us define the polyhedron 
$J_{j} =S^{n-\frac{n-m_{\sigma}}{4}-8j+9} \times S^{8j+7}$
The spheres
(the factors of this Cartesian product) $S^{n-\frac{n-m_{\sigma}}{4}-8j+9}$, $S^{8j+7}$
 are re-denoted
by $J_{2j,1}$, $J_{2j,2}$ correspondingly. Using this notation we
have again the formula $(\ref{Jj})$.

The standard inclusion ${\rm i}_{J_{j}}: J_{j,1} \times J_{j,2}
\subset S^{\frac{n-m_{\sigma}}{4}+9} \times S^{\frac{n-m_{\sigma}}{4}+9}$
is well-defined, each factor embeds in the sphere-image as the
standard subsphere. The union $\bigcup_{j=0}^{j'_{max}}
\mathrm{Im}(i_{J_j})$ of the images of these inclusions will be denoted by
$J_Y  \subset
J_X \subset S^{\frac{n-m_{\sigma}}{8}+9} \times S^{\frac{n-m_{\sigma}}{8}+9}$.
The polyhedron $J_Y$ is well-defined.

Define the standard 4-sheeted covering with ramification 
\begin{eqnarray}\label{varphiX}
\varphi_{X_{d \times \dd}}: X_{b \times \bb} \to J_X
\end{eqnarray}
the standard 2-sheeted covering with ramification 
\begin{eqnarray}\label{varphiY}
\varphi_{Y_{b \times \bb}}: Y_{b \times \bb} \to J_Y.
\end{eqnarray}

Consider the
$\I_{b \times \bb}$--covering $\bar X_{b \times \bb} \to X_{b \times \bb}$. 
The total space $\bar X_{b \times \bb}$ is the union of products of pairs of spheres, each product
is the 4-sheeted covering over the corresponding pair of the projective spaces in the formula
 $(\ref{Xj})$.  Define the standard $16$-sheeted covering with ramification
 $\bar X_{b \times \bb} \to J_X$. This covering 
over each product of a pair of spheres is defined as the Cartesian product of  joins
 of the corresponding cyclic
 $\J_a$ и $\JJ_a$--coverings over  $S^1$. The $\I_{b \times \bb}$--action on 
$\bar X_{b \times \bb}$ is commuted with 
the covering with ramification, constructed above. Therefore the covering $(\ref{varphiX})$ is well defined by 
factorization of the considered caverning
with ramification with respect to considered action.
The definition $(\ref{varphiY})$ is analogous. 

The polyhedra
$J_X$ and $J_Y$ are equipped by involutions  $\chi_{J_X}$,  $\chi_{J_Y}$, which are defined analogously to 
the involution
$\chi^{[2]}$, $\chi^{[3]}$. The cylinders of this involutions are well defined, 
denote those cylinders by 
$J_Y \int_{\chi} S^1$,  $J_X \int_{\chi} S^1$. 
The covering with ramification
 $(\ref{varphiX})$ commutes with the action of the involutions 
 $\chi_{J_X}$, $\chi^{[2]}$.
Therefore the following covering with ramification 
\begin{eqnarray}\label{varphiXS}
c_X: X_{b \times \bb} \int_{\chi^{[2]}} S^1 \to J_X \int_{\chi} S^1,
\end{eqnarray}
is well defined, this covering admits a $4$-sheeted factor, which is denoted by 
 $\hat c_X$.

The covering with ramification
 \begin{eqnarray}\label{varphiYS}
c_Y: Y_{b \times \bb} \int_{\chi^{[3]}} S^1 \to J_Y \int_{\chi} S^1,
\end{eqnarray}
is defined analogously to $(\ref{varphiXS})$. 
This covering admits a natural $2$-sheeted factor, which is denoted by  $\hat c_Y$.

\begin{lemma}\label{embJ}
There exists an embedding
\begin{eqnarray}\label{iJ}
{\rm{i}}_{J_X \int_{\chi} S^1}:  (J_X \int_{\chi} S^1) \times D^{10} \subset D^{n-11} \times S^1 \times D^{10} \subset D^{n-10}
\times D^{10} \subset \R^n,
\end{eqnarray}
where $D^i$ is the standard $i$-dimensional disk (of a small radius)
$D^{n-11} \times S^1 \subset D^{n-10}$ is the standard embedded solid torus. 
\end{lemma}

\subsubsection*{Proof of Lemma $\ref{embJ}$}
Let us define the collection of $j_{max}+1$ standard coordinate
subspaces of dimension $\frac{7n + m_{\sigma}}{8} + 6$ in the Euclidean space
$\R^{n-13}$, each coordinate space in the collection contains the origin. Let
us consider the Euclidean space $\R^{2n-\frac{n-m_{\sigma}}{4} + 20}$,
and let us fix the Cartesian product structure
$\R^{n-\frac{n-m_{\sigma}}{8} + 10} \times \R^{n-\frac{n-m_{\sigma}}{8} + 10}$.
 Inside the first factor of
this Cartesian product let us take the collection of $\frac{j_{max}}+1$
 coordinate subspaces $V_{2j}$,
$j=0, \dots, j_{max}$, $\dim(V_{j})=8j+8$. The space $V_{j}$ is the
standard coordinate codimension 2 subspace in the space $V_{j+2}$,
the space $V_{j_{max}}$ coincides with the first factor of the
Cartesian product. Inside the second factor of this Cartesian
product let us take the collection of
$\frac{j_{max}}+1$ coordinate
subspaces
$$W_{2j}, \quad j=0,2, \dots, j_{max},$$
$\dim(W_{j})=n-\frac{n-m_{\sigma}}{8} + 10 - 8j$.
 The space  $W_{j+1}$ is the standard coordinate
subspace in the space $W_{j}$  of the codimension 8, the space $W_{0}$
coincides with the second factor of the Cartesian product.

Let us consider the following collection of subspaces:
\begin{eqnarray}\label{coll}
\{V_j \times W_j \subset \R^{n- \frac{n-m_{\sigma}}{8} + 10} \times \R^{n-
\frac{n-m_{\sigma}}{8} + 10}\},
\end{eqnarray}
$\dim(V_{j} \times W_{j}) = \frac{7n + m_{\sigma}}{8} + 18$, which contain the origin.  
Let us consider the subspace  $U=\bigcup_{j=0}^{j_{max}} V_{j} \times W_{j}$, $U
\subset \R^{\frac{7n + m_{\sigma}}{4} + 20}$, in this formula the union is taken over the spaces of the collection  $(\ref{coll})$. 

Take a one-parametric family of orthogonal projections
$\pi(t): \R^{n-\frac{n-m_{\sigma}}{8} + 10} \times
\R^{n-\frac{n-m_{\sigma}}{8} + 10} \to \R^{n-11}$, $t \in S^1$,
which  satisfies the following condition.

--1. The following equation 
$\pi(t + 180^o)=I_{antidiag} \circ \pi(t)$ is satisfied, where
$t \in S^1$, $ I_{antidiag}: \R^{n-\frac{n-m_{\sigma}}{8} + 10} \times
\R^{n-\frac{n-m_{\sigma}}{8} + 10} \to \R^{n-\frac{n-m_{\sigma}}{8} + 10} \times
\R^{n-\frac{n-m_{\sigma}}{8} + 10}$ is an orthogonal mapping, which is 
the identity on the diagonal and is antipodal on the antidiagonal.

--2. The kernel $Ker(\pi)$ of the projections $\pi(t)$ for each $t$ is a
linear subspace, which is denoted by $L(t) \subset
\R^{n-\frac{n-m_{\sigma}}{8} + 10} \times \R^{n-\frac{n-m_{\sigma}}{8} + 10}$,
$\dim(L(t))=\frac{3n + m_{\sigma}}{4} +  26 $. In the family of projections with the boundary condition
  $\pi(t)$ the space $L(t)$ for an arbitrary $t$ intersects each subspace 
 $V_j \times W_j$ of the collection  $(\ref{coll})$ only at the origin. 

Evidently, the family of projections $\pi(t)$, with the required properties exists.
For example, we may take first
$\pi(0)$  as an embedding on the antidiagonal and is the identity on the diagonal.
Because the dimension of the diagonal is even, this condition determines a family
 $\pi(t)$, which satisfies Condition 1. Then take a small alteration of the family, keeping the boundary condition.

 By the general position argument in the case  $\sigma \ge 5$
 the following equality is satisfied:
 $\dim(V_j \times W_j) + \dim(L(t)) = \frac{7n + m_{\sigma}}{8} + 18  +
\frac{3n + m_{\sigma}}{4} + 26 \le 2   (\frac{7n + m_{\sigma}}{8} + 10) - 1$. 
Therefore there exists a family $\pi(t)$, which is satisfies Condition 2 for each $t$.  
 
 Let us denote the constructed family of embeddings by 
 $i_U(t): U \subset \R^{n-11}$. Let us denote the family of embeddings of the standard unite disk,
 which is associated with the family  $i_U(t)$ by 
$i_{\bar{U}}(t): \bar{U}(t) \subset D^{n-13}$.  In this formula by  $\bar{U}(t)$
is defined the union of the image of the standard unite disk with the center at the origin, which is 
associated with the union $U(t)$ of the vector spaces of the collection 
$(\ref{coll})$ for the prescribed value of the parameter $t$. 
By the involution 
 $I_{antidiag}$ the induced involution 
$\bar{I}_{antidiag}: \bar{U}(0) \to \bar{U}(0)$ is well-defined.

Consider  the embedding
$\rm{i}: S^1 \times D^{n-11} \subset \R^{n-10}$. Take the composition of the one-parameter family of  
embeddings $i_{\bar{U}}(t)$, $t \in [0,180^0]$,
with the one-parameter family
 $\rm{i}$, this composition determines an embedding, denoted by $i_{\bar{U}}(t)$, $t \in [0,180^0]$.

The required embedding $(\ref{iJ})$ is defined by the composition
$${\rm{i}}_J = ({\rm{i}}_U \times {\rm{id}}_{\R^{10}}) \circ ({\rm{i}}_{J,U}
\times {\rm{id}}_{D^{10}}): (J_X \int_{\chi} S^1) \times D^{10} \subset $$
$$(\bar{U} \int_{\bar{I}_{antidiag}} S^1) \times D^{10} \subset
\R^{n-10} \times \R^{10} = \R^n.$$ 
 In this formula ${\rm{i}}_{J_X,U}: J_X \int_{\chi} S^1
\subset \bar{U} \int_{\bar{I}_{antidiag}} S^1$ is the
embedding, which is constructed by means of the collection of the standard
embeddings: $\{Z_{j,1} \times Z_{j,2} \subset V_{j} \times
W_{j}\}$, ${\rm{id}}_{D^{10}}: D^{10} \subset \R^{10}$ is the standard embedding,
$\rm{id}_{\R^{10}}: \R^{10} \subset \R^{10}$ is the standard
diffeomorphism. Lemma $\ref{embJ}$ is proved.

\subsection*{An $\H_{b \times \bb}$--structure for a formal mapping $d_X^{(2)}$}

Let us consider an arbitrary equivariant generic
${\rm{PL}}$--mapping
\begin{eqnarray}\label{X2}
 d_X^{(2)}: (X_{b \times \bb} \int_{\chi^{[2]}} S^1)^2  \to \R^n \times \R^n.
\end{eqnarray} 
Such a (equivariant) mapping let us called a {\it formal} mapping.

 Let us define an (open) polyhedron
of a (formal) self-intersection  of the mapping  $d_X^{(2)}$ as the quotient of the inverse image of the diagonal
$diag(\R^n \times \R^n) \subset \R^n \times \R^n$ without points on the diagonal of the origin space.
  Denote this polyhedron by  $\N^{(2)}=\N^{(2)}(d^{(2)}_X)_{\circ}$. The closure  $\N^{(2)}$ of the polyhedron $\N^{(2)}_{\circ}$ contains boundary, this boundary will be
denoted by
 $\partial
\N^{(2)}_X)$ (an analogous construction is in [A1], the formula (44)).

In the case the mapping  $(\ref{X2})$ is defined by the  extension of a ${\rm{PL}}$-mapping
\begin{eqnarray}\label{X}
 d_X: X_{b \times \bb} \int_{\chi^{[2]}} S^1 \to \R^n,
\end{eqnarray}  
the polyhedron
$\N^{(2)}_{\circ}$ coincides with the polyhedron $\N(d_X)$ of self-intersection of the
mapping $d_X$. In this case the boundary $\partial(Cl(\N^{(2)}))$ coincides with the polyhedron of critical points of the mapping $d_X$.

Consider the subpolyhedron 
 \begin{eqnarray}\label{Xbb}
 X_{b \times \bb} \subset X_{b \times \bb} \int_{\chi^{[2]}} S^1,
 \end{eqnarray}
 which is defined as the fiber over the marked point of the fibration
$X_{b \times \bb} \int_{\chi} S^1 \to S^1$.
 The restriction of the equivariant mapping  $(\ref{X2})$ to the subpolyhedron  $(\ref{Xbb})$
denote by  
 \begin{eqnarray}\label{X22}
 d_X^{[2]}: X_{b \times \bb}^2  \to  \R^n \times \R^n.
\end{eqnarray} 
The  self-intersection  polyhedron of the formal mapping
 $d_X^{[2]}$ denote by $\N^{[2]}=\N^{[2]}(d_X^{[2]})_{\circ} \subset \N^{(2)}(d_X^{(2)})_{\circ}$.

 Suppose the polyhedron  $\N^{(2)}_{\circ}$ contains a marked closed component, which is denoted by 
 \begin{eqnarray}\label{N(2)}
 \N_{\H_{b \times \bb}}^{(2)} \subset \N^{(2)}(d_X^{(2)})_{\circ}.
 \end{eqnarray}
 Then the polyhedron
 $\N^{[2]}_{\circ}$ also contains a marked closed, which is denoted by 
 \begin{eqnarray}\label{NbbNd2}
 \N_{\H_{b \times \bb}}^{[2]} \subset \N^{[2]}(d_X^{[2]})_{\circ}.
 \end{eqnarray}
The following inclusion 
$$\N_{\H_{b \times \bb}}^{[2]} \subset \N_{\H_{b \times \bb}}^{(2)}$$
is well-defined.

The structure mapping
$$ \zeta_{\circ}: \N^{(2)}_{\circ} \to K(\Z/2^{[3]},1), $$
which is analogous  to  the structure mapping  
 ([A1], formula (43) is well-defined.
 Let us consider the restriction of the structure mapping  $\zeta_{\circ}$ on the marked component
 of the self-intersection polyhedron: 
 $$ \zeta: \N^{(2)}_{\H_{b \times \bb}} \to K(\Z/2^{[3]},1). $$
 Assume that this mapping admits a reduction, given by a mapping
\begin{eqnarray}\label{xNint}
 \zeta_{\H_{b \times \bb}\int} : \N^{(2)}_{\H_{b \times \bb}} \to K(\H_{b \times \bb} \int_{\chi^{[3]}} \Z,1). 
 \end{eqnarray}
 In this case the restriction of the structure mapping
$$ \zeta: \N^{[2]}_{\H_{b \times \bb}} \to K(\Z/2^{[3]},1) $$ 
on the marked component   $(\ref{NbbNd2})$ admits a reduction by the mapping 
 \begin{eqnarray}\label{xN}
 \zeta_{\H_{b \times \bb}} : \N^{[2]}_{\H_{b \times \bb}} \to K(\H_{b \times \bb},1). 
 \end{eqnarray}

%Assume that the formal mapping  $(\ref{X2})$ has a holonomic self-intersection along the marked component
%$(\ref{N(2)})$ in the sense of [A1,Definition 24]. Then the canonical double covering  $\bar \N^{(2)}_{\H_{b \times \bb}}$ over the polyhedron %$(\ref{N(2)})$  admits an immersion into $X_{b \times \bb} \int_{\chi^{[2]}} S^1$ 
%(denote an immersed neighborhood of this immersed subpolyhedron by $U(\bar \N_{\H_{b \times \bb}}^{(2)}) \looparrowright X_{b \times \bb} %\int_{\chi^{[2]}} S^1$)
%and the polyhedron $(\ref{N(2)})$
%is a closed component
%of the polyhedron of the self-intersection of an immersion $U(\bar \N_{\H_{b \times \bb}}^{(2)}) \looparrowright \R^n$.

\subsubsection*{Generic homology classes} 

Let us formulate a homological condition.
Let $(s_1,s_2)$ be an an arbitrary pair of integers
satisfying the following conditions: $s_1 = 1 \pmod{2}$,  $s_2 = 1
\pmod{2}$, $s_1 + s_2 = n - \frac{n-m_{\sigma}}{8}$.

For an arbitrary pair $(s_1,s_2)$, described below, we shall
construct a manifold $X(s_1,s_2)= \RP^{s_1} \times
\RP^{s_2}$ and consider the
embedding $i_{X(s_1,s_2)}: X(s_1,s_2)
\subset X_{b \times \bb} \subset X_{b \times \bb} \int_{\chi} S^1$, 
which is defined as the extension of the embedding $X(s_1,s_2) \subset X_{b \times
\bb}$, which is the Cartesian product of the two standard coordinate inclusions

By the construction for each pair $(s_1,s_2)$ the following inclusion 
\begin{eqnarray}\label{XIbb}
{\rm{i}}_{X(s_1,s_2)}: X(s_1,s_2) \subset X_{b \times \bb} \stackrel{\eta_{X_{b \times \bb}}}{\longrightarrow} K(\I_{b \times \bb} \int_{\chi^{[2]}} \Z,1)
\end{eqnarray}
is well-defined.

Additionally, define the manifold  $\RP^{n - \frac{n-m_{\sigma}}{8}-1} \times S^1$, which we denote by  $X_{\infty}$
and define the embeddings 
\begin{eqnarray}\label{iXinf}
i_{X_{\infty}}: X_{\infty} \subset X_{b \times \bb} \int_{\chi^{[2]}} S^1,
\end{eqnarray}
\begin{eqnarray}\label{iXXinf}
{{\rm{i}}}_{X_{\infty}}: X_{\infty} \subset X_{b \times \bb} \int_{\chi^{[2]}} S^1 \stackrel{\eta_{X_{b \times \bb}}}{\longrightarrow} K(\I_{b \times \bb} \int_{\chi^{[2]}} \Z,1).
\end{eqnarray}

Consider the standard 2-sheeted covering 
 $p_{X_{b \times \bb} \times S^1}: X_{b \times \bb} \times S^1 \to  X_{b \times \bb} \int_{\chi^{[2]}} S^1$, 
for which the embedded circle 
 $pt \times S^1 \subset X_{b \times \bb} \times S^1$ 
is mapped to the circle 
 $p_{X_{b \times \bb} \times S^1}(pt \times S^1)$,  the image of this circle by the standard projection 
  $X_{b \times \bb} \int_{\chi} S^1 \to S^1$ is projected as the standard $2$--sheeted covering
   $S^1 \to S^1$.  Define the embedding   
 $\RP^{n - \frac{n-m_{\sigma}}{8}-1} \times S^1 \subset  X_{b \times \bb} \times S^1$ as the product 
of the standard embedding 
$\RP^{n - \frac{n-m_{\sigma}}{8}-1}  \subset  X_{b \times \bb}$ on the first factor with the identity mapping
 $S^1 \to S^1$. Define the mapping  $i_{X_{\infty}}: \RP^{n - \frac{n-m_{\sigma}}{8}-1} \times S^1 \to X_{d \times \dd} \int_{\chi} S^1$
as the composition of this embedding with the projection
 $p_{X_{d \times \dd} \times S^1}$. The embedding   $(\ref{iXinf})$ is well-defined.

Analogously to the case of framing immersions (see Definition of the homology class $(\ref{etadd16})$)
let us prove that  the images of the fundamental classes by means of the mappings $(\ref{XIbb})$, $(\ref{iXinf})$
determine the following homology classes
\begin{eqnarray}\label{collectionx}
\{[X(s_1,s_2)] =  {{\rm{i}}}_{X(s_1,s_2);\ast}[X(s_1,s_2)], \quad
[X_{\infty}] = i_{X_{\infty},\ast}: [X_{\infty}]  \in
\end{eqnarray}
$$H_{n - \frac{n-m_{\sigma}}{8}}(X_{b \times \bb} \int_{\chi^{[2]}} S^1;\Z/2[\Z/2]).$$

Consider  the following compositions:
$$\{ {\rm{i}}_{J_X \int_{\chi} S^1} \circ {\rm{i}}_{X(s_1,s_2)}: X(s_1,s_2) \to X_{b \times \bb} \subset X_{b \times \bb} \int_{\chi^{[2]}} S^1 \to S^1 \times D^{n-1} \subset \R^n\},$$
$$\{ {\rm{i}}_{J_X \int_{\chi} S^1} \circ {\rm{i}}_{X_{\infty}}: X_{\infty} \to X_{b \times \bb} \to X_{b \times \bb} \int_{\chi^{[2]}} S^1 \to S^1 \times D^{n-1} \subset 
\R^n\}.$$

Assume that the formal mapping $(\ref{X22})$ is the formal extension of a mapping
$d_{X_{b \times \bb} \int_{\chi} S^1}$. In this case the image of the polyhedron $J_X \int_{\chi} S^1$ 
by this mapping is a immersed framed manifold with singularities. 
This observation allow to define the elements
 $(\ref{collectionx})$ analogously to the regular case. 
 For a formal mapping $(\ref{X22})$ this construction is analogous.

The homology class
$$[X_{\infty}] \in H_{n - \frac{n-m_{\sigma}}{8}}(K(\I_{b \times \bb} \int_{\chi^{[2]}} \Z,1);\Z/2[\Z/2])$$
is not in the subgroup
$(\ref{DZ2lok})$, the last classes of the collection  $(\ref{collectionx})$ belong to this subgroup.
The collection  $(\ref{collectionx})$ corresponds to the standard basis of the subgroup
$$\mathrm{Im}(H_i(K(\I_{b \times \bb} \int_{\chi^{[2]}} \Z,1);\Z[\Z/2]) \to
H_{i}(K(\I_{b \times \bb} \int_{\chi^{[2]}} \Z,1);\Z/2[\Z/2])),$$
which is described in Lemma $\ref{obstr}$.

Assuming that the formal mapping
 $(\ref{X22})$ is holonomic. In this case each manifold $X(s_1,s_2)$ determines the polyhedron
 of the self-intersection  of the mapping
 $d_{X_{b \times \bb} \int_{\chi} S^1} \circ {\rm{i}}_{X(s_1,s_2)}$. 
In a general case the analogous polyhedron of the is well-defined, denote this polyhedron by 
$NX(s_1,s_2)_{\circ} = NX(s_1,s_2)(d^{[2]})_{\circ}$. By the construction the following inclusion is
well-defined:
\begin{eqnarray}\label{NX}
NX(s_1,s_2)_{\circ} \subset N(d_X^{[2]})_{\circ}.
\end{eqnarray}
Define the closed subpolyhedron
\begin{eqnarray}\label{NXbb}
NX(s_1,s_2)_{\H_{b \times \bb}} \subset N^{[2]}_{\H_{b \times \bb}}; \quad NX(s_1,s_2)_{\H_{b \times \bb}} \subset NX(s_1,s_2)_{\circ}
\end{eqnarray}
by the intersection of the polyhedron
 $(\ref{NX})$ with the component  $(\ref{NbbNd2})$.

The restriction 
$\zeta_{b \times \bb} \vert_{NX(s_1,s_2)_{\H_{b \times \bb}}}$%\zeta_{b \times \bb} : N^{(2)}_{b \times \bb}
of the structure mapping  $(\ref{xN})$ to the subpolyhedron  $(\ref{NXbb})$ denote by
\begin{eqnarray}\label{zetaNX}
\zeta_{\H_{b \times \bb}, NX(s_1,s_2)}: NX(s_1,s_2)_{\H_{b \times \bb}} \to K(\H_{b \times \bb},1).
\end{eqnarray}

Consider the fundamental class
$$[NX(s_1,s_2)_{\H_{b \times \bb}}]
 \in H_{\dim(NX)}(NX(s_1,s_2)_{\H_{b \times \bb}}),$$
where $\dim(NX(s_1,s_2))=n-\frac{n-m_{\sigma}}{4}$.

The manifold $X(s_1,s_2)$ is oriented and the codimension of the formal mapping
$d^{[2]}_{X(s_1,s_2)}$ is even. Therefore the following collection of the homology classes is well-defined:
\begin{eqnarray}\label{relativeX}
\{ \zeta_{\ast,\H_{b \times \bb}}([NX(s_1,s_2)_{\H_{b \times \bb}}]) \in
\end{eqnarray}
$$
\mathrm{Im} D_{n-\frac{n-m_{\sigma}}{4}}(\H_{b \times \bb};\Z[\Z/2]) \to D_{n-\frac{n-m_{\sigma}}{4}}(\H_{b \times \bb};\Z/2[\Z/2])\}.
$$

The transfer homomorphism
\begin{eqnarray}\label{transferX}
!: D_{n-\frac{n-m_{\sigma}}{4}}(\H_{b \times \bb};\Z[\Z/2]) \to D_{n-\frac{n-m_{\sigma}}{4}}
(\I_{b \times \bb};\Z[\Z/2]),
\end{eqnarray}
which is associated with  the double
covering $K(\I_{b \times \bb} \int_{\chi^{[2]}} \Z,1) \to
K(\H_{b \times \bb}\int_{\chi^{[3]}} \Z,1)$ is well-defined. 
The
collection of permanent homology classes (modulo 2)
\begin{eqnarray}\label{absoluteX}
\zeta_{\H_{b \times \bb},\ast}^{!}([NX(s_1,s_2)_{\H_{b \times \bb}}]) \in
 D_{n-\frac{n-m_{\sigma}}{4}}(\I_{b \times \bb};\Z/2[\Z/2])\}
\end{eqnarray}
 is well-defined. 
 
Define an extra collection of homology classes of the group 
$D_{n-\frac{n-m_{\sigma}}{4}}(\I_{b \times \bb};\Z/2[\Z/2])$.
Let us consider the polyhedron $(\ref{NX})$, generally speaking, with a boundary.  
The standard compactification $С\bar NX(s_1,s_2)$ of the canonical 2-sheeted covering
 $\overline{NX}(s_1,s_2)$ over $\overline{NX}(s_1,s_2)$ is well-defined. 
The polyhedron $С\overline{NX}(s_1,s_2)$ contains the following closed subpolyhedron
 $\overline{NX}(s_1,s_2)_{b \times \bb}$ as a marked component.
The polyhedron $С\overline{NX}(s_1,s_2)$ itself, is not closed, is equipped with the natural embedding
into the polyhedron
$X(s_1,s_2) \subset X_{b \times \bb}$.
 
Consider the mapping 
\begin{eqnarray}\label{etaXbb}
\eta_{X_{b \times
\bb}}:  X_{b \times \bb} \to K(\I_{b \times \bb},1),
\end{eqnarray}
which is defined by means of the mapping
$(\ref{etaX})$. The induced mapping
 $\eta_{X_{b \times
\bb}}\vert_{C\overline{NX}(s_1,s_2)}$ determines the permanent homology class 
\begin{eqnarray}\label{absX}
\eta_{b \times \bb,\ast}([C \overline{NX}(s_1,s_2)]) \in
D_{\dim(NX))}(\I_{b \times \bb};\Z/2[\Z/2]). 
\end{eqnarray}
It is easy to calculate 
 $\dim(NX))=n-\frac{n-m_{\sigma}}{4}$.  Obviously, the homology class
$(\ref{absX})$ is defined as the reduction modulo 2 of the corresponding integer 
homology class.

The following equation for each pair $(s_1,s_2)$ is well defined:
\begin{eqnarray}\label{cyclic,aaX}
\eta_{b \times \bb,\ast}([C \overline{NX}(s_1,s_2)])
= \zeta_{\H_{b \times \bb},\ast}^{!}([NX(s_1,s_2)_{\H_{b \times \bb}}]).
\end{eqnarray}
This equation determines the identity of the two homology
classes  $(\ref{absX})$ and $(\ref{absoluteX})$ in the group 
$D_{\dim(NX)}(\I_{b \times \bb}\Z/2[\Z/2])$.

\begin{definition}\label{IaIdrel}
Let us say that the formal mapping
 $d_X^{(2)}$ $(\ref{X2})$ %with holonomic self-intersection along a closed component
%$(\ref{N(2)})$ 
admits an
$\H_{b \times \bb}$--structure, if there exists a mapping  $(\ref{xNint})$,
which determines a reduction of the restriction of the structure mapping on the marked component.
Moreover, the restriction
 $(\ref{xN})$ of the mapping $(\ref{xNint})$ is elated with the mapping  $(\ref{etaX})$
 by the family of equations $(\ref{cyclic,aaX})$.
\end{definition}

\begin{proposition}\label{22}
There exists a formal mapping $(\ref{X2})$, for which the polyhedron of the formal self-intersection
contains a closed component  
$(\ref{N(2)}$ along which the self-intersection is holonomic. An $\H_{b \times \bb}$--structure
in the sense of Definition  $\ref{IaIdrel}$ is well-defined, which is given by the mapping
$(\ref{xN})$. 
\end{proposition}
\[  \]

The sketch of the proof of Proposition
$\ref{22}$ is in [A3,Section 7].
\[ \]

\subsection*{An $\J_a \times \JJ_a$--structure for a formal mapping $d_Y^{(2)}$}
%with holonomic self-intersection along the marked component} 

Definition of a $\J_a \times \JJ_a$--structure is analogous to Definition $(\ref{IaIdrel})$
of an $\H_{b \times \bb}$--structure. In this definition polyhedra are replaced by and their factors, and structure mappings are replaced corresponding quadratic extension. 

\begin{definition}\label{IaIarel}
Let us say that a formal mapping
 \begin{eqnarray}\label{dY}
d_Y^{(2)}: (Y_{b \times \bb} \int_{\chi^{[3]}} S^1)^2  \to  \R^n \times \R^n 
\end{eqnarray} 
with holonomic self-intersection along a closed marked component 
\begin{eqnarray}\label{Naa}
\N_{\J_a \times \JJ_a}^{(2)} \subset \N^{(2)}(d_Y^{(2)})_{\circ}
\end{eqnarray}
admits an
$\J_a \times \JJ_a$--structure, if there exists a mapping 
\begin{eqnarray}\label{zetNaa}
\zeta_{\J_a \times \JJ_a} : \N^{[2]}_{\J_a \times \JJ_a} \to K(\J_a \times \JJ_a,1), 
\end{eqnarray} 
which determine a reduction of the restriction of the structure mapping on the marked component, and which is related with the mapping
 $(\ref{etaY})$ by the following family of equations:
$$\eta_{b \times \bb,\ast}([C \overline{NY}(s_1,s_2)])
= \zeta_{\J_a \times \JJ_a,\ast}^{!}([NY(s_1,s_2)_{\J_a \times \JJ_a}]).$$
%$(\ref{cyclic,aaX})$.
\end{definition}

\begin{proposition}\label{23}
There exists a formal mapping
 $(\ref{dY})$, for which the polyhedron of formal self-intersection contains a closed marked component
$(\ref{Naa})$.%, along which the self-intersection is holonomic. 
Additionally, a
$\J_a \times \JJ_a$--structure in the sense of Definition $\ref{IaIarel}$ is well-defined by the mapping
$(\ref{zetNaa})$. 
\end{proposition}

\subsection*{Construction of 
$\Z/2^{[2]}$--framed immersion, which admits an $\H_{b \times \bb}$--structure, proof of Lemmas
 $\ref{Hbb}$ and
$\ref{Jaa}$}

Let us prove Lemma  $\ref{Hbb}$.
Let $y=[(g,\eta_N,\Psi)]$ be given, where $g$ is a $\Z/2^{[2]}$--framed
immersion, $g: N^{n-\frac{n-m_{\sigma}}{8}}
\looparrowright \R^n$, $\eta_{N}: N^{n-\frac{n-m_{\sigma}}{8}} \to
K(\Z^{[2]},1)$ is the
characteristic class of the $\Z/2^{[2]}$--framing $\Psi$.

By the condition of the lemma there exists a mapping
\begin{eqnarray}\label{povtorX}
\eta_{b \times \bb}: N^{n-\frac{n-m_{\sigma}}{8}} \to K((\I_{b \times
\bb}) \int _{\chi^{[2]}} \Z,1),
\end{eqnarray}
which determines a reduction of the characteristic mapping $\eta_N$. This mapping  satisfies the Conditions 1 and 2
of Definition $\ref{ad-struct}$.
The mapping $\eta_{b \times \bb}$ determines (up to a homotopy)
the mapping 
\begin{eqnarray}\label{etabbX}
\eta_{b \times
\bb,X}: N^{n-\frac{n-m_{\sigma}}{8}} \to X_{b \times \bb} \int_{\chi^{[2]}} S^1, 
\end{eqnarray}
because the polyhedron
 $X_{b \times \bb}\int_{\chi^{[2]}} S^1$ is a subspace of the Eilenberg-Mac Lane space 
   $K(\I_{b \times \bb} \int_{\chi^{[2]}} \Z,1)$, and 
   this subspace contains the standard skeleton of the dimension
$n-\frac{n-m_{\sigma}}{8}+1 = \dim(N)+1$.

Analogously to [Theorem 23,A1] let us define a $\Z/2^{[2]}$-framed immersion
$(g, \Psi, \eta_N)$, which is determined in the group
$Imm^{\Z/2^{[2]}}(n-\frac{n-m_{\sigma}}{8},\frac{n-m_{\sigma}}{8})$
the given element $y$. Additionally, the manifold $N^{n-2k}$ contains a closed marked component 
$N^{n-2k}_{b \times \bb} \subset N^{n-2k}$, and the self-intersection manifold  $L^{n-\frac{n-m_{\sigma}}{4}}_{b \times \bb}$ of the restriction immersion
$g_{b \times \bb} = g \vert_{N^{n-2k}_{b \times \bb}}$,
contains a closed component 
 $L^{n-\frac{n-m_{\sigma}}{4}}_{\H_{b \times \bb}}$, which is the first component in the formula  $(\ref{compL3a})$.

By  geometrical arguments a natural projection 
\begin{eqnarray}\label{WNN}
L^{n-4k}_{\H_{b \times \bb}} \to N^{(2)}_{\H_{b \times \bb}},
\end{eqnarray}
$n-4k=n-\frac{n-m_{\sigma}}{4}$.
is well-defined.

Define the mapping
 $(\ref{1zetaad})$, which is required in Definition  $\ref{ad-struct}$. Consider the mapping
$(\ref{xNint})$, which is constructed in Proposition  $\ref{22}$ and define the mapping
$(\ref{1zetaad})$ by the composition of the mapping
$(\ref{WNN})$ with $(\ref{xNint})$.

The following Lemma is analogous to Lemma
 33 in [A1].

\begin{lemma}\label{self-int}
Assume  $x \in Imm^{[2]}(n-2k,2k)$ be an arbitrary element, which is represented by a triple
$(g',\eta_{N'},\Xi')$. (We will use the assumption that the characteristic mapping
$\eta_{N'}$ admits an $\I_{b \times \bb} \int_{\chi^{[2]}} \Z$--reduction, given a mapping $\eta_{b \times \bb}$.) Let $y \in Imm^{[3]}(n-4k,4k)$ be an arbitrary element, which is represented by a triple $(h,\zeta_L,\Lambda)$.
Then there exists another triple $(g,\eta_N,\Xi)$, $g: N^{n-2k} \looparrowright \R^n$, which represents an element $x$, for which the immersion  $g$ is self-transversal. The self-intersection manifold of $g$, is an $\Z/2^{[3]}$--immersed manifold, represented the element
$\delta^{[3]}(x)$, contains a closed component $L^{n-4k}$. Moreover, the characteristic mapping
$\zeta_L$ admits an $\H_{b \times \bb} \int_{\chi^{[3]}} \Z$--reduction by a mapping
 $\zeta_{\H_{b \times \bb}}$, for which the canonical 2-sheeted covering
mapping $\bar \zeta_{\H_{b \times \bb}}: \bar L^{n-4k} \to K(\I_{b \times \bb} \int_{\chi^{[2]}} \Z ,1)$ 
is induced from the mapping
 $\eta_{b \times \bb}$ by the immersion $\bar L^{n-4k} \looparrowright N^{n-2k}$.
\end{lemma}

By construction the immersion
 $g \vert_{VN^{n-2k}}$, and therefore, the immersion   $(\ref{WNR})$ is a $\Z/2^{[2]}$--framed immersion,
 and its characteristic mapping admits a  $\I_{b \times \bb} \int_{\chi^{[2]}} \Z$--reduction.
 Let us apply Lemma $\ref{self-int}$ and construct a
$\Z/2^{[2]}$--framed immersion $(g,\eta_N,\Psi)$, for which the characteristic mapping 
admits an
$\I_{b \times \bb} \int_{\chi^{[3]}} \Z$ -- reduction, and for which the $\Z^{[3]}$--framed
 immersion of the self-intersects manifold contains a closed component   $(\ref{WNN})$, for which the mapping
 $(\ref{1zetaad})$ determines an $\H_{b \times \bb}$--reduction of the characteristic mapping.

Let us check that the immersion 
 $(g,\eta_N,\Psi)$ admits an $\H_{b \times \bb}$--structure. Let us check Condition 3 in Definition  П
$\ref{ad-struct}$. The proof of this condition is analogous to the calculation of the degree of the mapping 
$\kappa_0$ in the proof of [Proposition 28,  A1].

Let us decompose the image of the fundamental class of the manifold 
$N^{n-2k}_{b \times \bb}$ by the mapping
$(\ref{etabbX})$ over the basic homology classes   $(\ref{collectionx})$. 
The homology class  $[X_{\infty}]$ is not required, because of Property  2
of Definition $\ref{ad-struct}$.

Each homology class
 $(\ref{collectionx})$ satisfies the equation $(\ref{cyclic,aaX})$. Therefore the image of the fundamental class by the mapping  $(\ref{etabbX})$  also satisfies by the same relation $(\ref{dd!})$. The Condition 3 from Definition
$\ref{ad-struct}$ is proved.

Condition 1 follows by construction, because the manifold
 $N^{n-2k}$ contains the only component $N^{n-2k}_{b \times \bb}$. Let us proof Condition 2
 from Definition $\ref{ad-struct}$.

 %Достаточно доказать свойство 2 из Леммы $\ref{ad}$.
Consider the restriction of the mapping  $(\ref{etaX})$ on the subpolyhedron  $X_{\infty}$, 
which is defined by the formula $(\ref{iXinf})$.  Consider the mapping
$${\rm{i}}_{X(s_1,s_2)}: X_{\infty} \to K(\I_{b \times \bb} \int_{\chi^{[2]}} \Z,1),$$ 
which is defined by the formula  $(\ref{XIbb})$. Consider the image of the homology class
\begin{eqnarray}\label{im}
{\rm{i}}_{X(s_1,s_2),\ast}([((\bar{w}_{\frac{n-m_{\sigma}}{8}})^{7})^{op}]) \in
H_{m_{\sigma}}(K(\I_{b \times \bb} \int_{\chi^{[2]}} \Z,1);\Z/2[\Z/2]),
\end{eqnarray}
where by $\bar{w}_{\frac{n-m_{\sigma}}{8}}$  is denoted the normal Stiefel-Whitney
class of the manifold
$X_{\infty}$ of the dimension $\frac{n-m_{\sigma}}{8}$.

\begin{lemma}\label{imX}
The homology class $(\ref{im})$ is not in the subgroup 
$D_{m_{\sigma}}(\I_{b \times \bb};\Z/2[\Z/2]) \subset H_{m_{\sigma}}(K(\I_{b \times \bb} \int_{\chi^{[2]}} \Z),1);\Z/2[\Z/2])$, in particular, this homology class in non-trivial.
\end{lemma}

\subsubsection*{Proof of Lemma $\ref{imX}$}

Let us prove that the total obstruction  $(\ref{ox})$ for the homology class $(\ref{im})$ is non-trivial.
Apply to the considered homology class the forgetful homomorphism  $(\ref{forg})$ and consider this homology class
with local coefficients as an integer homology class. 

Let us calculate the homology class $(\ref{im})$. At the firs step let us calculate 
the homology class
\begin{eqnarray}\label{opw}
((\bar{w}_{\frac{n-m_{\sigma}}{8}})^{7})^{op} \in H_{m_{\sigma}}(\RP^{n - \frac{n-m_{\sigma}}{8}-1} \times S^1;\Z/2[\Z/2]).
\end{eqnarray}
A direct calculation shows (in this calculation the inequality $\sigma \ge 5$ is used: the integer
$\frac{n-m_{\sigma}}{8}$ is divided by $4$), that the cohomology class 
$\bar{w}_{\frac{n-m_{\sigma}}{8}} \in H^{\frac{n-m_{\sigma}}{8}}(\RP^{n - \frac{n-m_{\sigma}}{8}-1} \times S^1;\Z/2)$
is the generic class of the first factor. Therefore, the cohomology class
 $(\bar{w}_{\frac{n-m_{\sigma}}{8}})^{7}$ is also the generic class of the first factor. Consider
the homology class 
$(\ref{opw})$ and let us divide this class by the 1-dimensional cohomology class, which is induced from
the generic cohomology class
 $H^1(S^1;\Z)$ by the standard projection $X_{\infty} \to K((\I_b \times \II_b) \int_{\chi^{[2]}} \Z) \to S^1$. 
%которая ассоциирована с гомоморфизмом
%$(\ref{epid})$.
Consider the image of this class by the mapping
${\rm{i}}_{X(s_1,s_2)}$. This image coincides with the class  $(\ref{tdtdd})$ for $i=m_{\sigma}-1$ and is non-trivial.
Lemma $\ref{imX}$ is proved.
\[  \]

\subsubsection*{The last step of the proof of Lemma $\ref{Hbb}$. The proof of Property  2 in Lemma $\ref{ad}$}

Let us assume that the homology class
$(\ref{etaadlok})$ is not in the subgroup  $(\ref{DZ2lokad})$, $i=n-4k$.
This means that the image of the fundamental class of the manifold 
$N^{n-2k}_{b \times \bb}$ by the mapping
 $(\ref{etabbX})$ is decomposed over
 the basis of this homology  group, such that the basic class
$[X_{\infty}]$ is not involved. Then by Lemma 
$\ref{imX}$, Property 2 of Definition $\ref{ad-struct}$ is not satisfied.
%Но это свойство доказано в Лемме $\ref{dd}$.
Statement 2 of Lemma  $\ref{ad}$ is proved.
Lemma $\ref{Hbb}$ is proved. 

The proof of Lemma $\ref{Jaa}$ is analogous.

\subsection*{Proof of Theorem  $\ref{th13}$}

%STOP

Let us start to prove Theorem $\ref{th13}$. Let us define the space
\begin{eqnarray}\label{ZZQ}
ZZ_{\I_a \times \II_a}= S^{\frac{n+m_{\sigma}}{2}+1}/\i \times
S^{\frac{n+m_{\sigma}}{2}+1}/\i.
\end{eqnarray}
(The space 
$(\ref{ZZIaIa})$ is not used below).
Obviously,  $\dim(ZZ_{\I_a \times \II_a}) = n + m_{\sigma}+2 > n$.

Let us define a family  $\{Z_1,  \dots,  Z_{j_{max}} \}$,
of standard submanifolds in the manifold $ZZ_{\I_a \times \II_a}$, where
\begin{eqnarray}\label{jmaxQ}
j_{max} = \frac{n+m_{\sigma}+4}{m_{\sigma}+2},
\end{eqnarray}
by the following formula:
\begin{eqnarray}\label{Zj}
Z_1 = S^{n-\frac{n-m_{\sigma}}{2}+1} \times
S^{\frac{m_{\sigma}}{2}}/\i \subset
S^{n-\frac{n-m_{\sigma}}{2}+1}/\i \times
S^{n-\frac{n-m_{\sigma}}{2}+1}/\i, \dots
\end{eqnarray}
$$
Z_j = S^{n-\frac{n-m_{\sigma}}{2}+ 1-(\frac{m_{\sigma}}{2}+1)j}/\i
\times S^{(\frac{m_{\sigma}}{2}+1)j-1}/\i \subset
S^{n-\frac{n-m_{\sigma}}{2}+1}/\i \times
S^{n-\frac{n-m_{\sigma}}{2}+1}/\i,  \dots
$$
$$
Z_{j_{max}} = S^{\frac{m_{\sigma}}{2}}/\i \times
S^{n-\frac{n-m_{\sigma}}{2}+1}/\i \subset
S^{n-\frac{n-m_{\sigma}}{2}+1}/\i \times
S^{n-\frac{n-m_{\sigma}}{2}+1}/\i.
$$

In the formulas   $j_{max}$ is defined by the formula
$(\ref{jmaxQ})$. The subpolyhedron  
\begin{eqnarray}\label{Z}
Z_{a \times \aa} \subset ZZ_{\I_a \times \II_a}
\end{eqnarray}
 is defined as the union of the standard
submanifolds of the family  $(\ref{Zj})$. Obviously, $\dim(Z_{a
\times \aa})=      \frac{n+2m_{\sigma}+2}{2}$.

The standard involution
\begin{eqnarray}\label{ZZ3}
\chi^{[4]}: ZZ_{\I_a \times \II_a} \to ZZ_{\I_a \times \II_a}, 
\end{eqnarray}
id defined analogously to the involution
 $(\ref{ZZ2})$, $(\ref{ZZ1})$. The subpolyhedron
 $(\ref{Z})$ is invariant with respect to the involution $(\ref{ZZ3})$.

Let us consider the standard cell decomposition of the space
$K(\I_a \times \II_a,1)= K(\I_a,1) \times K(\II_a ,1)$,
 this decomposition is defined as the Cartesian product of the standard cell decompositions of the factors.
 
The following standard inclusion is well defined 
\begin{eqnarray}\label{etaIaIaZ}
Z_{a \times \aa} \subset K(\J_a
\times \JJ_a,1).
 \end{eqnarray}

 The skeleton 
 $(\ref{etaIaIaZ})$ is invariant with respect to the involution  $(\ref{ZZ3})$. Therefore the inclusion
\begin{eqnarray}\label{ostov}
Z_{a \times \aa} \int_{\chi^{[4]}} S^1 \subset Z_{a \times \aa} \int_{\chi^{[4]}} S^1 \subset K((\I_a \times \II_a) \int_{\chi^{[4]}} \Z,1), 
\end{eqnarray}
which is an analog of the inclusions
 $(\ref{etaX})$, $(\ref{etaY})$ is well defined.

Let us define a polyhedron  $J_Z$, this polyhedron is the analog
of the polyhedron $J_X$, see $(\ref{JX})$. Let us denote by
$JJ_{\I_a}$ the space of join of  $j_{max}+1$ copies of the standard 
lens spaces $S^{\frac{m_{\sigma}}{2}}/\i$. Let us denote by
$JJ_{\I_a}$ the space of join of  $j_{max}+1$ copies of the
standard lens spaces $S^{\frac{m_{\sigma}}{2}}/\i$. In this
formulas  $j_{max}$ is defined by  $(\ref{jmaxQ})$.

Define a
subpolyhedron $J_{\bar{Z},j} \subset
JJ_{\I_a} \times JJ_{\II_a}$
 by the
formula
$$J_{Z,j} = JJ_{\Q,j} \times JJ_{\I_a,j}, 1 \le i \le j_{max},$$
where $JJ_{\I_a,j} \subset JJ_{\I_a}$  is the subjoin with the
coordinates  $0 \le i \le j_{max}-j$, $JJ_{\I_a,i} \subset
JJ_{\I_a}$ is the subjoin with the coordinates  $1 \le i \le j$.
Define a
subpolyhedron $J_{Z,j} \subset
JJ_{\Q} \times JJ_{\QQ}$  by the
formula
$$JJ_{Z,j} = JJ_{\Q,j} \times JJ_{\QQ,j}, 1 \le i \le j_{max},$$
where $JJ_{\Q,j} \subset JJ_{\Q}$ is the subjoin with the
coordinates $0 \le
i \le j_{max}-j$, $JJ_{\QQ,j} \subset JJ_{\QQ}$ is the subjoin with the
coordinates  $1 \le i \le j$.

Let us define $J_Z$ by the formula:
\begin{eqnarray}\label{JZ}
J_Z = \bigcup_{j=1}^{j_{max}} J_{Z,j} \subset JJ_{\Q} \times JJ_{\I_a}.
\end{eqnarray}
Let us define $J_{\bar{Z}}$ by the formula:
\begin{eqnarray}\label{barJZ}
J_{\bar{Z}} = \cup_{j=1}^{j_{max}} J_{\bar{Z},j} \subset JJ_{\I_a} \times JJ_{\II_a}.
\end{eqnarray}

On the polyhedra
  $Z_{\I_a \times \II_a}$, $J_Z$ a free involutions are well defined. The both involutions are denoted by
 $T_{\Q}$, the involutions corresponds to the quadratic extension  $(\ref{QZ4})$.
 The polyhedron 
$(\ref{JZ})$ is invariant with respect to the involution $T_{\Q}$, the inclusion of the factorpolyhedra denote by 
\begin{eqnarray}\label{QJZ}
(J_Z)/T_{\Q} = \cup_{j=1}^{j_{max}} J_{Z,j}/T_{\Q} \subset JJ_{\I_a}/T_{\Q} \times JJ_{\II_a}/T_{\Q}.
\end{eqnarray}

The following standard 4-sheeted covering with ramification is well defined:
\begin{eqnarray}\label{nakr}
c_Z: Z_{a \times \aa} \to  J_Z, 
\end{eqnarray}
This covering is factorized to 1-sheeted covering with ramification
 $\hat c_Z$.  The covering $(\ref{nakr})$ is defined by the composition of the 2-sheeted covering
with ramification
\begin{eqnarray}\label{nakrTQ}
Z_{a \times \aa}/T_{\Q} \to (J_Z)/T_{\Q} 
\end{eqnarray} 
and the standard 2-sheeted covering
$Z_{a \times \aa} \to Z_{a \times \aa}/T_{\Q}$.
The covering
 $(\ref{nakr})$ is invariant with respect to the standard involution, which corresponds to the involution  $(\ref{ZZ3})$,
 this involution is not re-denoted. 

The following embeddings
 $S^{\frac{m_{\sigma}}{2}}/\Q \subset
\R^{m_{\sigma}-1}$, $S^{\frac{m_{\sigma}}{2}}/\i \subset
\R^{m_{\sigma}-1}$ are well defined. This embeddings in the case
$\sigma=5$ was constructed by Hirsch
 [Hi] (see also  [Ro]). The join of the several copies of the considered embeddings
determines the following embedding:
\begin{eqnarray}\label{Jemb}
J_{Z,j} \subset \R^{m_{\sigma}(j_{max}+2)-2}.
\end{eqnarray}
The following lemma is the analog of Lemma $\ref{embJ}$.

\begin{lemma}\label{embJZ}
There exists an embedding
\begin{eqnarray}\label{iJZ}
\rm{i}_{J_Z}: J_Z \int_{\chi} S^1  \subset
\D^{n-2} \times S^1 \subset \R^{n-1} \subset \R^n.
\end{eqnarray}
\end{lemma}

Denote by 
 $\tilde Z_{a \times \aa} \subset Z_{a \times \aa}$ the inclusion of the standard skeleton of dimension  $n-\frac{n-m_{\sigma}}{2}$.
 Let us consider  a generic ${\rm{PL}}$--mapping $d_{\tilde Z}: \tilde Z_{a \times \aa} \to D^{n-1} \times S^1 \subset \R^n$. Let us assume that
 the mapping $d_{\tilde Z}$ is generic. Let us consider a polyhedron of self-intersection of the map $d_{\tilde Z}$ and let us denote
 this polyhedron by $\N(d_{\tilde Z})$. Obviously, $\dim(\N(d_{\tilde Z}))=m_{\sigma}$.
 In particular, in the case  $n=2^{11}$
 we get $\dim(\N(d_{\tilde Z}))=30$. The polyhedron
 $\N(d_{\tilde Z})$ contains the      subpolyhedron of critical points of the mapping $d_{\tilde Z}$.
 This polyhedron is called the boundary of the polyhedron $\N(d_{\tilde Z})$  and is denoted by $\partial \N(d_{\tilde Z})
 \subset \N(d_{\tilde Z})$.
 Hence, the standard inclusion  $\partial \N(d_{\tilde Z}) \subset \tilde Z_{a \times \aa}$ is well-defined.
 Because of the upper row of the diagram $(\ref{a,Qaa})$, the inclusion of the subgroup
 $a \times \aa \subset \Z/2^{[4]}$ is well-defined.

 The following mapping called the structure mapping is well defined:
 \begin{eqnarray}\label{dZ}
 \zeta_{\N(d_{\tilde Z})}: (\N(d_{\tilde Z}),\partial \N(d_{\tilde Z})) \to (K(\Z/2^{[5]},1),K((\I_a \times \II_a) \int_{\chi^{[4]}} \Z,1)),
 \end{eqnarray}
 the restriction of this mapping to the boundary of the polyhedron
of self-intersection of the mapping $d_{\tilde Z}$ is the
composition of the standard inclusion
 $\partial \N(d_{\tilde Z}) \to K(\I_a \times \II_a)$ and the mapping $(\ref{etaIaIaZ})$.
 The space
 $K(\I_a \times \II_a,1)$ is equipped with a mapping (this mapping is an inclusion in the homotopy category)
into $K(\Z/2^{[4]},1)$, this mapping is defined by the composition of the mapping $K(\I_a \times \II_a,1) \subset K(\Z/2^{[4]},1)$ and the diagonal embedding $K(\Z/2^{[4]},1) \times K(\Z/2^{[4]},1) \subset K(\Z/2^{[5]},1)$.

 \begin{definition}\label{structdZ}
 Let us say that the mapping  $d_{\tilde Z}$ admits  $\Q \times \Z/4$--structure, if for the mapping  $(\ref{dZ})$
 the following conditions are satisfied.

 --1. The polyhedron  $\N(d_{\tilde Z})$ is decomposed into two components:
\begin{eqnarray}\label{compdZ}
\N(d_{\tilde Z}) = \N_{\Q \times \Z/4} \cup \N_{\I_a \times \II_a \times \Z/2},
\end{eqnarray}
where the boundary  $\partial \N(d_{\tilde Z})$ is contained in the component
$\N{\I_a \times \II_a \times \Z/2}$.

--2. The restriction of the mapping  $(\ref{dZ})$ to the component
$\N_{\Q \times \Z/4}$ admits a reduction to a
mapping into the subspace $K((\Q \times \Z/4) \int_{\chi^{[5]}} \Z,1)$ i.e. this
restriction is given by the following composition:
$$  \zeta_{\Q \times \Z/4}: \N_{\Q \times \Z/4} \to K((\Q \times \Z/4) \int_{\chi^{[5]}} \Z,1)
\to K(\Z/2^{[5]},1),$$
where the mapping of the classifying spaces  $K((\Q \times \Z/4) \int_{\chi^{[5]}} \Z,1) \to K(\Z/2^{[5]},1)$ is induced by the homomorphism
$(\ref{iQaa})$.

 --3. The restriction of the mapping $(\ref{dZ})$ to the component
 $\N_{\I_a \times \II_a \times \Z/2}$
is given by the following composition:
$$  \zeta_{\I_a \times \II_a \times \Z/2}:
(\N_{\I_a \times \II_a \times \Z/2},\partial
\N_{\I_a \times \II_a \times \Z/2})$$
$$ \longrightarrow (K((\I_a \times \II_a \times
\Z/2) \int_{\chi^{[5]}} \Z,1),K((\I_a \times \II_a) \int_{\chi^{[4]}} \Z,1),$$
where the
mapping of the classifying spaces $K(\I_a \times \II_a \times
\Z/2,1) \to K(\Z/2^{[5]},1)$ is induced by the homomorphism
$(\ref{iaaa})$.
 \end{definition}

The following lemma is analogous to Lemmas $\ref{22}$,
$\ref{23}$. The proof of this lemma is analogous to the proof of
Lemma 33 of [A1], see part ${\rm{III}}$ of the paper.

\begin{lemma}\label{134}
There exists a mapping $d_{\tilde Z}: \tilde Z_{a \times \aa}
\to D^{n-1} \times S^1 \subset \R^n$, admitting (a relative) $\Q \times \Z/4$--structure
(see  Definition
$(\ref{structdZ})$).

The mapping $d_{\bar{Y}}$ is defined as a generic alteration of the following composition:
\begin{eqnarray}
\tilde Z_{a \times \aa} \int_{\chi^{[4]}}S^1 \subset Z_{a \times \aa} \int_{\chi^{[4]}}S^1
\stackrel {\varphi}{\longrightarrow} J_Z \int_{\chi} S^1
\stackrel{\rm{i}_{J_Z}}{\subset} D^{n-2} \times S^1 \subset \R^{n-1} \subset \R^n,
\end{eqnarray}
where $c_Z:  Z_{a \times \aa} \int_{\chi^{[4]}} S^1 \to J_Z \int_{\chi} S^1$  is the 4--sheeted covering with ramification,
defined by means of the covering $(\ref{nakr})$,
$\tilde Z_{a \times \aa} \int_{\chi^{[4]}}S^1 \subset Z_{a \times \aa} \int_{\chi^{[4]}}S^1 $
 $J_{\tilde Z} \subset J_Z$ is the standard  embedding,
determined by means of the inclusion $(\ref{etaIaIaZ})$, ${\rm{i}}_J : J_Z \subset D^{n-2} \times S^1 \subset \R^n$
is the inclusion, constructed in Lemma
$\ref{embJZ}$.
\end{lemma}

\subsubsection*{Proof of Theorem  $\ref{th13}$ }

Consider the mapping $(\ref{n-8k})$, which is determined a reduction of the restriction of the characteristic mapping 
$\eta_N$ on the marked component. Without loss of the
generality, we may assume that the image of the map $\eta_{a
\times \aa}$ is contained in the subspace $\tilde Z_{a \times \aa} \int_{\chi^{[4]}} S^1 \subset K((\J_a
\times \JJ_a) \int_{\chi^{[4]}} \Z,1)$. Let us define an $\Z/2^{[4]}$--immersion 
$g_{a \times \aa}: N^{\frac{n+m_{\sigma}}{2}}_{a \times \aa} \looparrowright \R^n$, which is in
the regular homotopy class of the immersion
$g \vert_{N^{\frac{n+m_{\sigma}}{2}}_{a \times \aa}}$, and the self-intersection of this immersion
satisfies the required conditions from Definition 
 $\ref{Q-struct}$.

Let us consider the composition $N^{\frac{n+m_{\sigma}}{2}}_{a \times \aa}
\stackrel{\eta_{a \times \aa}}{\longrightarrow} \tilde
Z_{a \times \aa} \int_{\chi^{[4]}} S^1 \stackrel{d_{\tilde Z}}{\longrightarrow}
\R^n$.  Let us consider a $C^0$--small alteration of this
composition into an immersion $g_1$ in the regular homotopy class
of the immersion $g$.  Let us denote the self-intersection
manifold of the immersion $g_1$ by $L^{m_{\sigma}}_{a \times \aa}$.  The
caliber of the deformation above is chosen so small that the
manifold $L^{m_{\sigma}}_{a \times \aa}$ is decomposed into two disjoint
components as in the formula $(\ref{LQII})$.
%\begin{eqnarray}\label{comLL}
%L_1^{m_{\sigma}}_{a \times \aa} =   L^{m_{\sigma}}_{\Q \times \Z/4} \cup
%L^{m_{\sigma}}_{\I_a \times \II_a \times \Z/2}.
%\end{eqnarray}

The component $L^{m_{\sigma}}_{\Q \times \Z/4}$ is contained in
a small regular neighborhood of the image of the first component
$\N_{\Q \times \Z/4}$ from the decomposition
 $(\ref{compdZ})$. The component  $L^{m_{\sigma}}_{\I_a \times \II_a \times \Z/2}$
 is represented as the union of the two submanifolds with boundaries along the common boundary.
 The first of the manifolds with boundary is immersed into
 a small regular neighborhood of the image of the second components
$\N_{\I_a \times \II_a \times \Z/2}$ from the
decomposition $(\ref{compdZ})$. The second manifold with boundary
in immersed into a small neighborhood of regular values of the
mapping $d_{\tilde Z}$. The common boundary is immersed into a
small regular neighborhood of critical values of the map
$d_{\tilde Z}$. Therefore, for the manifold $L^{m_{\sigma}}_{a \times \aa}$  the formula
$(\ref{compon})$ is satisfied and the induced mapping
$(\ref{LQII})$ is well-defined.

 The immersion $g_{a \times \aa}$ is a $\Z/2^{[4]}$--framed immersion in
 the prescribed cobordism class of the restriction of the triple $(g,\eta_N, \Psi)$ on the marked component.
 By this condition the pair of mappings  $(\ref{n-8k})$, $(\ref{LQII})$
is well-defined, the mapping 
$\eta_{a \times \aa}:
N^{\frac{n+m_{\sigma}}{2}}_{a \times \aa} \to K((\J_a \times \JJ_a) \int_{\chi^{[4]}} \Z,1)$,
%удовлетворяющее уравнению $(\ref{etaIaIa})$, 
%пара отображений
corresponds to the mapping  $(\lambda=\zeta_{\Q \times \Z/4} \cup \zeta_{\J_a \times \JJ_a \times \Z/2}$, and the Conditions 1 and 2 from Definition 
$\ref{Q-struct}$ are satisfied. 
Theorem $\ref{th13}$ is proved.

\section{Compression Theorem}

Let us prove the Соmрression Theorem $\ref{comp}$.

% \begin{theorem}{Compression Theorem}\label{comp}

%For an arbitrary positive integer $d$ there exists a positive
%integer  $l=l(d)$ such that an arbitrary element in the cobordism
%group $Imm^{sf}(2^{l'}-3,1)$, при $l' \ge l$, admits a compression
%of the order $d-1$ (see Definition $7$ of [A2]).
%\end{theorem}

\subsubsection*{Remark}
A sketch of an alternative proof of the Compression Theorem $\ref{comp}$ was proposed  by D.
Ravenel, see [R].

%\subsubsection*{Remark}
%A necessary and sufficient  conditions of the existence of a
%compression of the order $d$ is reformulated in Corollary
%$\ref{totalobstr}$ below (in a more general situation). This allow
%to reformulate Theorem $\ref{comp}$  as follows. For an arbitrary
%positive integer $d$ there exists a positive integer $l=l(d)$ such
%that for an arbitrary  $l' \ge l$ the homomorphism $ J_{sf}^d:
%Imm^{sf}(2^{l'}-3,1) \to Imm^{sf}(d,2^{l'}-2-d)$ is trivial.
\[  \]

Let $M^{n-k}$ be a closed manifold of dimension $(n-k)$, $\varphi:
M^{n-k} \looparrowright \R^m$ be an immersion of this manifold
into $\R^n$ in the codimension $k$, $\Xi_M$ is a skew-framing of
the immersion $\varphi$ with the characteristic class $\kappa_M$
of this skew-framing. Additionally, let us assume that the
manifold $M^{n-k}$ is equipped with the family of 1-dimensional
cohomology classes modulo 2:
\begin{eqnarray}\label{A_M}
A_j=\{\kappa_i\}, \quad \kappa_i \in H^1(M;\Z/2), \quad 0 \le i
\le j, \quad \kappa_0 = \kappa_M.
\end{eqnarray}
This collection of cohomology classes is represented by the
collection of classifying maps:
\begin{eqnarray}\label{A'_M}
A_{j}=\{\kappa_i: M^{n-k} \to \RP^{\infty} \}, \quad i=0,1,
\dots, j, \quad \kappa_0 = \kappa_M.
\end{eqnarray}

\begin{definition}\label{cob}

The cobordism group of immersions $Imm^{sf;A_j}(n-k;k)$ is
represented by triples $(\varphi,\Xi_M,A_j)$, where:

 ---$\varphi: M^{n-k} \looparrowright \R^n$--is an immersion of a closed $(n-k)$-dimensional manifold into Euclidean space,

--- $\Xi_M$ is skew-framing of the immersion $\varphi$,

---$A_j$ is a collection  of cohomology classes, described in $(\ref{A'_M})$.

The cobordism relation of triples  is the standard.
\end{definition}

\subsubsection*{Remark}
In the case $j=0$ the cobordism group $Imm^{sf;A_j}(n-k;k)$
coincides with the cobordism group $Imm^{sf}(n-k;k)$ of
skew-framed immersions.
\[  \]

A natural homomorphism
\begin{eqnarray}\label{102}
J^{k}_{sf;A_j}: Imm^{sf;A_j} (n-1,1) \to Imm^{sf;A_j}(n-k,k)
\end{eqnarray}
is defined as follows. Let us assume that the triple $(\varphi_0,
\Xi_{M_0}, A_j(M_0))$ represents an element in the cobordism group
 $Imm^{sf;A_j}(n-1,1)$. Let us consider the following triple $(\varphi, \Xi_M, A_j(M))$,
where the immersion $\varphi: M^{n-k} \looparrowright \R^n$ define
as follows. The manifold $M^{n-k}$ is a submanifold in the
manifold $M_0^{n-1}$, the fundamental class of this submanifold
represents the homological Euler class of the bundle
$(k-1)\kappa_{M_0}$, the immersion $\varphi$ is defined as the
restriction of the immersion $\varphi_0$ on $M^{n-k}$. The
skew-framing $\Xi_M$ of the immersion $\varphi$ is defined by the
standard construction like in the case $j=0$, the collection
$A_j(M)$ of cohomology classes is the restriction of the
collection $A_j(M_0)$ on the submanifold $M^{n-k} \subset
M_0^{n-1}$.
\[  \]

Let us generalize Definition $\ref{5}$ for the cobordism group
$Imm^{sf;A_j}(n-k,k)$.

\begin{definition}\label{comp-5}
Let $[(\varphi,\Xi_M,A_j)] \in  Imm^{sf;A_j}(n-k,k)$. We shall say
that the element  $[(\varphi,\Xi_M,A_j)]$ admits a compression of
the order  $q$, if in its cobordism class there exists a triple
 $(\varphi',\Xi'_M,A'_j)]$, such that the pair  $(M^{n-k}, \kappa_0)$
admits a compression of the order $q$ is the sense of the
Definition $\ref{5}$.
\end{definition}

Let us define the transfer homomorphism
\begin{eqnarray}\label{transfer}
r_j^!:Imm^{sf;A_j}(n-k,k) \to Imm^{sf;\hat A_{j-1}}(n-k,k)
\end{eqnarray}
 with respect to the cohomology
class $\kappa_j$.

Let $x \in Imm^{sf;A_j}(n-k,k)$ be an element represented by a
triple $(\varphi, \Xi_M, A_j)$. Let us define the 2-sheeted cover
$$p_j:\hat M^{n-k} \to M^{n-k}$$
as the regular cover with the characteristic class
 $w_1(p_j)=\kappa_j \otimes \kappa_M \in H^1(M^{n-k};\Z/2)$.
(We will denote below by $p_j$ the linear bundle over $M^{n-k}$
with the characteristic class $w_1(p_j)$ and also the
characteristic class $w_1(p_j)$ itself.)

Let us define a skew-framing
 $\Xi_{\hat M}$. Let us consider the immersion
$\varphi: M^{n-k} \looparrowright \R^n$. Let us denote the
immersion $\varphi \circ p_j$ by $\hat \varphi$. Let us denote the
normal bundle of the immersion $\varphi$ by $\nu_{M}$. Let us
denote the normal bundle of the immersion $\varphi \circ p_j$ by
$\nu_{\hat M}$. Let us define the skew-framing $\Xi_{\hat M}$ of
the immersion $\varphi \circ p_j$ by the formula $\Xi_{\hat M} =
p_j^{\ast}(\Xi_M)$.

Let us define the collection of the cohomology classes $\hat
A_{j-1}$ by the following formula:
$$\hat A_{j-1} = \{ \hat \kappa_j = p_j \circ \kappa_j, \quad j=0, \dots, j-1 \}. $$

 We will define
$$r_j^!(x)=[(\hat \varphi,\Xi_{\hat M},\hat A_{j-1})].$$

\begin{example}\label{Eccles}
In the case  $j=0$ the transfer homomorphism $(\ref{transfer})$ is
given by the following formula:
$$r_{0}: Imm^{sf}(n-k,k) \to
Imm^{fr}(n-k,k)$$ The properties of this homomorphism was
considered in [A-E].
\end{example}

Let us consider the homomorphism, given by the composition of the
$(j_2-j_1)$ transfer homomorphisms:
\begin{eqnarray}\label{103}
r_{j_1+1, \dots, j_2}^!: Imm^{sf;A_{j_2}}(n-k,k) \to Imm^{sf;\hat
A_{j_1}}(n-k,k), \quad j_2 \ge j_1.
\end{eqnarray}
In the case $j_1=0$, the homomorphism $(\ref{103})$ will be
denoted by
\begin{eqnarray}\label{103tot}
r_{tot}^!: Imm^{sf;A_{j_2}}(n-k,k) \to Imm^{sf}(n-k,k).
\end{eqnarray}

Let us describe the homomorphism $(\ref{103})$ explicitly. Let us
assume that an element $x \in Imm^{sf;A_{j_2}}(n-k,k)$ is given by
a triple $(\varphi, \Xi_M, A_{j_2})$. Let us consider the the
following subcollection $\{\kappa_{j_1+1}, \dots, \kappa_{j_2} \}$
of the last $(j_2-j_1)$ cohomology classes of the collection
$A_{j_2}$. Let us define the $2^{j_2-j_1}$--sheeted cover $p: \hat
M^{n-k} \to M^{n-k}$, given by the following collection of the
cohomology classes:
$$\{\kappa_{j_1+1} \otimes \kappa_0, \dots, \kappa_{j_2} \otimes \kappa_0 \}.$$
The immersion $\hat \varphi: \hat M^{n-k} \looparrowright \R^n$ is
given by the following composition:
$$ \hat \varphi = \varphi \circ p. $$
The normal bundle of the immersion  $\hat \varphi$ is equipped
with the skew-framing $p^{\ast}(\Xi_M)=\Xi_{\hat M}$. The
collection of cohomology classes
$$\hat A_{j_1} = \{ \hat \kappa_0
= \kappa_0 \circ p, \dots \hat \kappa_{j_1}=\kappa_{j_1} \circ p
\}$$
is well-defined by the classifying maps.

The cobordism class of the triple $(\hat \varphi, \Xi_{\hat M},
\hat A_{j_1})$ determines the element $r_{j_1+1, \dots,
j_2}^!(x)$.

\begin{proposition}\label{totaltrans}
For an arbitrary even positive integer $k$, $k=0 (mod 2)$,
$2^l-2=n>k>0$, there exists a positive integer $\psi = \psi(k)$,
such that  the total transfer homomorphism $(\ref{103tot})$ for
$j_2=j$:
\begin{eqnarray}\label{104}
r_{tot}: Imm^{sf;A_{\psi}}(n-k,k) \to Imm^{sf}(n-k,k)
\end{eqnarray}
is trivial.
\end{proposition}

\subsubsection*{Proof of the Proposition $\ref{totaltrans}$}

Let us start with the following lemma.
\begin{lemma}\label{P-Th}

The cobordism group $Imm^{sf;A_j}(n-k,k)$ is a finite $2$-group if
$n-k >0$.
\end{lemma}

\subsubsection*{Proof of Lemma $\ref{P-Th}$}

The cobordism group $Imm^{sf}(n-k,k)$ by the Pontrjagin-Thom
construction is the stable homotopy group $\Pi_{n-k}(P_{k-1})$,
where $P_{k-1} = \RP^{\infty}/\RP^{k-1}$, see f.ex. [A-E]. By the
standard arguments the cobordism group $Imm^{sf;A_j}(n-k,k)$ is a
stable homotopy group $\Pi_{n-k}(P_{k-1} \times \prod_{i=1}^{j}
\RP^{\infty}_j))$. The space $P_{k-1} \times \prod_{i=1}^{j}
\RP^{\infty}_j)$ is a space of the 2-homotopy type. The Lemma
$\ref{P-Th}$ is proved.

Let us define the following sequence of positive integers
 $s_{n-k}, \dots,
s_1$, where the indexes decrease from $(n-k)$ to 1:
$$ s_{n-k} =
ord(\Pi_{n-k}), \quad s_{n-k-1}= ord(\Pi_{n-k-1}), \quad \dots,
\quad s_1= ord(\Pi_1).$$
 In this formula we denote by
$ord(\Pi_i)$ the logarithm of the maximal order of an element in
the 2-component of the $i$-th stable homotopy group of spheres.
Then let us define
 $\psi(i)=\sum_{j=i}^{n-k} s_j$ and let us define the required integer
 by the following formula:
 \begin{eqnarray}\label{psi}
\psi = \psi(1)+ \sigma,
\end{eqnarray}
where $\sigma = ord(Imm^{sf;A_{\psi(1)}}(n-k,k))$ (see Lemma
$\ref{P-Th}$). Therefore we have $\psi \ge \psi(1) \ge \psi(2) \ge
\dots \ge \psi(n-k)$.

Let  $x \in Imm^{sf;A_{\psi}}(n-k,k)$ be an  element represented
by a triple $(\varphi, \Xi_M, A_{\psi})$. Let us consider the
element
\begin{eqnarray}\label{element}
r_{\psi(1)+1, \dots, \psi}^!(x),
\end{eqnarray}
represented by a triple $(\hat \varphi, \Xi_{\hat M}, \hat
A_{\psi-\psi(1)})$, given by the transfer homomorphism with
respect to the subcollection $\{ \kappa_{\psi(1)+1}, \dots,
\kappa_{\psi} \}$ of the cohomology classes in the collection
$A_{\psi}(M)$ (the last $\sigma$ classes in the collection
$A_{\psi}$). Let us prove that the element $(\ref{element})$
admits a compression of the order 0, see the Definition
$\ref{comp-5}$.

Let us denote the product $\RP^{\infty}_0 \times
\prod_{j=1}^{j=j_0} \RP^{\infty}(j)$ by  $X(j_0)$ and let us
consider the map
\begin{eqnarray}\label{lambdaM}
\lambda_M: M^{n-k} \to X(s),
\end{eqnarray}
were $s=\psi_1$, defined as the
direct product of  the classifying maps of the collection
$A_{\psi-\psi(1)}$ of cohomology classes (the classes
$\kappa_{\psi(1)}+1, \dots, \kappa_{\psi}$ are omitted). Let us
denote again the space $X(\psi(1))$ by $X$ for short.

Let us consider a natural filtration
\begin{eqnarray}\label{105}
\dots \subset X^{(n-k+1)} \subset X^{(n-k)} \subset \dots \subset
X^{(1)} \subset X.
\end{eqnarray}
This filtration is the direct product of the two standard
coordinate filtrations. Each stratum  $X^{(i)} \setminus
X^{(i+1)}$, $i=1, \dots n-k$ is an union of open cells of the
codimension $i$. Each cell is determined by the corresponded
multi-index $\mu=(m_1, \dots, m_{\psi(n-k)})$, $m_1+ \dots +
m_{\psi(n-k)}=i$, $m_i \ge 0$, each coordinate of the multi-index
shows the codimension of the skeleton of the corresponded
coordinate projective space that contains the given cell.

Let us assume that the map $\lambda$ is in a general position with
respect to the filtration  $(\ref{105})$. Let us denote by $L^0
\subset M^{n-k}$ a $0$--dimensional submanifold in $M^{n-k}$,
defined as the inverse image of the stratum  $X^{n-k}$ of this
filtration.

Let us consider the triple $(\hat \varphi, \Xi_{\hat M}, \hat
A'_{\psi-\psi(1)})$. Let us prove that this triple is cobordant to
a triple $(\hat \varphi', \Xi_{\hat M'}, \hat A'_{\psi-\psi(1)})$
such that the mapping
 $\hat \lambda': \hat M'^{n-k} \to X$,
 constructed by means of the collection of the cohomology classes
 $A_{\hat M'}$, satisfies the following property:
 \begin{eqnarray}\label{106}
  \hat \lambda_1^{-1}(X^{(n-k)}) = \emptyset.
\end{eqnarray}

By the assumption  $k$ is even. Therefore, because $n$ is even,
the manifold $M^{n-k}$ is oriented. Let us consider an arbitrary
choose  the orientation of each cell in $X^{(n-k)} \setminus
X^{(n-k+1)}$. Let us  denote by $lk(\mu)$ the integer coefficient
of self-intersection of the image $\lambda(M^{n-k})$ with the
oriented cell of a multi-index  $\mu$. Let us denote analogically
by $lk(\hat \mu)$ the integer coefficient of self-intersection of
the image $\hat \lambda(\hat M^{n-k})$ with the oriented cell of a
multi-index $\mu$. Obviously, the collection of the integer
$\{lk(\hat \mu)\}$ is obtain from the collection of the
corresponded  integers $\{lk(\mu)\}$ by the multiplication on
$2^{\psi-\psi(1)}$.

Let $A'_{\psi(1)}$ be the subcollection of $A_{\psi}$, consists of
the first $\psi(1)$ cohomology classes. By the construction, the
exponent of the group $Imm^{sf;A_{\psi(1)}}(n-k,k)$ is equal to
$2^{\psi-\psi(1)}$. Therefore the disjoint union of the $2^{\psi -
\psi(1)}$ copies of the triple $(\varphi, \Xi_{M}, A'_{\psi(1)})$
determines the trivial element in the cobordism group
$Imm^{sf;A_{\psi(1)}}(n-k,k)$.

Let us consider the triple $(2^{\sigma})(-\varphi, -\Xi_{M},
A'_{\psi(1)})$. This triple is define as the disjoin union of
$2^{\sigma}= \psi - \psi(1)$ copies of the triple  $(- \varphi,
-\Xi_{M}, A'_{\psi(1)})$ (the orientation of $M^{n-k}$ is changed
and the immersion $\varphi$ and the skew-framing $\Xi_M$ are
changed to the opposite). The collection of the coefficients for
the triple $(2^{\sigma})(-\varphi, -\Xi_{M}, A'_{\psi(1)})$ will
be denoted by $\{lk(2^{\sigma}(-\mu))\}$. Obviously,
$\{lk(2^{\sigma}(-\mu))\} = - 2^{\sigma}\{lk(\mu)\}$.

Let us consider the triple  $(\hat \varphi', \Xi_{\hat M'}, \hat
A'_{\psi(1)})$, defined as the disjoint union of the triple
 $(\hat
M^{n-k}, \Xi_{\hat M}, \hat A_{\psi(1)})$ with the triple
$(2^{\sigma})(-\varphi, -\Xi_{M}, A'_{\psi(1)})$. The new triple
$(\hat \varphi', \Xi_{\hat M'}, \hat A'_{\psi(1)})$ and the triple
 $(\hat
M^{n-k}, \Xi_{\hat M}, \hat A_{\psi(1)})$ represent the common
element in the cobordism group $Imm^{sf;\hat A_{\psi(1)}}(k,n-k)$.
The mapping $\hat \lambda' : \hat M'^k \to X$, constructed by
means of the collection $\hat A'_{\psi(1)}$ of the cohomology
classes is well-defined. The collection of the intersection
coefficients, defined for the mapping  $\hat \lambda' $ will be
denoted by $\hat lk'(\mu)$. Obviously for an arbitrary multi-index
$\mu$ we have $\hat lk'(\mu)=0$.

A normal surgery of the triple $(\hat \varphi', \Xi_{\hat M'},
\hat A'_{\psi(1)})$  to a triple $(\hat \varphi'', \Xi_{\hat M''},
\hat A''_{\psi(1)})$ by 1-handles is defined such that the the map
$\lambda''$ defined by means of the collection  $\hat
A''_{\psi(1)})$ of cohomology classes satisfy the condition
$(\ref{106})$:
$$\hat \lambda''^{-1}(X^{(n-k)}) =  \emptyset.$$
This gives the first step of the proof.

Let us describe the next steps of the proof. Let us denote the
triple $(\hat \varphi'', \Xi_{M''}, A''_{\psi(1)})$ by $(\varphi,
\Xi_M, A_{\psi(1)})$ again. The map $\lambda: \hat M^{n-k} \to X$,
constructed by means of the collection $A_{\psi(1)}$ of cohomology
classes satisfy the condition $(\ref{106})$. Let us consider the
triple $(\hat \varphi, \Xi_{\hat M}, \hat A_{\psi(1)})$, given by
the following element
$$r^!_{\kappa_{\psi(2)+1}, \dots, \kappa_{\psi(1)}}(\varphi,
\Xi_M, A_{\psi(1)})$$ in the cobordism group $Imm^{sf; \hat
A_{\psi(2)}}(n-k,k)$.

Let us define the space $X(\psi(2))$ by the Cartesian product of
infinite--dimensional projective spaces with indexes $(0,1, \dots,
\psi(2))$:
$$X(\psi(2)) = \RP^{\infty}_0 \times
\prod_{j=1}^{j=\psi(2)} \RP^{\infty}(j).$$ The standard coordinate
inclusion
\begin{eqnarray}\label{incl2}
 i_{\psi(2)}: X(\psi(2))
\subset X(\psi(1))
\end{eqnarray}
is well-defined. The space $X(\psi(2))$ is equipped with the
following standard stratification:
\begin{eqnarray}\label{15strat}
\dots \subset X^{i}(\psi(2)) \subset X^{i-1}(\psi(2)) \subset
\dots \subset X(\psi(2)).
\end{eqnarray}
The  inclusion $(\ref{incl2})$ is agree with the stratifications
$(\ref{105})$, $(\ref{15strat})$ of the origin and the target
spaces.

The collection $A_{\psi(2)}$ of cohomology classes determines the
map  $\hat \lambda: \hat M^{n-k} \to X(\psi(2))$. The condition
$(\ref{106})$ implies the following analogical condition:
\begin{eqnarray}\label{106bis}
\hat \lambda^{-1}(X^{(n-k)}(\psi(2))) = \emptyset.
\end{eqnarray}
Let us denote by $\hat L^1 \subset \hat M^{n-k}$ a 1-dimensional
submanifold in $\hat M^{n-k}$ given by the following formula:
 $$\hat L^1= \hat \lambda^{-1}(X^{(n-k-1)}(\psi(2))).$$
The restriction of the cohomology classes of the collection $\hat
A_{\psi(2)}$ on the submanifold $\hat L^1$ is trivial. In
particular, the submanifold  $\hat L^1$ is framed.

The components of the manifold $\hat L^1$ are equipped with the
collection of the multi-indexes corresponded to the top cells of
the space   $X^{n-k-1}(\psi(2))$. A fixed multi-index determines a
disjoint collection of  $2^{s_1}$ copies of 1-dimensional framed
manifold (probably, non-connected) and the copies are pairwise
diffeomorphic as a framed manifolds.

A framed 2-dimensional manifold $\tilde K^2$ with a framed
boundary $\partial (\tilde K^2)= (\tilde L^1)$ is well-defined.
This framed manifold determines the body of a handle for the
normal surgery of the triple $(\hat \varphi, \Xi_{\hat M}, \hat
A_{\psi(2)})$ to a triple $(\hat \varphi', \Xi_{\hat M'}, \hat
A'_{\psi(2)})$ such that the collection $\hat A'_{\psi(2)}$ of
cohomology classes determines the map
 $$\hat \lambda': \hat
M'^{n-k} \to X(\psi(2)),$$ satisfied the condition
$(\ref{106bis})$.

The next steps of the proof are analogical to the step 1. The
parameter $i$ denoted the dimension of the obstruction is changed
from $2$ up to $n-k$. In each step the analogical condition to the
conditions $(\ref{106})$, $(\ref{106bis})$ is considered. At the
last step of the proof we have a framed manifold
$(M^{n-k},\Xi_M)$, equipped with a collection $\hat A_{\psi(n-k)}$
of the trivial cohomological classes. The framed manifold
represented by $\psi(n-k)$ disjoint copies of the framed manifold
$(M^{n-k},\Xi_M)$ is a framed boundary (and therefore a
skew-framed boundary). The Proposition $\ref{totaltrans}$ is
proved.
\[  \]

Let us describe an algebraic obstruction for the compression of
the given order.

\begin{lemma}\label{codcompr}
An arbitrary element $x \in Imm^{sf;A_j}(n-k,k)$ admits a
compression of the order $i$, $i \le n-k$, if and only if the
element
$$J_{sf}^{k'}(x) \in  Imm^{sf;A_j}(n-k',k')$$
($i \le n-k'
\le n-k$) admits a compression of the same order $i$.
\end{lemma}

\begin{theorem}\label{totalobstr}
For an arbitrary element $x \in Imm^{sf;A_j}(n-k,k)$ the total
obstruction for a compression of an order $q$  ($0 \le q \le n-k$)
is given by the element $J_{sf;\A_j}^{q}(x) \in
Imm^{sf;A_j}(q,n-q)$.
\end{theorem}

To prove  Lemma $\ref{codcompr}$ and Theorem $\ref{totalobstr}$,
let us formulate an auxiliary lemma. Let as assume that a triple
$(\varphi, \Xi_M, A_j)$ represents an element
 $x \in Imm^{sf;A_j}(n-k,k)$. Let us additionally assume that
this element admits a compression of the order $(i-1)$. This means
that in the triple $(\varphi, \Xi_M, A_j)$ can be taken in its
cobordism class such that the
 characteristic class
$\kappa_M \in H^1(M^{n-k};\Z/2)$ of the skew-framing $\Xi_M$ is
given by the following composition:
$$\kappa_M : M^{n-k} \to \RP^{n-k-i} \subset \RP^{\infty},$$
$i < n-k$. We shall denote the map $M^{n-k} \to \RP^{n-k-i}$
described above again by $\kappa_M$.

 Let us consider the manifold
 $Q^i \subset M^{n-k}$, given by the formula:
\begin{eqnarray}\label{Q}
Q^i = \kappa_M^{-1}(pt), \quad pt \in \RP^{n-k-i}.
\end{eqnarray}
The manifold is equipped with the natural framing $\Psi_Q$,
because the restriction of the skew-framing $\Xi_M$ over the
submanifold $Q^i \subset M^{n-k}$ is a framing.

Moreover, the restriction of cohomology classes of the collection
$A_j$ over the submanifold $Q^i \subset M^{n-k}$ determines the
collection $A_j(Q)$ of cohomology classes on $Q^i$. Note that the
class $\kappa_0$ in the collection $A_Q$ is the trivial class. The
immersion $\varphi_Q: Q^i \looparrowright \R^n$ is defined as the
restriction of the immersion $\varphi$ over the submanifold $Q^i
\subset M^{n-k}$. A triple $(\varphi_Q, \Psi_Q, A_j(Q))$
determines an element $J_{sf;\A_j}^i(x)=y \in
Imm^{sf;A_j}(i,n-i)$.

\begin{lemma}\label{compspes}

The element  $x=[(\varphi, \Xi_M, A_M)] \in Imm^{sf,A_j} (n-k,k)$, which
admits a compression of the order $(i-1)$, admits  a
compression of the order $i$ if and only if the element
$J_{sf;\A_j}^i(x)=y =[(\varphi_Q, \Xi_Q, A_j(Q))] \in
Imm^{sf,A_j}(i,n-i)$ is trivial.
\end{lemma}

\subsubsection*{Proof of Lemma $\ref{compspes}$}

At the first step let us prove that if $y=0$ then $x$ admits a
compression of the order $i$. Let us consider a skew-framed in the
codimension $(n-i)$ $(i+1)$--dimensional manifold
$(P^{i+1},\Xi_P)$ with boundary $\partial P^{i+1} = Q^i$, equipped
with the collection $A_j(P)$ of cohomology classes, such that the
restriction of the skew-framing $\Xi_P$ over the boundary $Q^i$ is
a framing coincided with the framing $\Psi_Q$ and the restriction
of the collection $A_j(P)$ over the boundary $Q^i$ coincides with
the collection $A_j(Q)$.

Let us describe a normal surgery of the skew-framed manifold
$(M^{n-k},\Xi_M)$ into a skew-framed manifold $(T^{n-k},\Xi_T)$.
Let us construct a manifold with boundary called the body of a
handle. Let us consider the manifold $P^{i+1}$ and let us denote
the $(n-i)$--dimensional normal bundle over $P^{i+1}$  by $\nu_P$.
The normal bundle $\nu_P$ is equipped with the skew-framing
$\Xi_P$, i.e. the bundle map (an isomorphism on each fiber)
\begin{eqnarray}\label{nuP}
\nu_P \to (n-i)\kappa_P
\end{eqnarray}
is well-defined. Let us denote by $U_P$ the disk bundle over $P$
spanned by the first $(k-i)$  factor in the Whitney sum
$(\ref{nuP})$. The manifold $U_P$ is a skew-framed manifold in the
codimension $(n-k)$ manifold with boundary will be called the body
of a handle.

The boundary $\partial U_P$ of the body of the handle contains a
submanifold $Q^i \times D^{n-k-i}$,  the total space of the disk
bundle over the manifold $Q^i$. Let us consider the Cartesian
product $M^{n-k} \times I$ of the manifold $M^{n-k}$ and the unite
segment $I=[0,1]$.
 A second copy of the manifold
 $Q^i \times D^{n-k-i}$ is embedded into the submanifold  $M^{n-k} \times \{1\} \subset \partial(M^{n-k}
\times I)$, this is a regular neighborhood of the submanifold $Q^i
\times \{1\} \subset M^{n-k} \times \{1\}$. Let us define the
manifold  $T^{n-k}$ by the following formula:
\begin{eqnarray}\label{T}
 T^{n-k}= \partial^+((M^{n-k} \times I) \cup_{Q^i \times D^{n-k-i}}  U_P),
\end{eqnarray}
where by $\partial^+$ is denoted the "upper" component of the
cobordism i.e. the component that contains the last part $\partial
U_P \setminus (Q^i \times D^{n-k-i})$  of the boundary of the body
$U_P$.

After the standard operation called "smoothing the corners" the
$PL$--manifold  $T^{n-k}$ becomes a smooth closed smooth manifold.
The immersion $\varphi_T: T^{n-k} \looparrowright \R^n$ (this
immersion is well-defined up to a regular homotopy), the skew
framing  $\Xi_T$  with the characteristic class $\kappa_T$ (i.e.
the bundle fiberwise isomorphism $\nu_T \to k\kappa_T$) are
well-defined. The manifold $T^{n-k}$ is equipped with the
collection $A_j(T)$ of the collection of characteristic classes,
each class in the collection is determined by the gluing of the
corresponded classes of the two components in the decomposition
$(\ref{T})$. The class $\kappa_0(T)$ of the collection $A_j(T)$
coincides with the characteristic class $\kappa_T$ of the
skew-framing $\Xi_T$. The triple $(\varphi_T, \Xi_T, A_j(T))$
determines an element in the cobordism group $
Imm^{sf;A_j}(n-k,k)$ and by the construction $[(\varphi, \Xi_M,
A_j)] = [(\varphi_T, \Xi_T, A_j(T))]$.

Let us prove that the element $[(\varphi_T, \Xi_T, A_j(T))]$
admits a compression of the order   $i$. We will prove that the
characteristic class $\kappa_T$ is represented by a classifying
map $\kappa_T: T^{n-k} \to \RP^{n-k-i-1} \subset \RP^{\infty}$.
Take a positive integer $b$ big enough and let us consider
$\kappa_M : M^{n-k} \to \RP^b$, such that
$$\kappa_M^{-1}(\RP^{b-n+k+i}) = Q^i ,$$
$Q^i \subset M^{n-k}$, $\RP^{b-n+k+i} \subset \RP^b$.

Let us consider the mapping $g: P^{i+1} \to \RP^{b-n+k+i}$, the
restriction of this mapping to the component of the  boundary
$\partial P^{i+1} = Q^i$ coincides with the map $\kappa_M
\vert_{Q^i}$.  Let us consider the "thickening"  $h: U_P \to
\RP^b$ of the map $g$, this map $h$ is defined by the standard
extension of the map $g$ to the body of the handle $U_P$.

The map  $g': M^{n-k} \cup_{Q^i \times D^{n-k-i}} U_P \to \RP^b$,
$g' \vert_{U_P} = g$ is well-defined and the restriction
 $g' \vert_{T^{n-k} \subset
M^{n-k} \cup U_P}$ does not meet the submanifold
 $\RP^{b-n-k+i} \subset \RP^b$.
 The space
$\RP^{b} \setminus \RP^{b-n-k+i}$ is retracted to the its subspace
$\RP^{n-k-i-1}$ by a deformation, the required compression of the
map $\kappa_T$ of the order  $i$ is constructed. We have proved
that the element $x$ admits a compression of the order $i$.

Let us prove the inverse statement: assume that the element
$x=[(\varphi, \Xi_M, A_j)]$  admits a compression of the order
$i$, then the triple $(\varphi_Q,\Psi_Q, A_j(Q)))$, $Q^i$ is given
by the equation $(\ref{Q})$ determines the trivial element in the
cobordism group $Imm^{sf;A_j}(i,n-i)$.

Let  $(\varphi_W,\Xi_W,A_j(W))$ be a triple, where $W^{n-k+1}$ is
a manifold with boundary, $\partial W^{n-k+1} = M^{n-k} \cup
M_1^{n-k}$; $(\varphi_W, \Xi_W)$ is a skew-framed immersion of the
manifold $W$ into $\R^n \times I$; $A_j(W)$ is a collection of
characteristic classes. Moreover, the triple
$(\varphi_W,\Xi_W,A_j(W))$ determines a cobordism between the
triples $(\varphi_M, \Xi_M, A_j(M))$ and $(\varphi_{M'}, \Xi_{M'},
A_j(M'))$, where the pair $(M'^{n-k}, \kappa_{M'})$ (the
cohomology class $\kappa_{M'}$ is the characteristic class of the
skew-framing $\Xi_{M'}$ and this class is included into the
collection $A_j(M')$) admits a compression of the order $i$, i.e.
the classifying map $\kappa_{M'} = \kappa_W \vert_{M'}$ is given
by the following composition:
$$ \kappa_{M'}: M'^{n-k} \to \RP^{n-k-i} \subset \RP^{\infty}.$$

Let us consider the standard submanifold  $\RP^{b-n-k+i} \subset
\RP^b$, this submanifold intersects the submanifold  $\RP^{n-k-i}
\subset \RP^b$ at a point $pt \in \RP^{n-k-i}\setminus
\RP^{n-k-i-1}$  and does not intersect the standard  submanifold
$\RP^{n-k-i-1} \subset \RP^{n-k-i}$. The image
$Im(\kappa_M(M^{n-k}))$ is in the submanifold $\RP^{n-k-i} \subset
\RP^b$, the image  $Im(\kappa_{M'}(M'^{n-k}))$ is in the
submanifold $\RP^{n-k-i-1} \subset \RP^{n-k-i}$.

Let us denote by  $P^{i+1}$ the submanifold
$F^{-1}(\RP^{b-n+k+i})$ (we assume that $F$ is transversal along
the submanifold $\RP^{b-n+k+i} \subset \RP^{b}$). By the
construction $\partial P^{i+1} = Q^i$. Let us define a
skew-framing $\Xi_P$ in the codimension $(n-i)$ as the direct sum
of a skew-framing of the submanifold $P^{i+1} \subset W^{n-k+1}$
and the skew-framing $\Xi_W$, restricted to the submanifold
$P^{i+1} \subset W^{n-k+1}$.

The restriction of the skew-framing $\Xi_W$ on $\partial W^{n-k+1}
=Q^{n-k}$ coincides with the skew-framing $\Psi_Q$ with the
trivial characteristic class $\kappa_Q$ (i.e. the skew-framing
$\Psi_Q$ is the framing). The restriction of the collection
$A_j(P)$ of cohomology classes on $\partial W^{n-k+1} =Q^{n-k}$
coincides with the collection $A_j(Q)$. This proves that the
triple $(\varphi_Q,\Psi_Q, A_Q)$, is a boundary. Lemma
$\ref{compspes}$ is proved.

\subsubsection*{Proof of Theorem $\ref{totalobstr}$}

Let us assume that a compression of the order $(i-1)$, $i<q$ for
an element $x \in Imm^{sf;A_j}(n-k,k)$, $x=[(\varphi, \Xi_M,
A_M)]$ is well-defined. By Lemma $\ref{compspes}$, the obstruction
to a compression of the order $i$ of the element $x \in
Imm^{sf;A_j}(n-k,k)$ represented by the same triple $(\varphi,
\Xi_M, A_M)$, coincides with the obstruction of a compression of
the same order $i$ for the element $J_{k'}^{sf}(x) \in
Imm^{sf;A_j}(n-k',k')$.  Therefore, by induction over $i$, the
total obstruction for a compression of the order $q$ for the
element $x$ is trivial if and only if the total obstruction for a
compression of the order $q$ for the element $J_{k'}^{sf}(x)$ is
trivial. Theorem $\ref{totalobstr}$ is proved.
\[  \]

To prove the Compression Theorem $\ref{comp}$ the following
construction by U.Koshorke of the total obstruction for a homotopy
of a bundle map into a bundle monomorphism on each fiber of the
bundles (see [K]) is required.

Let
 $\alpha \to Q^q$, $\beta \to Q^q$ be a pair of the vector bundles over the smooth manifold
 $Q^q$ (we do not assume that the manifold $Q^q$ is closed)
  $dim(\alpha)=a$,
$dim(\beta)=b$, $dim(Q^q)=q$,  $2(b-a+1)<q$. Let $u: \alpha \to \beta$ be a generic vector bundle morphism.
let us denote by
$\Sigma \subset Q^q$ a submanifold, given by the formula:
\begin{eqnarray}\label{Sigma}
\Sigma = \{ x \in Q^q \vert Ker(u_x: \alpha_x \to \beta_x) \ne 0
\}.
\end{eqnarray}
This manifold $\Sigma$ is the singular manifold of the bundle morphism $u$.
Note that under the presented dimensional restrictions, for a generic vector bundle morphism $u$
we have
$rk(u) \ge a-1$.The codimension of the submanifold
$\Sigma \subset Q^q$ is equal to
$b-a+1$.

Let us describe the normal bundle  of the submanifold
$(\ref{Sigma})$, this bundle will be denoted by $\nu_{\Sigma}$. Let us denote by $\lambda: E(\lambda) \to \Sigma$ the linear subbundle,
determined as the subbundle of kernels of the morphism $u$ over the singular submanifold
 $\Sigma \subset Q^q$. Therefore, the inclusion of the bundles over $\Sigma$
  $\varepsilon: \lambda \subset \alpha$ is well-defined.
Let us denote by $\Lambda_{\alpha}$ the bundle over $\Sigma$, this
bundle is the orthogonal complement to the subbundle
$\varepsilon(\lambda) \subset \alpha$. A natural vector-bundle
morphism over $\Sigma$ (isomorphism of fibers)
 $v: \Lambda_{\alpha} \subset \beta$ is well-defined.
Let us define the bundle $\Lambda_{\beta}$ over $\Sigma$ as the
orthogonal complement to the subbundle $v(\Lambda_{\alpha})$ in
the bundle $\alpha \vert_{\Sigma}$. The normal bundle
$\nu(\Sigma)$ is determined by the following formula:
\begin{eqnarray}\label{016}
\nu(\Sigma) = \lambda \otimes \Lambda_{\beta}.
\end{eqnarray}

If the manifold
$Q^q$ has a boundary $\partial Q$ and the vector bundles morphism  $u$ is the
morphism of the bundles over the manifold with boundary, then the singular submanifold
$\partial \Sigma \subset \partial Q$ of the restriction $u \vert_{\partial Q}$
is a boundary of the submanifold $\Sigma \subset Q^q$ with the normal bundle, given by the same formula $(\ref{016})$.

In the paper [K] (in this paper there is a reference to the
previous papers by the same author) a cobordism group of
embeddings of of manifolds in $Q^q$ (in this construction the
manifold $Q^q$ is closed) of codimension
 $b-a+1$ with an additional structure of the normal bundle, given by the equation  ($\ref{016}$)
is defined. For an arbitrary generic vector bundle morphism $u:
\alpha \to \beta$ an element in this cobordism group is
well-defined. This element is the total obstruction of a homotopy
of the vector bundle morphism $u$ to a fiberwise monomorphism.

Let $\hat \nu: E(\hat \nu) \to \RP^{2^k-1}$ be a vector bundle, $\dim(\hat \nu)=
n+1-2^k$, $2^k < n+2$ over the standard projective space, isomorphic to the following Whitney sum:
$\hat \nu \equiv (n+1-2^k)\kappa_{\RP}$, where  $\kappa_{\RP}$ is the canonical line bundle over
$\RP^{2^k-1}$. Let us denote the Whitney sum $\hat \nu \oplus \kappa$ by
$\nu$, $\dim(\nu)= (n-2^k+2)$. The standard projection $\pi: \nu \to
\hat \nu$ with the kernel $\kappa$ is well-defined.
The bundle $\nu_{\RP}$  in the case
\begin{eqnarray}\label{b}
b(2^k) \le n+2
\end{eqnarray}
(the positive integer $b(r)$, $r=2^k$ is equal to the corresponded
power of $2$, see  [A-E]) is isomorphic to the normal bundle
of the projective space $\RP^{2^k-1}$.

Let us define an admissible $s$--family of sections (singularities
in the family of sections are possible) of the bundle $\hat
\nu_{\RP}$.

\subsubsection*{Definition of an admissible family
of sections of the bundle $\hat \nu_{\RP}$}

We shall say that a generic $s$--family of sections
$$\hat \psi =
\{\hat \psi_1, \dots, \hat \psi_s\}, \quad \hat \psi: s\varepsilon \to \hat \nu_{\RP}$$
of the bundle $\hat \nu_{\RP}$ is admissible, if
there exists a regular $s$--family
$$\psi = \{
\psi_1, \dots, \psi_s\}, \quad \psi: s\varepsilon \to \nu_{\RP}$$
of sections of the bundle  $\nu_{\RP}$ satisfied the condition:
$\pi \circ \psi=\hat \psi$.
\[  \]

\begin{lemma}\label{addmiss}
Let us assume that $s \le n+2-2^{k+1}-k-1$ and $n \equiv -2
\pmod{2^{2^k}}$, $n >0$. Then the bundle $\hat \nu_{\RP}$ has an
admissible $s$-family of sections.
\end{lemma}

\subsubsection*{Proof of Lemma $\ref{addmiss}$}

By the Davis table  the projective space $\RP^{2^k-1}$ is
immersable into the Euclidean space $\R^{2^{k+1}-k-1}$ (this is
not the lowest possible dimension of the target Euclidean space of
immersions). By the equation $(\ref{b})$, and because $b(k) \le
2^k$ for $n=2^{2^k}$, the bundle $\nu_{\RP}$ is the normal bundle
over the projective space $\RP^{2^k-1}$. Therefore the bundle
$\nu_{\RP}$ admits a generic regular $s$--family of section,
denoted by $\psi$. The projection $\hat \psi = \pi \circ \psi$ of
this regular family is the admissible $s$--family of sections of
the bundle $\hat \nu_{\RP}$. Lemma $\ref{addmiss}$ is proved.

Let us consider an admissible generic $s$--family of sections $\hat \psi$ of the bundle $\hat \nu_{\RP}$.
Let us denote by $\Sigma \subset \RP^{2^k-1}$ the singular manifold of the family $\hat \psi$.
This denotation corresponds to $(\ref{Sigma})$, if we take $\alpha \equiv s \varepsilon$,
$\beta \equiv \hat \nu$, $u=\psi$. In the following lemma we will describe the normal bundle
$\nu_{\Sigma}$ of the submanifold $\Sigma \subset \RP^{2^k-1}$.

\begin{lemma}\label{7.5}
Let us assume that
\begin{eqnarray}\label{s}
s = n+2-2^{k+1}-k-1, \quad k \ge 2,
\end{eqnarray}
where $n = 2^{2^k}-2$. Then the singular submanifold $\Sigma
\subset \RP^{2^k-1}$ of an admissible $s$--family of sections of
the vector-bundle $\hat \nu_{\RP}$ is a smooth submanifold of
dimension $k$, the normal bundle $\nu_{\Sigma}$ is equipped with a
skew-framing $\Xi_{\Sigma}: (2^k-1-k) \kappa_{\RP} \equiv
\nu_{\Sigma}$, the characteristic class  $\kappa_{\Xi}$ of this
skew-framing coincides with the restriction
$\kappa_{\RP}\vert_{\Sigma \subset \RP^{2^k-1}}$.
 \end{lemma}

\subsubsection*{Proof of Lemma $\ref{7.5}$}

 Let us describe the normal bundle $\nu_{\Sigma}$ by means of Koschorke's Theorem.
Moreover, let us define a skew-framing of this bundle.
Let us denote by
$\lambda \subset s\varepsilon $ the subbundle of the kernels of the family $\hat \psi$
over the submanifold $\Sigma$. By the assumption the the family $\hat \psi$ is admissible,
therefore:
\begin{eqnarray}\label{lambda}
\lambda = \kappa_{\Sigma},
\end{eqnarray}
where by $\kappa_{\Sigma}$ the vector bundle $\kappa_{\RP} \vert_{\Sigma}$ is denoted.
Let us denote the orthogonal complement to $\kappa_{\Sigma} \subset s \varepsilon$ over $\Sigma$ by  $s\varepsilon -\kappa_{\Sigma}$.
Let us denote the orthogonal complement to $\hat \psi(s\varepsilon -\kappa_{\Sigma})$ in the bundle
$\hat \nu_{\RP}\vert_{\Sigma}$ by $\Lambda$.
By the construction:
\begin{eqnarray}\label{barPsi}
(s\varepsilon - \kappa_{\Sigma}) \oplus \kappa_{\Sigma} = s \varepsilon.
\end{eqnarray}
Let us prove that the bundle  $\nu_{\Sigma}$ satisfies the equation:
\begin{eqnarray}\label{lambda1}
 \nu_{\Sigma} \equiv (n+2 - 2^k -k-1) \kappa_{\Sigma}.
 \end{eqnarray}

 By the Koschorke Theorem the following isomorphism of vector bundles is well-defined:
 $$\Lambda \otimes \kappa_{\Sigma} \equiv \nu_{\Sigma}.$$
 This equation is equivalent to the equation:
  $$(s\varepsilon)\otimes \kappa_{\Sigma} \equiv \nu_{\Sigma} \oplus \kappa_{\Sigma}.$$
 This proves the equation $(\ref{lambda1})$.
   The isomorphism $(\ref{lambda1})$ defines a skew-framing $\Xi_{\Sigma}$
of the bundle $\nu_{\Sigma}$   with the characteristic class $\kappa_{\Sigma}$.
The Lemma $\ref{7.5}$ is proved.

 \subsubsection*{Remark}
Because the restriction of the normal bundle  $\nu_{\RP}$ over
$\Sigma$ is isomorph to the trivial bundle by the canonical
isomorphism, the skew-framing $\Xi_{\Sigma}$ determines the
skew-framing of the normal bundle (in the Euclidean space) of the
manifold $\Sigma$.
\[  \]

Let us consider an element
$$ x \in Imm^{sf,A_j}(2^{k}-2,n-2^{k}+2), \quad n=-2(mod 2^{k}),$$
given by the cobordism class of a triple $(\varphi, \Xi_M, \A_j(M))$. Put $m=2^{k}-2$.
Let us consider the vector bundle
$\nu \to \RP^{2^k-1}$, $\dim(\nu)=n-m$.
The normal bundle
$\nu_M$ of the manifold $M^{m}$ is given by the formula:
$$\nu_M = \kappa_M^{\ast}(\nu_{\RP}).$$
Let us consider the map
\begin{eqnarray}\label{lam}
\lambda_M: M^{m} \to \RP^{2^k-1} \times \prod_{i=1}^j \RP^{\infty},
\end{eqnarray}
see $(\ref{lambdaM})$, constructed by means of the collection
$A_j$. Let us consider the standard projection  $\pi_0:
\RP^{2^k-1} \times \prod_{i=1}^j \RP^{\infty} \to \RP^{2^k-1}$ on
the factor $\RP^{2^k-1}$. The composition  $\pi_0 \circ \lambda_M:
M^m \to \RP^{2^k-1}$ coincides with the map $\kappa_M$.

Let us define the subbundle $\hat \nu_M \subset \nu_M$ of the
codimension 1 (i.e. of the dimension $\dim(\hat \nu_M)=n+1-2^{k}$)
by the formula:
$$\hat \nu_M = \kappa_M^{\ast}(\hat \nu_{\RP}).$$
Let us define a family of sections  $\hat \xi_M=\{\hat \xi_1,
\dots, \hat \xi_s\}$, $s=n+3-2^{k+1}+k$ of the bundle $\hat
\nu_M$. This collection is the pull-back image of an admissible
collection of sections $\hat \psi = \{\hat \psi_1, \dots \hat
\psi_s\}$ of the bundle $\hat \nu_M$ by the map $\kappa_M$. Let us
denote by
\begin{eqnarray}\label{NsubM}
N^{k-1} \subset M^m
\end{eqnarray}
the submanifold  of singular sections. It is not hard to prove
that $\dim(N^{k-1})=(k-1)$. The manifold $N^{k-1}$ is equipped
with the collection $A_j(N)$ of cohomology classes, a class of
$A_j(N)$ is defined as the restriction of the corresponded class
of the collection $A_M$ over the submanifold  $N^{k-1} \subset
M^m$. The class $\kappa_M \vert_N$ in the collection $A_j(N)$ is
denoted by $\kappa_N$. Let us denote the immersion $\varphi
\vert_N$ by $\varphi_N$.

Let us denote by $\nu_N$ the normal bundle of the immersed
(embedded by the general position arguments) manifold
$\varphi_N(N^{k-1})$ in the Euclidean space $\R^n$. The normal
bundle $\nu_N$ is isomorph to the Whitney sum  $\nu_N = \nu_M
\vert_N \oplus \nu_{N \subset M}$, where by $\nu_{N \subset M}$ is
denoted the normal bundle of the submanifold $N^{k-1} \subset M^m$
inside the manifold $M^m$.

By Lemma $\ref{7.5}$ and by the transversality of the map
$\kappa_M$ along the submanifold $\Sigma \subset \RP^{2^k-1}$, the
bundle $\nu_{N \subset M}$ is equipped with the skew-framing
$\Xi_{N \subset M}$ with the characteristic class $\kappa_N$. The
bundle $\nu_M \vert_N$ is also equipped with a skew-framing with
the same characteristic class (see an analogical Remark after
Lemma $\ref{7.5}$). This gives a skew-framing  $\Xi_N$ of the
immersion $\varphi_N$ of codimension $(n-k+1)$. Let us denote by
$A_j(N)$ the collection of cohomology classes from the group
$H^1(N^{k-1};\Z/2)$, this collection is the restriction of the
collection  $A_j$ over the submanifold $N^{k-1} \subset M^m$.

\begin{lemma}\label{7.6}
The triple $(\varphi_N, \Xi_{N}, A_j(N))$ determines an element
$x_{k-1} \in Imm^{sf,A_j}(k-1,n-k+1)$, this element is the total
obstruction of a compression of the order $(k-1)$ for the element
$x \in Imm^{sf,A_j}(2^{k}-2,n-2^{k}+2)$.
\end{lemma}

\subsubsection*{Proof of Lemma $\ref{7.6}$}

Let us consider the submanifold $\Sigma^k \subset \RP^{2^{k}-1}$
of singularities of the family of sections  $\hat \psi$, let us
re-denote this manifold by $\Sigma_0^{k}$. This manifold is
equipped by the following natural stratification (a filtration):
\begin{eqnarray}\label{017}
\emptyset \subset \Sigma_{k}^0 \subset \dots \subset
\Sigma_{1}^{k-1} \subset \Sigma^{k}_0 \subset \RP^{2^{k}-1}.
\end{eqnarray}

The submanifold $\Sigma_i$, $dim(\Sigma_i)=k-i$ in $(\ref{017})$
is defined as the singular submanifold of the subfamily of the
first $(s-i)$ sections in $\hat \psi$. By the straightforward
calculations follows that the fundamental class
 $[\Sigma_i]$ of the corresponded submanifold in the filtration in the group
$H_{k-i}(\RP^{2^{k}-1};\Z/2)$ represents the only generator of
this group: this homology class is dual to the characteristic
class $\bar w_{n+1-2^{k}-k+i}(n+1-2^{k})\kappa_{\RP}$ of the
bundle $\hat \nu_{\RP}$.

Without loss of the generality we assume that the map
$\kappa_M : M^m \to \RP^{2^{k}-1}$ is transversal along
the stratification  $(\ref{017})$.
Let us denote the inverse image of the stratification
 $(\ref{017})$ by
\begin{eqnarray}\label{018}
 N_{k-1}^0 \subset N_{k-2}^1 \subset \dots \subset N_{0}^{k-1}
\subset M^m.
\end{eqnarray}
The top manifold $N_0^{k-1}$ of the filtration  $(\ref{018})$  coincides with manifold
$N^{k-1}$, defined above.

Let us prove the lemma by the induction over the parameter $i$,
$i=0, \dots, k-1$. Let us assume that the image of the map
$\kappa_M : M^m \to \RP^{2^{k}-1}$ is in the standard projective
subspace $\RP^{2^{k}-2-i} \subset \RP^{2^{k}-1}$.  In this case
$N^{i-1}_{k-i}=\emptyset$. By the standard argument we may assume
that the stratum $\Sigma_{k-i}^i$  intersects in the general
position the standard submanifold $\RP^{2^{k}-1-i}$ of the
complementary dimension at the only point. (The index of
self-intersection of this two submanifolds in the manifold
$\RP^{2^{k}-1}$ is odd and well-defined modulo 2.)

The framed manifold  $N_{k-i-1}^i$ is the regular preimage of the
marked point by the map $\kappa_M$ (the image of this map is in
the submanifold
 $\RP^{2^{k}-2-i}$). Let us denote an element represented by the cobordism class of the triple $(\varphi_{N_{k-i-1}},\Xi_{N_{k-i-1}},
A_j({N_{k-i-1}}))$ in the cobordism group $Imm^{sf;A_j}(i,n-i)$ by $x_i$. By Lemma $\ref{7.5}$ the condition
 $x_i=0$, $i \le k-2$ is satisfied if and only if there exists a normal cobordism
of the map $\kappa_M$ to the map $\kappa_{M'}: M'^{m} \to
\RP^{2^{k}-2-i} \subset \RP^{2^{k}}$.
%From the equation $x_{k-1}=0$ we have the equation $x_{i+1}=0$.
Therefore a
compression of the order $i+1$, $i+1 \le k-1$ of an element $x$ is
well defined. If we put $i+1 = k-1$ we have a compression of the
order $(k-1)$. Lemma $\ref{7.6}$ is proved.
 \[  \]

\[  \]

The following Proposition is the main in the proof of Theorem $\ref{comp}$.

\begin{proposition}\label{38}
Let $n=-2(mod(2^k))$, $n>2^k$, $x \in Imm^{sf;A_j}(2^k-2,n-2^k+2)$ be an arbitrary element in the kernel of the homomorphism 
$(\ref{Omega})$ (in this formula $b(k)=2^k$). Let $x_{k-1} \in Imm^{sf;A_j}(k-1,n-k+1)$
be the total obstruction for a compression of the order $(k-1)$ of the element
$x$. Then the element
$x_{k-1}$ is in the image of the transfer homomorphism, i.e. there exists an element $y_{k-1} \in
Imm^{sf,A_{j+1}}(k-1,n-k+1)$, for which
$r_{j+1}(y_{k-1})=x_{k-1}$. Assuming $k-1$ is even, the element $y_{k-1}$ is in the kernel of the homomorphism $(\ref{Omega})$ (where we assume that $b(k)=k+1$)
\end{proposition}

\subsubsection*{Proof of Proposition $\ref{38}$}

Let us assume that the cobordism class of the element $x$ is given
by a triple $(\varphi_M,\Xi_M,A_j(M))$, where $\varphi_M: M^{m}
\looparrowright \R^n$ is an immersion, $\dim(M^{m})=m=2^k-2$. 
Because $codim(\varphi_M)$ is odd, $M^m$ is oriented.  Assume that $[(\varphi_M,\Xi_M,A_j(M))]$ is in the kernel of $(\ref{Omega})$ (in this formula $b(k)=2^k$).
(This assumption is satisfied for $j=0$, because characteristic numbers  for $M^m$ are trivial.)

Let
us consider the normal bundle $\nu_M \cong (n-2^{k}+2)\kappa_M$
over the manifold $M^{m}$ and the subbundle  $\hat \nu_M \subset
\nu_M$ of the codimension 1, $\hat \nu_M \equiv
(n-2^{k}+1)\kappa_M$. Let us prove that there exists a regular
$s$--family of sections $\hat \psi$ of the bundle $\hat \nu_M$,
$s=n+2+k - 2^{k+1}$.

Let us denote by $D(\kappa_M)$ a manifold with boundary, the total
space of the disk bundle, associated with line bundle $\kappa_M$.
The vector bundle $\hat \nu_M$ is lifted to the vector bundle over
$D(\kappa_M)$ (we will denote this lift again by  $\hat \nu_M$).
By the R.Cohen theorem [C] there exists an immersion  $D(\kappa_M)
\looparrowright \R^{2^{k+1}-2-k}$, because
$\dim(D(\kappa_M))=2^k-1$ and $\alpha(2^k-1)=k$.
The proof of the statement is
in Theorem $\ref{39}$.

Equivalently, the
bundle $\hat \nu_M$ admits an $s$--family of regular sections. The
tautological lift, denoted by $\psi$, of the regular $s$--family
$\hat \xi$ of the bundle $\hat \nu_M$ to a regular $s$--family of
sections of the bundle $\hat \nu_M$ is defined.

Let us consider the admissible $s$--family of sections $\hat \psi$
of the bundle $\hat \nu_M$. This family of sections is defined as
the pull-back of admissible $s$--family of sections of the bundle
$\hat \nu_{\RP}$, see Lemma $\ref{addmiss}$. A regular $s$--family
of sections $\psi$ of the bundle $\nu_M$ is defined as the lift of
the admissible $s$--family $\hat \psi$.

Let us consider the manifold $M^m \times I$ and let us define the
bundle $\nu_{M \times I}$ by the formula $\nu_{M \times I} = p_M^{\ast}(\nu_M)$, where $p_M: M^m \times I \to M$
is the standard projection on the second factor.
Let us define the
bundle $\hat \nu_{M \times I}$ by the formula $\hat \nu_{M \times I} = p_M^{\ast}(\hat \nu_M)$.
Let us consider a generic $s$--family of sections
$X= \{\chi_1, \dots \chi_s\}$ of the bundle  $\nu_{M \times I}$
with the following boundaries conditions:
\begin{eqnarray}\label{220}
X = \psi \quad over \quad M^{m} \times \{1\},
\end{eqnarray}
\begin{eqnarray}\label{221}
X = \xi \quad over \quad M^{m} \times \{0\}.
\end{eqnarray}

Let us denote by $V^{k-1} \subset M^{m} \times I$ the singular subset of the family $X$.
By the general position argument this subset is a closed submanifold in $M^{m} \times I$,
because over the boundary $\partial(M^{m} \times I)$ the family $X$ is regular.
Let us denote by
$$\hat X \{\hat \chi_1, \dots, \hat \chi_s \}$$
the projection of the $s$--family $X$ into the $s$--family of sections of the vector bundle
$\hat \nu_{M \times I}$. The $s$--family $\hat X$ satisfies the following boundary conditions:
\begin{eqnarray}\label{222}
\hat X = \hat \psi \quad over \quad M^{m} \times \{1\},
\end{eqnarray}
\begin{eqnarray}\label{223}
\hat X = \hat \xi \quad over \quad M^{m} \times \{0\}.
\end{eqnarray}

Let us denote by $\hat K^{k} \subset M^{m} \times I$ the subset of
singular sections of the $s$--family $\hat X$. This subset $\hat
K^k$ is a $k$--dimensional manifold with boundary. The only
component of the boundary $\partial \hat K^k$ is a submanifold of
$M^{m} \times \{1\})$, this component coincides with the
submanifold $N^{k-1} \subset M^m \times \{1\}$ of the
singularities of the admissible family $\psi$, see
$(\ref{NsubM})$.

By Lemma $\ref{7.5}$, the triple $(\varphi_N, \Xi_N, A_N)$ is
well-defined. Here $\varphi_N = \varphi \vert_N$, $\Xi_N$ is the
skew-framing of this immersion, $A_j(N)$ is the restriction of the
collection  $A_j$ on $N^{k-1} \subset M^m \times \{1\}$.
 The triple $(\varphi_N, \Xi_N, A_N)$ represents the total obstruction $x_{k-1} \in Imm^{sf;A_j}(k-1,n-k+1)$
for a compression of the order  $(k-1)$ of the element
$x=[(\varphi, \Xi_M, \kappa_M)]$.

Let us use the formula
$(\ref{016})$ to calculate the normal bundle of the submanifold
$\hat K^{k} \subset M^m \times I$ and of the line normal bundle of the submanifold
$V^{k-1} \subset \hat K^{k}$.

Let us denote by $\lambda$ the line bundle over $V^{k-1}$ of kernels of the $s$--family of sections $X$.
Let us prove that the normal bundle $\nu_{V \subset M \times I}$ of the submanifold
$V^{k-1}
\subset M^{m} \times I$ is given by the formula:
\begin{eqnarray}\label{normV}
\nu_{V \subset M \times I} \equiv \varepsilon \oplus  (m-k+1)\lambda, \quad m=2^k-2.
\end{eqnarray}
The restriction of the normal bundle $\nu_{M \times I} \vert_{V}$
over the submanifold $V^{k-1}$ is isomorph to the bundle
$(n+2-2^{k})\kappa_{M \times I}$, and therefore, because $b(k-1)
\le 2^k$, is isomorph to the trivial bundle:
\begin{eqnarray}\label{norm}
(n+2-2^{k})\kappa_{M \times I} \equiv (n+2-2^{k})\varepsilon.
\end{eqnarray}
This isomorphism is canonical, i.e. does not depends of $M^m$ and
of $V^{k-1}$. Let us define the vector bundle $\Lambda$ as the
orthogonal complement to the line subbundle $\lambda$ in the
trivial bundle of linear combinations of the base sections. The
bundle $\Lambda$ represents the stable vector bundle $-\bar
\lambda$. Therefore the orthogonal complement of the subbundle
$\Xi(\Lambda)$ in the vector bundle $\nu_{M \times I}$ is the
subbundle in the vector bundle $\nu_M \vert_V$ isomorph to the
vector bundle $\lambda \oplus (2^{k} - k -1)\varepsilon$. By the
Koschorke Theorem, the normal bundle of the submanifold $V^{k-1}
\subset M^{2^{k}-2} \times I$ is isomorph to the vector bundle:
$$\nu_{V \subset M \times I} \equiv \lambda \otimes(\lambda \oplus (2^{k}-k-1)\varepsilon)
\equiv \varepsilon \oplus (2^{k}-k-1)\lambda.$$ This vector bundle
$\nu_{V \subset M \times I}$ also represents the stable normal
bundle of the manifold $V^{k-1}$ because of the equation
$(\ref{norm})$. The formula $(\ref{normV})$ is proved.

The restriction of the immersion $\varphi \times Id\vert_V$ is
regular homotopic to an immersion
 $\varphi_V: V^{k-1} \looparrowright \R^n \times \{1\}$. By the computation the normal bundle of the immersion
 $\varphi_V$ is equipped with a skew-framing, denoted by $\Xi_V$. The collection of cohomology classes $A_j(V)$
 is defined by the formula $A_j(V) = A_j(M \times I) \vert_V$, where $A_j(M \times I)$ is induced from the given collection $A_j(M)$ of cohomology classes on $M^m$ by the projection $p_M$.
Let us define the collection of cohomology classes
$A_{j+1}(V)$ by the addition to the collection $A_j(V)$ the last cohomology class  $\kappa_{j+1} = \lambda \otimes \kappa_{M \times I}$.
The triple
 $y_{k-1}=(\varphi_V, \Xi_V, A_{j+1}(V))$ determines an element in the cobordism group $Imm^{sf;A_{j+1}}(k-1,n-k+1)$.

Let us denote by $U_V \subset \hat K^{k}$ a small closed regular
neighborhood of the submanifold  $V^{k-1} \subset \hat K^{k}$. The
line bundle of kernels of the family $\hat X$ of sections over the
submanifold (with boundary) $U_V \subset M^m \times I$ is isomorph
to the line bundle $p_{U_V,V}^{\ast}(\lambda)$, where $p_{U_V,V}:
U_V \to V$ is the projection of the neighborhood on the central
submanifold. Let us denote the line bundle
$p_{U_V,V}^{\ast}(\lambda)$  by $\hat \lambda$. The orthogonal
complement to the subbundle $\hat \lambda$ in the vector bundle of
the linear combinations of the base sections over $U_V$ will be
denote again by $\hat \Lambda$.

Let us consider the subbundle $\hat X(\hat \Lambda)$ in the vector
bundle $\hat \nu_{M \times I}$. By the analogical calculations,
using the canonical isomorphisms of the vector bundles over $U_V$:
$(2^{k-1})\lambda \equiv (2^{k-1})\kappa_M \equiv
(2^{k-1})\varepsilon$, the orthogonal complement in $\hat \nu_{M
\times I}$ of the subbundle $\hat X(\hat \Lambda)$, denoted by
$-\hat X(\hat \Lambda)$, is isomorph to the bundle $\lambda \oplus
(2^{k-1}-k-1)\varepsilon \oplus (2^{k-1}-1)\kappa_M$.  By the
Koschorke Theorem, the stable isomorphism class of the normal
bundle $\nu_{U_V \subset M \times I}$ of the submanifold $U_V
\subset M^{m} \times I$ is given by the formula:

$$\nu_{U_V \subset M \times I} \equiv (-k-1)\hat \lambda \oplus (-\kappa_M \otimes
\hat \lambda).$$

In particular, from this calculation follows that the line normal bundle of the submanifold
$V^{k-1} \subset K^{k}$ is isomorph to the line bundle  $\lambda \otimes \kappa_M$.

Let us denote $\partial U_V$ by $Q^{k-1}$. The space $Q^{k-1}$ is a closed manifold, $\dim(Q^{k-1})=k-1$.
The normal bundle $\nu_{Q \subset M\times I}$ of the submanifold $Q^{k-1} \subset M^m \times I$ is given
by the formula:
$$ \nu_{Q \subset M\times I} \equiv \varepsilon \oplus (m-k+1) \kappa_M. $$
The restriction of the immersion $\varphi \times id: M^m \times I \looparrowright \R^n \times I$
to the submanifold $Q^{k-1} \subset M^m \times I$ is regular homotopic to a skew-framed immersion
$\hat \varphi_Q: Q^{k-1} \looparrowright \R^n \times \{1\}$ оf the codimension $(n-k+1)$
with the sklew-framing, denoted by $\hat \Xi_Q$, and with the characteristic class of this skew-framing $\hat \kappa_Q = \kappa_{M \times I} \vert_Q$. The manifold $Q^{k-1}$ is equipped by the collection $\hat A_j(Q)$ of cohomology classes, $\hat A_j(Q) = A_j(M \times I) \vert_{Q \subset M \times I}$. The triple
 $(\varphi_Q, \Xi_Q, A_j(Q))$ determines an element in the cobordism group  $Imm^{sf;A_j}(k-1,n-k+1)$.

The manifold $K^{k} \setminus U_V$ has the boundary consists of
the two components:
 $\partial(K^{k} \setminus U_V)= Q^{k-1} \cup K^{k-1}$.
 The restriction of the immersion $\varphi \times Id \vert_{(K^{k} \setminus U_V)}$ is regular homotopic to an immersion
 $\varphi_K: K^k \looparrowright \R^n \times I$ with the following the boundary conditions:
 $\varphi_K \vert Q = \varphi_Q: Q^{k-1} \looparrowright \R^n \times \{1\}$,
 $\varphi_K \vert N = \varphi_N: N^{n-1} \looparrowright \R^n \times \{0\}$.
 The immersion $\varphi_K$ is a skew-framed immersion with a skew-framing $\Xi_K$ and with the characteristic class $\kappa_K = \kappa_{M \times I} \vert_K$ of this skew-framing. The manifold $K^{k}$ is equipped by the collection of cohomology classes $A_j(K) = A_j(M \times I) \vert_{K \subset M \times I}$.

The triple
 $x_{k-1}=(\varphi_N, \Xi_N, A_j(N))$ determines an element in the cobordism group  $Imm^{sf;A_j}(k-1,n-k+1)$.
 The element $x_{k-1}$ is the total obstruction to a compression of the element $x=[(\varphi,\Xi_M,A_j)]$ of the order $(k-1)$.
The element $[(\varphi_Q, \hat \Xi_Q, \hat A_j(Q))]$ is the image by the transfer homomorphism $r_{j+1}^!$ of the element
 $y_{k-1}=[(\varphi_V, \Xi_V, A_{j+1}(V))]$. The elements   $x_{k-1}=[(\varphi_N, \Xi_N, A_j(N))]$ and $r^!_{j+1}(y_{k-1})=[(\varphi_Q, \hat \Xi_Q, \hat A_j(Q))]$ are equal. The cobordism between the elements $x_{k-1}$ and $r^!_{j+1}(y_{k+1})$ is given by the triple $(\varphi_K, \Xi_K, \A_j(K))$.

Assume that $k-1$ is even. Let us prove that $y_{k-1}$ is in the kernel of $(\ref{Omega})$ (in this formula we assume that $b(k)=k+1$). The construction of  $y_{k-1}$ 
from $x$ is generalized for the image of $x$ by $(\ref{Omega})$ (in this formula $b(k)=2^k$). Because $x$ is in the kernel of $(\ref{Omega})$,
$y_{k-1}$ is also in the kernel of $(\ref{Omega})$.
Proposition $\ref{38}$ is proved.

\subsubsection*{Proof of the Compression Theorem $\ref{comp}$}

Let us define a positive integer $\psi = \psi(q)$, by the formula
 $(\ref{psi})$ for $n-k=q-1$. By Proposition
$\ref{totaltrans}$ the total transfer homomorphism $(\ref{104})$, which is defined
on the group $Imm^{sf;A_{\psi}}(d-1,n-d+1)$ is the trivial
homomorphism. Let us define a positive integer  $l(d) =
\exp_2(\exp_2 \dots \exp_2(d)\dots +1)$, where the number of the
iterations of the function $\exp_2(x+1) = 2^{x+1}$ is equal to
$\psi$ and the initial value is $x=d-1$.

Let  $l'$ be an arbitrary power of $2$, $l' \ge l(d)$. Let us
define  $n=l'-2$. Let us prove that an arbitrary element in
$Imm^{sf}(n-1,1)$ admits a compression of the order $d-1$.

Let us define $n_0=l(d)-2$, by the assumption $n_0 \le n$. Let us
define the following sequence of $\psi$ integers: $2n_1 =
log_2(n_0+2)-2$, $2n_2= log_2(n_1+2)-2$, $\dots,$
 $2n_{\psi}=log_2(n_{\psi-1} +2)-2$. All this integers are positive and $n_{\psi}=d-1$.

Let $x_0 \in Imm^{sf}(n_0,n-n_0)$ be the image an arbitrary element
$x \in Imm^{sf}(n-1,1)$ by $J^{sf}: Imm^{sf}(n-1,1) \to Imm^{sf}(n_0,n-n_0)$.
Denote by $x_{d-1}$  the total obstruction of a compression of
the order $d-1$ for the element $x$. 
The element $x_{d-1}$ coincides with the total obstruction of a
compression of the order  $(d-1)$ for the element $x_{n_0}$.

Let us consider the total obstruction $x_{n_1} \in
Imm^{sf}(n_1,n-n_1)$ of the retraction of the order $n_1$ for the
element $x_{n_0}$. The element $x_{n_1}$ is in the kernel of the homomorphism 
$(\ref{Omega})$ (in this formula we put $b(k)-2=n_1$ is even), because characteristic classes of $M^{n_1}$ are trivial).
By Proposition $\ref{38}$, there exists an
element $y_{n_1} \in Imm^{sf;\{\kappa_1\}}(n_1,n-n_1)$, such that
the image of the element $y_{n_1}$ by the transfer homomorphism is
equal to $x_{n_1}$.

Let us consider the total obstruction $y_{n_2} \in
Imm^{sf;\kappa_1}(n_2,n-n_2)$ of a compression of the order $n_2$
for an element $y_{n_1}$. By the Proposition $\ref{38}$ there
exists an element  $z_{n_2} \in Imm^{sf; \kappa_1,
\kappa_2}(n_2,n-n_2)$ such that the element $r_2(z_{n_2})$ is the
total obstruction of a compression of the order  $n_2$ for the
element $y_{n_1}$. The element $r_{tot}(z_{n_2})=r_1 \circ
r_2(z_{n_2})=J^{sf}_{n_2}(x_{n_0})=x_{n_2}$ is the total
obstruction for a compression of the order  $n_2$ of the element
$x_{n_1}$. The same element $x_{n_2}$ is the total obstruction for
a compression of the order  $(n_2)$ for the initial element $x$.

By the induction we prove that the total obstruction of a
compression of the order  $(d-1)$ for the element $x$ is in the
image of the total transfer homomorphism of the multiplicity
 $\psi$, i.e. this total obstruction is equal to
 $r_{tot}(z)$,  $z \in Imm^{sf;A_{\psi}}(d-1,n-d+1)$. By Proposition
$\ref{totaltrans}$ we have  $r_{tot}(z)=0$. Therefore,
$x_{d-1}=0$. The Compression Theorem $\ref{comp}$ is proved.
\[  \]

%\subsection*{Appendix}

In the proof of Proposition $\ref{38}$ we used the following fact.
Denote $b=b(k)=2^{2k}$.
Assume that $n=2^l-2$, $l \ge b(k)$. %Let us denote the integer
% $(n-2b(k)+k+3)$ by $s$.

\begin{theorem}\label{39}
Let  a cobordism class $x \in
Imm^{sf,A_r}(b(k)-2,n-b(k)+2)$ is represented by a skew-framed immersion $(f: M^{b(k)-2} \looparrowright \R^n,\Xi_M)$,
$\dim(M^{b(k)-2})=b(k)-2$, where the manifold $M^{b(k)-2}$ is
equipped with a mapping $(\ref{lam})$.
Assume that the element $x$ belongs to the kernel of the forgetful homomorphism
\begin{eqnarray}\label{Omega}
Imm^{sf,A_r}(b(k)-2,n-b(k)+2) \to \Omega_{b(k)-2}(\prod_{i=0}^r \RP^{\infty}_i).
\end{eqnarray}
%into orientable bordism group is trivial.
%Hurewicz image
%\begin{eqnarray}\label{Hur}
%(\lambda_M)_{\ast}([M]) \in H_{b(k)-2}(\prod_{i=0}^r \RP^{\infty}_i;\Z)
%\end{eqnarray}
% with integer coefficients is trivial.
Then in the regular cobordism class $[x]$ there exists an element for which the manifold $M^{b(k)}$ admits an immersion into the Euclidean
space $\R^{2b(k)-3-k}$ with a non-degenerate skew cross section given by a linear bundle $\kappa_M$.
\end{theorem}

%\begin{theorem}\label{40} An arbitrary element $\alpha \in
%Imm^{sf}(b(k)-1,n-b(k)+1)$ (the case $r=0$), is represented by a
%skew-framed manifold $(M^{b(k)-1},\Xi_M)$ with a family $\{\xi_1,
%\dots, \xi_{s}\}$ of $s$--sections of the normal bundle
%$\nu(M)=(n-b(k)+1)\kappa_M$, that is a regular section over a
%submanifold $N^{b(k)-2} \subset M^{b(k)-1}$, dual to the
%cohomological class $\kappa_M$.
%\end{theorem}

%We shall prove this theorem using the
%main result [A-E] and the following result: the standard real
%projective space $\RP^{b(k)-1}$ is immersed into the Euclidean
%space $\R^{2b(k)-k-2}$.
%\[  \]
%The condition that the homomorphism $(\ref{Omega})$ is trivial  in the assumption in Theorem $\ref{39}$ gives no restriction for a proof of %Proposition $\ref{38}$. 
%Assume that a manifold, equipped with the collection $A_j$, represented the element $x$ in this proposition is an oriented boundary.
%Then a manifold, which is equipped with the collection $A_{j+1}$, and which represents the element $y_{k-1}$, is an oriented boundary.
%Evidently, the manifold  $M^{2^{l'}-l'}$, equipped with cohomology class $\kappa_M \in H^1(M^{2^{l'}-l'};\Z/2)$, which represents an element in  %$Imm^{sf}(2^{l'}-l',l'-2)$
%is an oriented boundary.

\subsection*{Remark}
 Theorem $\ref{39}$
is a corollary of the R.Cohen's Immersion Theorem [C]. 
\[  \]
 
The main step of the proof of Theorem $\ref{39}$ is the following lemma.

\begin{lemma}\label{Cohen}
Let $(N^{j}, \Xi_N)$, $j=0, \dots, k-1$, be a framed (in particular, oriented) manifold, equipped with a mapping
$$\lambda_N: N^j \to \prod_{i=1}^r \RP^{\infty}_i.$$
Let us assume that the Hurewicz image
\begin{eqnarray}\label{Hur}
(\lambda_N)_{\ast}([N]) \in H_{j}(\prod_{i=1}^r \RP^{\infty}_i;\Z)
\end{eqnarray}
 with integer coefficients % system, which is associated with the trivial mapping into the first factor $\RP^{\infty}_{0}$
% of the product, 
 is trivial. Then there exists a skew-boundary $(W^{k},\Psi_W, \lambda_W)$ in codimension $b(k)-1$,
 $\partial(W^{k}_W) = N^{k-1}$, $\Psi_W \vert_{\partial W^{k}} = \Xi_N$, $\lambda_W \vert {\partial W=N^k} = \lambda_N$, where
$$\lambda_W: W^{k} \to \RP^{\infty}_0 \times \prod_{i=1}^r \RP^{\infty}_i, $$
 and $\kappa_W = w_1(\Psi_W)$ coincides with the projection of the mapping $\lambda_W$ on the factor $\RP^{\infty}_0$.
\end{lemma}
\[  \]

\subsubsection*{A sketch of the proof of Lemma $\ref{Cohen}$}
 To prove Lemma  $\ref{Cohen}$ let us consider the Atiyah–Hirzebruch spectral sequence  for the
cobordism group
of framed immersions of the dimensions $0, \dots, k-1$ (the stable homotopy group) of the space $\prod_{i=1}^r \RP^{\infty}_i$.
 Let us consider the Atiyah–Hirzebruch spectral sequence  for
cobordism groups  
of skew-framed immersions of dimensions $0, \dots, k-1$ in the codimension $b(k)-1$ of the space $\RP^{\infty}_{0} \times \prod_{i=1}^r \RP^{\infty}_i$.
There is a natural mapping of the first spectral sequence to the second spectral sequence.
By the main result of [A-E] all higher coefficients  in $E_2$-terms of the kernel  are trivial. Therefore the kernel of $E^{\infty}$--therm is totally described by
the Hurewicz image $(\ref{Hur})$ of a cobordism class of a corresponding mapping of a framed manifold. This Hurewicz image is trivial by the assumption. Therefore an arbitrary framed manifold $(N^{k}, \Xi_N)$ is a skew-framed boundary in codimension $b(k)-1$.
 Lemma Lemma $\ref{Cohen}$ is proved.

\subsubsection*{A sketch of the proof of Theorem $\ref{39}$}

By the assumption the normal bundle
 $\nu_M$ is isomorphic to the Whitney sum  $Cb(k)\kappa_M$, where $\kappa_M$ is the given line bundle over  $M^{b(k)-2}$,
$b(k)=2^{2k}$, $C$ in a positive  integer. 
Let us calculate the Koschorke construction to prove that the normal bundle $Cb(k)\kappa_M$ of $M^{b(k)-2}$ 
admits a regular skew-section  $\kappa_M \oplus k\varepsilon \oplus (C-1)b(k)\varepsilon$.
%Это свойство сохраняется в процессе построения. 
%Докажем, что расслоение $\nu_M$ допускает послойный мономорфизм расслоения $k\varepsilon \oplus \kappa_M$. 
The construction is given by the induction over the index $j=0, \dots, k$.

The bundle
 $\nu_M$ admits a regular family of  $2k+(C-1)b(k)$ sections. 
 Denote by 
 $I \subset \nu_M$ the orthogonal complement to the subbundle  $(C-1)b(k)\varepsilon  \subset \nu_M$,
 which consists of the first
$(C-1)b(k)$ sections of the given family. At a $j$-th step of the induction 
$j=0, \dots, k$, let us assume that there exists a morphism (this morphism could be tot a fiberwise isomorphism) 
 $\rho_{j}$ of the bundle $(2k-j+1)\varepsilon \oplus \kappa_M$
into the subbundle $I \subset \nu_M$, such that the following conditions are satisfied:

--1. Denote by  $\hat \rho_{j}$ the restriction of the morphism  $\rho_{j}$ on the subbundle 
$(j+1) \varepsilon \oplus \kappa_M \subset (2k-j+1)\varepsilon \oplus \kappa_M$. It is required that  $\hat \rho_{j}$ 
is regular (a fiberwise isomorphism). 

--2. Denote by  $\tilde \rho_j$ the restriction of the morphism $\rho_{j}$ at the subbundle  $(2k-j+1)\varepsilon \subset (2k-j+1)\varepsilon \oplus \kappa_M$. It is required that 
$\tilde \rho_{j}$ is regular.

Let us prove the base of the induction $j=0$. Consider the manifold $M^{b(k)-2}$ as a manifold with the prescribed orientation, equipped with the collection $A_r$ of cohomology classes. This manifold determines an element in the cobordism group
$\Omega_{b(k)-2}(\prod_{i=0}^r \RP^{\infty}_i)$, which is given by the image of the element  $x$ by the homomorphism $\Omega$.
The obstruction of the existence of a regular morphism $\hat \rho_{0}$ for such a manifold is well-defined and is trivial, because  $x$ 
belongs to the kernel of $(\ref{Omega})$. An extension of the morphism
$\hat \rho_0$ to a morphism  $\rho_0$, which satisfies the condition  --2, is well defined.

Let us prove the step
$j \mapsto j+1$ of the induction. Consider the morphism $\rho_{j}$ and denote the restriction of this morphism 
on the subbundle
$(j+2) \varepsilon \oplus \kappa_M \subset (2k-j+1)\varepsilon \oplus \kappa_M$ by $\hat \rho^{\circ}_j$ 
(one more section then in the family
$\hat \rho^{\circ}_j$). The morphism $\hat \rho^{\circ}_j$, generally speaking, is not regular. Denote by  $N^{j+1} \subset M^{b(k)-2}$
the singular submanifold. The conditions  --1 and --2 imply that the restriction of the cohomology class $\kappa_M$ on
$N^{j+1}$ is trivial and that  $N_{j+1}$ is a framed submanifold in  $M^{b(k)-2}$, and, therefore,  a framed manifold. 
Let us prove that this framed manifold satisfies conditions of Lemma 
 $\ref{Cohen}$.

Let us prove that the  Hurewicz image $(\ref{Hur})$ is trivial.
There are the two cases: 

--a. $j$ is even; 

--b. $j$ is odd.

Let us consider the case a. Take the collection of section $\rho_{j+2}$, which is restricted to the subbundle
$(j+3) \varepsilon \oplus \kappa_M \subset (2k-j+1)\varepsilon \oplus \kappa_M$, denote this restriction by  $\hat \rho^{\circ \circ}_{j}$
(two extra sections with respect to $\hat \rho^{\circ}_j$). Denote the singular manifold of the morphism
 $\hat \rho^{\circ \circ}_{j}$ by  $K^{j+2}$. Evidently, 
$w_1(K^{j+2})=\kappa_M \vert_{K}$. The local fundamental class of the manifold manifold $K^{j+2}$, equipped by the restriction of the collection $A_r$, determines an element $y$
in  $H_{j+2}(\RP^{\infty}_0 \times\prod_{i=1}^r \RP^{\infty}_i;\Z^{tw})$, where the coefficients are integers and twisted with respect to the cohomology class $\kappa_M$.  Because the element 
$x$ belongs to the kernel of  $(\ref{Omega})$, $y=0$.  Therefore the homology class  $[K^{j+2}] \cap \kappa_M = [L^{j+1}]$, 
which is considered as an element in the group  $H_{j+1}(\RP^{\infty}_0 \times\prod_{i=1}^r \RP^{\infty}_i;\Z)$ is trivial. 
This proves that the statement in the case --a.

Let us consider the case b. Take the collection of section
$\rho_{j}^{\circ}$, which contains a regular subfamily 
$\rho_{j}$. Because the restriction  $\kappa_M$ on the manifold  $N^{j+1}$ is trivial, the tensor product of the morphism  $\rho_{j}^{\circ}$ 
and the line bundle  $\kappa_M$ is well defined as the morphism with the same singular manifold. In particular,
$\rho_{j}^{\circ} \otimes \kappa_M$ determines a morphism of the bundle  $(C-1)b(k)\kappa_M \oplus (j+1)\kappa_M \oplus \varepsilon$ 
into the bundle  $I \otimes \kappa_M$, which is the orthogonal complement of the subbundle   $(C-1)b(k)\kappa_M \subset Cb(k)\varepsilon$.
Because $j+1$ is even, and because $x$ is in the kernel of $(\ref{Omega})$, $N^{j+1}$ determines the element is the kernel of $(\ref{Hur})$. 
This proves that the statement in the case --b.

The obstruction of the existence of a regular morphism $\hat \rho_{j+1}$
with the trivial kernel over the singular manifold is given by a  cobordism class of a mapping $\lambda_N: N^j \to \prod_{i=1}^r \RP^{\infty}_i$ of a framed $j$-dimensional manifold, which is considered as a  skew-framed manifold in codimension $b(k)-1$ with local twisted coefficients system,
associated with $\kappa_M$  (note that $\kappa_M \vert_{N^{j+1}}$ is trivial, but in the regular cobordism class this property is not assumed).

Therefore, there exists 
 $\rho_{j+1}$, which satisfies the condition --1. Let us prove that there exists $\rho_{j+1}$, which satisfies the both conditions --1 and --2. 
 By the construction a homotopy of the morphism
$\hat \rho_{j}^{\circ}$ into the morphism  $\hat \rho_{j+1}$ has singularities with the trivial kernels. Evidently, there exists a homotopy 
$\tilde \rho_j$ into $\tilde \rho'_j$ of the bundle  $(2k-j+1)\varepsilon$, which has singularities with the trivial kernels, and
the restriction of $\tilde \rho'_j$ on the subbundle  $(j+1)\varepsilon$ coincides with the restriction  $\hat \rho_{j+1}$ 
on the considered subbundle. Therefore there exists a morphism 
$\rho'_{j+1}$ of the bundle $(2k-j+1)\varepsilon \oplus \kappa_M$ (one more sections then in
$\rho_{j+1}$), for which the condition --1 is satisfies, and instead of the condition  --2 the following condition is satisfied:
the restriction of $\rho'_{j+1}$ on the subbundle   $(2k-j+1)\varepsilon \subset (2k-j+1)\varepsilon \oplus \kappa_M$ has the trivial 
kernel. By general position arguments the restriction
 $\rho'_{j+1}$ on the trivial line bundle is regular. Let us restrict 
 $\rho'_{j+1}$ on the subbundle  $(2k-j)\varepsilon \subset (2k-j+1)\varepsilon \oplus \kappa_M$. There exists a small generic deformation 
 $\rho'_{j+1}\vert_{(2k-j)\varepsilon} \mapsto \tilde \rho_{j+1}$, for which the family $\tilde \rho_{j+1}$ is regular. 
This proves the condition 2. The induction is well-defined.

By the last step of the induction, there exists a regular section $\rho_{k-1}$.
Theorem $\ref{39}$ is proved.

\[  \]
\[  \]
Moscow Region, Troitsk, 142190, IZMIRAN, 
pmakhmet@mi.ras.ru

\end{document}